\documentclass[10pt,oneside,reqno]{amsart}
\usepackage{amsmath, amsfonts, amsthm}
\usepackage{algorithm,algorithmic}
\usepackage{cite}
\usepackage{booktabs}
\usepackage{geometry}
\usepackage{graphicx}
\usepackage{algorithm}
\usepackage{algorithmic}
\usepackage{epstopdf}
\usepackage{color}
\usepackage{subfigure}
\usepackage{pifont}
\usepackage{amsthm}
\usepackage{enumitem}
\usepackage{etoolbox}
\usepackage{amssymb}
\usepackage{amsmath,amsthm,mathtools}
\usepackage{caption}

\newtheorem{remark}{Remark}

\newcommand{\vect}[1]{\boldsymbol{#1}}
\newcommand{\mat}[1]{\mathbf{#1}}

\begin{document}

\title[A Stable SPH Discretization of the Elliptic Operator]{A Stable SPH Discretization 
of the Elliptic Operator with Heterogeneous Coefficients}

\author{Alexander~A.~Lukyanov$^{\dag}$}
\author{Kees~Vuik$^{\dag}$}
\keywords{Meshless method, Stability and Approximation, Laplace operator,  
Diffusive flow, Monotone scheme, Discrete maximum principle}

\subjclass[2010]{97N40, 65N75, 76M28}
\thanks{$^{\dag}$Delft University of Technology, Faculty of Electrical Engineering, 
Mathematics and Computer Science, Delft Institute of Applied Mathematics, 2628CN 
Delft, the Netherlands (aalukyanov1@gmail.com, c.vuik@tudelft.nl).}

\begin{abstract}
		Smoothed particle hydrodynamics (SPH) has been extensively used 
		to model high and low Reynolds number flows, free surface flows and 
		collapse of dams, study pore-scale flow and dispersion, elasticity, and 
		thermal problems. In different applications, it is required to have a stable 
		and accurate discretization of the elliptic operator with homogeneous and 
		heterogeneous coefficients. In this paper, the stability and approximation 
		analysis of different SPH discretization schemes (traditional and new) of 
		the diagonal elliptic operator for homogeneous and heterogeneous media 
		are presented. The optimum and new discretization scheme is also proposed. 
		This scheme enhances the Laplace approximation (Brookshaw's scheme 
		\cite{Brookshaw1985} and Schwaiger's scheme \cite{Schwaiger2008}) 
		used in the SPH community for thermal, viscous, and pressure projection 
		problems with an isotropic elliptic operator. The numerical results are 
		illustrated by numerical examples, where the comparison between different 
		versions of the meshless discretization methods are presented.
\end{abstract}


\date{\today}
\maketitle

\section{Introduction}
\label{Intro}
Smoothed particle hydrodynamics (SPH) was developed a few decades ago
to model inviscid fluid and gas flow dynamics in astrophysical problems
\cite{Lucy1977, GingoldMonaghan1977,GingoldMonaghan1982, Monaghan1992}.
The SPH is an interpolation-based numerical technique that can be used to
solve systems of partial differential equations (PDEs) using either Lagrangian
or Eulerian descriptions. The nature of SPH method allows to incorporate
different physical and chemical effects into the discretized governing equations
with relatively small code-development effort. In addition, geometrically complex
and/or dynamic boundaries, and interfaces can be handled without undue difficulty.
The SPH numerical procedure of calculating state variables (i.e., density, velocity,
and gradient of deformation) are computed as a weighted average of values in a local
region. Despite a few advantages of SPH method, this method is not free
from disadvantages. For example, for fluids, gases, or solids with non-trivial 
boundaries there is incompleteness of the kernel support 
combined with the lack of consistency of the kernel interpolation in conventional 
meshless methods which results in fuzzy boundaries. In some cases, this can be fixed 
by an automatic incorporation of the boundary condition \cite{Lucy1977}. 
The completeness of mesh free particle methods was discussed in
\cite{BelytschkoKrongauzDolbowGerlach1998}.
However, care must be taken to ensure that variables whose values 
do not approach zero at boundaries are accurately represented.

It has been observed in the literature that meshless methods (e.g., SPH) are not
free from instabilities, especially in the modeling solid mechanics 
problems. For example, the tensile instabilities were identified in 
\cite{SwegleHicksAttaway1995} by performing the Neumann 
analysis of the one-dimensional governing equations (conservation laws). 
Therefore, different stabilization techniques have been developed. It is 
important to note at this point that it is very difficult to perform the general 
Neumann analysis in two- and three-dimensions. Furthermore,
the high-frequency instability results from the low order discretization (rank
deficient) of the divergence operator. The tensile instability results from the
interaction between the second derivative of the Eulerian kernel (i.e., computed
kernel in the Eulerian coordinates) and the tensile stress. The tensile
instability only occurs for the Eulerian kernels because the Eulerian kernels
depend on both the stress and the second derivative of the kernel. It has been
shown \cite{BelytschkoGuoLiuXiao2000} that in the case where 
the kernel is a function of the Lagrangian coordinates (the Lagrangian kernel), 
tensile instability does not occur. A comprehensive analysis on this subject can 
be found in \cite{BelytschkoGuoLiuXiao2000}. Recently, some regularization and
stabilization of SPH schemes were proposed in \cite{ZhangHuAdams2017}, 
\cite{KhayyerGotohShimizu2017}.  

Since its introduction, SPH has been successfully used to model a wide range of
fluid flows and the behavior of solids subjected to large deformations.
For example, the SPH method was applied to simulate high energy explosions
\cite{LiuLiuLam2003} and impact \cite{Lukyanov2007,LukyanovPenkov2007}, most
notably free surface flows and collapse of dams \cite{Monaghan1994},
elastoplasticity 
\cite{JohnsonStrykBeissel1996,ChenBeraunJih1999b,ChenBeraunJih2001,LukyanovPenkov2007},
to model low Reynolds number flows \cite{MorrisFoxZhu1997,HolmesWilliamsTilke2010},
to study pore-scale flow and dispersion \cite{ZhuFoxMorris1999,ZhuFox2001},
and for thermal problems 
\cite{ChenBeraunJih1999b,ChaniotisPoulikakosKoumoutsakos2002,Schwaiger2008}.

In different applications including, but not limited to,
fluid flow related problems, it is required to have the 
stable and accurate discretization of the following operator 
(diagonal elliptic operator), which is the 
research subject of this paper:
\begin{equation}
  \label{LV_SPH_Lap}
    \begin{array}{lcl}
            \displaystyle
			\vspace{0.2cm}
			\mathbf{L}\left(\mathbf{u}\right) = -
            \nabla\left(\mathbf{M}\left(\vect{r}\right)\nabla \mathbf{u}
			\left(\vect{r}\right)\right) - g\left(\vect{r}\right), \
            \forall \vect{r}\in\Omega\subset\mathbb{R}^{n}, \\
            \mathbf{M}\left(\vect{r}\right) \in \{ \mbox{diag}[m(\vect{r})]: 
	        m(\vect{r})\in \mathbb{R}^{n}_{+}, \ 
	        m\left(\vect{r}\right)\in L_{2}\left(\Omega\right)\}, \ \mbox{diag}: 
	        \mathbb{R}^{n}\rightarrow \mathbb{R}^{n\times n},
	  \end{array}
\end{equation}
where $\mathbf{u}\left(\vect{r}\right)$ is the unknown scalar or vector variable
field, $\mathbf{M}^{\alpha\beta}\left(\vect{r}\right)$ is the 
diagonal matrix of the mobility field, e.g., one example includes 
$\mathbf{M}^{\alpha\beta}\left(\vect{r}\right) = 
m\left(\vect{r}\right)\delta_{\alpha\beta}, \  \alpha,\beta = 1,\ldots,n;$, where 
$m\left(\vect{r}\right)$ is the mobility scalar field, $\delta_{\alpha\beta}$ 
is the Kronecker symbol, $n=1,2,3$ is the spatial 
dimension. The sink/source term $g\left(\vect{r}\right)$ is 
assumed to be zero in this paper. Consider the operator in the expression
(\ref{LV_SPH_Lap}) with piecewise continuous coefficients
$\mathbf{M}\left(\vect{r}\right)$ in $\Omega$.
It has been noticed that some of the numerical methods for elliptic equations 
may violate the maximum principle (i.e. lead to spurious oscillations). Therefore,
proposed methods must satisfy a discrete maximum principle to avoid any spurious
oscillations. This is also applicable to meshless discretizations. Usually, the
oscillations are closely related to the poor approximation of the variable gradient 
$\nabla\mathbf{u}$ in the flux computation. In this paper, different numerical 
discretizations of the elliptic isotropic operator are analyzed. The objective of this 
paper is to develop numerical scheme satisfying the two-point flux approximation 
nature in the form
\begin{equation}
	\label{TPFA_FORM}
	\displaystyle
	\mat{L}\left(\mathbf{u}\right) \approx
	\sum_{J} V_{\vect{r}_{J}}
	\Psi\left(\left[\mathbf{u}\left(\vect{r}_{J}\right)-\mathbf{u}
	\left(\vect{r}_{I}\right)\right],
	\vect{r}_{J}-\vect{r}_{I}\right)-g\left(\vect{r}_{I}\right),
\end{equation}
where $J$ is the neighboring particle of the particle $I$, 
$V_{\vect{r}_{J}}$ is the $J$ - particle volume, $\Psi$ is the 
special kernel. This structure \eqref{TPFA_FORM} allows to have 
a better analysis of the stability and monotonicity. Furthermore, 
this will allow to apply upwinding strategy during the solution of 
nonlinear PDEs. The optimum discretization scheme of the shape 
\eqref{TPFA_FORM} is proposed in this paper. This scheme is based both on 
the Laplace approximation (Brookshaw's scheme \cite{Brookshaw1985}) 
and on a gradient approximation commonly used in the SPH community 
for thermal, viscous, and pressure projection problems. The proposed 
discretization scheme is combined with mixed corrections, which ensure 
linear completeness. The mixed correction utilizes Shepard Functions in 
combination with a correction to derivative approximations. In corrected 
meshless methods, the domain boundaries and field variables at the
boundaries are approximated with the improved accuracy comparing 
to the conventional SPH method. The resulting new 
scheme improves the particle deficiency (kernel support 
incompleteness) problem. 
The outline of the paper is as follows. In section 
\ref{SPH-discretization-Laplace-operator}, the existing 
discretizations of the Laplace operator with the building blocks necessary for 
these methods are discussed. In this section,  the SPH kernel and its gradient 
properties are also discussed. The description of meshless transmisibilities 
and their connections to the existing mesh-dependent discretization schemes 
is given in section \ref{MeshlessTransmissibilities} including the 
construction of a new meshless discretization scheme. The 
approximation, stability, and monotonicity analysis are performed 
in section \ref{ApproximationStabilityMonotonicity}. The numerical 
analysis of different boundary value problems is presented in section 
\ref{SolutionBoundaryValueProblems}. The paper is concluded by 
section \ref{Summary}.

\section{SPH discretization of the Laplace operator}
\label{SPH-discretization-Laplace-operator}
The meshless approximations (SPH approximations) to the operator
\eqref{LV_SPH_Lap} with heterogeneous and homogeneous 
coefficients are presented in this section. 
Let us consider a rectangle in 
$\mathbb{R}^{n}, \ n=1,2,3$:
\begin{equation}\label{LV_SPH_Box}
\displaystyle \Omega = \left
\{\vect{r}=\left\{x_{i}\right\}\in \mathbb{R}^{n} \left|\right. \
\ 0< \left|x_{i}-a_{i}\right| < l_{i}, \ l_{i}\in \mathbb{R}_{+},
\ \ \forall i=1,\ldots,n \right\}
\end{equation}
as the numerical domain. Here, $a_{i}$ are the center coordinates of the
rectangular and $l_{i}$ are the side lengths. The following norms are used 
to quantify the accuracy of different approximations for the entire 
numerical domain $\Omega$:
\begin{equation}\label{LukyanovVuik_SPH_Norm}
	\displaystyle
	\left\|\mathbf{f}\right\|^{{p}}_{\bar{\Omega}}=
	\left(
	\sum\limits_{\vect{\xi}_{k}\in \bar{\Omega}}V_{\vect{\xi}_{k}}
	\left(\left|\mathbf{f}\right|\right)^{p}\right)^{\frac{1}{p}}, \
	\left\|\mathbf{f}\right\|^{{p}}_{\tilde{\Omega}}=
	\left(
	\sum\limits_{\vect{\xi}_{k}\in \tilde{\Omega}}V_{\vect{\xi}_{k}}
	\left(\left|\mathbf{f}\right|\right)^{p}\right)^{\frac{1}{p}}
\end{equation}
where $\mathbf{f}$ is the approximated physical quantity, 
$V_{\vect{\xi}_{k}}$ is the volume of the particle $\vect{\xi}_{k}$,
$\bar{\Omega}$ denotes entire domain including boundary particles, 
$\tilde{\Omega}$ denotes only internal part of the domain ${\Omega}$.

The proposed discretization schemes should be compatible with 
a discontinuous $m\left(\vect{r}\right)$ (or piecewise 
function) coefficient of the operator \eqref{LV_SPH_Lap} 
since this coefficient cannot be differentiable in the 
classical sense.
The standard SPH spatial discretization of the Laplace operator \eqref{LV_SPH_Lap}
arises from the following relations ($g\left(\vect{r}\right)=0, \ \forall \vect{r}$):
\begin{equation}
   \label{LV_SPH_ID}
	 \begin{array}{lcl}
         \displaystyle
	     \vspace{0.2cm}
         \langle\mat{L}\left(\mathbf{u}\left(\vect{r}_{I}\right)\right)
         \rangle =\lim\limits_{\tilde{h}_{I}\rightarrow 0}
	     \int\limits_{\Omega_{\vect{r}_{I},\tilde{h}_{I}}}\mathbf{L}
		 \left(\mathbf{u}
	     \left(\vect{r}\right)\right)
	     W\left(\vect{r}-\vect{r}_{I},\tilde{h}_{I}\right)dV_{\vect{r}} \approx \\
         \displaystyle
	     \approx\lim\limits_{\tilde{h}_{I}\rightarrow 0}\left[
	     \int\limits_{\Omega_{\vect{r}_{I},\tilde{h}_{I}}}
	     \mathbf{M}\left(\vect{r}\right)\nabla  \mathbf{u}\left(\vect{r}\right)\cdot
	     \nabla W\left(\vect{r}-\vect{r}_{I},\tilde{h}_{I}\right)dV_{\vect{r}}
		 \right],
	 \end{array}
\end{equation}
where $W\left(\vect{r}-\vect{r}_{I},\tilde{h}_{I}\right)$ is the kernel that
weakly approximates the Dirac delta function $\delta\left(\vect{r}-\vect{r}_{I}\right)$ but with 
finite characteristic width $\tilde{h}_{I}$ around the particle $I$. The effective characteristic 
width $\tilde{h}_{I}$ will be defined in the upcoming section using real smoothing particle 
length $h_{I}$. However, it is important to note that, in the case of the homogeneous particle 
distribution, it is common to use $\tilde{h}_{I}=f\cdot h_{p}, \ f \geq 1$, where $h_{p}$ is the 
inter-particles distance. Hence, it is important that $\tilde{h}_{I}$ and 
$h_{p}$ are not mixed up during the analysis. The control volumes in the meshless 
discretization are the patches, which are interior to the support of the kernels
$W\left(\vect{r}-\vect{r}_{I},\tilde{h}_{I}\right)$,
i.e. $\Omega_{\vect{r}_{I},\tilde{h}_{I}} =
\mbox{supp} \ {W}\left(\vect{r}-\vect{r}_{I},\tilde{h}_{I}\right)$, 
$\tilde{h}_{I}$ is the diameter (or smoothing
length) of the particle $I$, and $\vect{r}$, $\vect{r}_{I}$ are points in
Euclidean space $\mathbb{R}^{n}$. Additionally, it can be required that the
kernels are radially symmetric and compactly supported ($\mbox{supp} \ W$) as:
\begin{equation}\label{LV_SPH_ID1}
   \displaystyle
   \Omega_{\vect{r}_I,\tilde{h}_{I}}=\mbox{supp} \ W =
   \lbrace \mathbf{r}_{J} \
   \vert \ W(\vect{r}_{J}-\vect{r}_{I},\tilde{h}_{IJ})\neq 0\rbrace,
\end{equation}
where $\tilde{h}_{IJ}$ is the effective smoothing length between particles
located at $\vect{r}_{I}$ and $\vect{r}_{J}$, which will be defined in the
upcoming section. From the definition of the $\Omega_{\vect{r}_{I},\tilde{h}_{I}}$,
it follows that there is an infinite cover of the numerical domain $\Omega$:
\begin{equation}\label{LV_SPH_ID2}
   \displaystyle
   \Omega = \bigcup
	 \limits_{\vect{r}_I\in\Omega}\Omega_{\vect{r}_I,\tilde{h}_{I}}.
\end{equation}
According to the Heine-Borel theorem, there is a finite subcover (since we
consider only compact numerical domain $\Omega$), that is
\begin{equation}\label{LV_SPH_ID3}
   \displaystyle
   \Omega = \mbox{span}
	 \left\{\Omega_{\vect{r}_I,\tilde{h}_{I}}/ I=1,...,N\right\},
\end{equation}
where $N$ is the number of particles in the numerical discretization. The SPH
spatial discretization of the integral \eqref{LV_SPH_ID} is 
defined over the control volumes $\Omega_{\vect{r}_I,\tilde{h}_{I}}$ to obtain 
the final discretization of the Laplace operator. The final step in the particle method 
is to approximate the integral relation on the right-hand side of the (\ref{LV_SPH_ID}) 
using Monte-Carlo expressions or any cubature rules \cite{FulkQuinn1996,ChenBeraunJih1999b,
BelytschkoKrongauzDolbowGerlach1998,LiuLiu2003} and which is known as 
the particle approximation step.

The assumption that the boundary term from the integration by parts is zero in
\eqref{LV_SPH_ID} is valid only in regions where the kernel has a full support, or 
the function, or the gradient of the function itself is zero. For particles near free 
surfaces or boundary, the neglect of these terms leads to significant errors for 
boundary value problems. Several techniques have been developed to address
these errors through various correction methods, e.g., by calculating the boundary
integrals \cite{LiuLiu2003}. In addition, this can be corrected as it will be shown 
below by applying normalized corrected meshless methods in the derivative
approximations. In the following section, the commonly used SPH kernels and its 
gradients are considered.

\subsection{SPH kernel and its gradient}
\label{SPH_Kernel_Grad}

A central point of the SPH formalism is the concept of the interpolating function
(or kernel) through which the continuum properties of the medium are recovered
from a discrete sample of $N$ points \eqref{LV_SPH_ID2} with prescribed mass
$m_{I}$ (for conventional Lagrangian methods) or volume $V_{I}$ (for fully
Eulerian methods). In the Lagrangian formulation, these points move according to
the specified governing laws, whereas these points are fixed in space for the
Eulerian formulation. A good interpolating kernel must satisfy a few basic
requirements: it must weakly tend to the delta function in the continuum limit
and has to be a continuous function with piecewise first 
derivatives at least. From a more practical point of view it is also advisable to 
deal with symmetric finite range kernels, the latter to avoid $N^2$ calculations. 
Cubic and quintic splines are the commonly used kernels in SPH formulations
\cite{LiLiu2002,LiuLiu2003}. Since the quintic spline does not provide the
numerical advantages, the cubic spline is used in this paper:
\begin{equation}\label{LV_SPH_K}
   \displaystyle
	 W(z,\tilde{h}) =\frac{\Xi}{\tilde{h}^{D}}\left\{
           \begin{array}{lcl}
              \vspace*{0.2cm}
		          \displaystyle   1 -\frac{3}{2}z^{2}+\frac{3}{4}z^{3}, &
							0\leqslant z \leqslant 1; \\
              \vspace*{0.2cm}
			  \displaystyle   \frac{1}{4}(2-z)^{3}, & 1 \leqslant z \leqslant 2; \\
			  \displaystyle   0, & z > 2;
		   \end{array}\right.
\end{equation}
where $\displaystyle z=\left\|\vect{r}_{J}-\vect{r}_{I}\right\|_{2}/\tilde{h}$ is
the dimensionless variable, $\tilde{h}=\tilde{h}_{IJ}$ is the effective smoothing
length between particles $I$ and $J$, and $\displaystyle\Xi$ is the normalization
factor equal to $\displaystyle {3}/{2}$, $\displaystyle {10}/{\left(7\pi\right)}$,
and $\displaystyle {1}/{\pi}$ in 1D, 2D, and 3D, respectively. The different choices
of computing effective smoothing length between particles used in this paper
will be discussed in the upcoming subsection.

Although, the kernel is normalized in continuous sense, it is important to note
that $W\left(\mathbf{r}_{J}-\mathbf{r}_{I},\tilde{h}_{IJ}\right)$ does not satisfy
the normalization condition in the discrete space
$$
\displaystyle
\sum\limits_{\mathbf{r}_{J}\in\Omega_{\mathbf{r}_{I},\tilde{h}_{I}}}
W\left(\mathbf{r}_{J}-\mathbf{r}_{I},\tilde{h}_{IJ}
\right)V_{\mathbf{r}_{J}}\neq 1
$$
due to the particle distribution and incomplete kernel support 
near the boundary and, hence, the discretized-normalized kernel function can be considered:
\begin{equation}\label{LV_SPH_DN}
    \begin{array}{lcl}
	    \displaystyle
	    \vspace{0.2cm}
	    \overline{W}\left(\vect{r}_{J}-\vect{r}_{I},\tilde{h}_{IJ}\right) = \\
	    \displaystyle =
	    \frac{\displaystyle
	    W\left(\vect{r}_{J}-\vect{r}_{I},\tilde{h}_{IJ}\right)}{\displaystyle\sum
		\limits_{\vect{r}_{J}\in\Omega_{\vect{r}_{I},\tilde{h}_{I}}}
		W\left(\vect{r}_{J}-\vect{r}_{I},\tilde{h}_{IJ}\right)V_{\mathbf{r}_{J}}} =
	    \frac{\displaystyle W
			\left(\vect{r}_{J}-\vect{r}_{I},\tilde{h}_{IJ}\right)}{\nu\left(\vect{r}_{I}
			\right)}
    \end{array}
\end{equation}
where $\nu\left(\vect{r}_{I}\right)$ is the specific volume of particle
$\vect{r}_{I}$ (i.e., it is approximately the inverse of the particle volume)
which has a larger value in a dense particle region than in a dilute particle
region. In regions of the high particle density, the denominator 
in (\ref{LV_SPH_DN}) is high resulting in lower values of the kernel
$\overline{W}\left(\vect{r}_{J}-\vect{r}_{I},\tilde{h}_{IJ}\right)$. Thus the
denominator normalizes the kernel function to ensure that the kernel
$\overline{W}\left(\vect{r}_{J}-\vect{r}_{I},\tilde{h}_{IJ}\right)$ forms a
local partitioning of unity
\begin{equation}\label{LukyanovVuik:eq:5}
	\displaystyle
	\sum\limits_{\vect{r}_{J}\in\Omega_{\vect{r}_{I},\tilde{h}_{I}}}
	\overline{W}\left(\vect{r}_{J}-\vect{r}_{I},
	\tilde{h}_{IJ}\right)V_{\vect{r}_{J}} = 1
\end{equation}
regardless of the particle distribution within the
$\Omega_{\vect{r}_{I},\tilde{h}_{I}} = \mbox{supp} \
{W}\left(\vect{r}-\vect{r}_{I},\tilde{h}_{I}\right)$.
The discretized-normalized kernel function
$\overline{W}\left(\vect{r}-\vect{r}_{I}, \tilde{h}\right)$
will also be used in the discretization schemes below. At this point, all
possible options of computing
$\overline{\nabla W}(\vect{r}_{J}-\vect{r}_{I},\tilde{h}_{IJ})$ are listed:
\begin{equation}\label{gradKernelOpt1}
	\overline{\nabla W}(\vect{r}_{J}-\vect{r}_{I},\tilde{h}_{IJ}) =
	\nabla_{\vect{r}_{J}} W(\vect{r}_{J}-\vect{r}_{I},\tilde{h}_{IJ}) = -
	\nabla_{\vect{r}_{I}} W(\vect{r}_{J}-\vect{r}_{I},\tilde{h}_{IJ}),
\end{equation}
\begin{equation}\label{gradKernelOpt2}
	\displaystyle
	\begin{array}{lcl}
		\vspace{0.2cm}\displaystyle
		\overline{\nabla W}(\vect{r}_{J}-\vect{r}_{I},\tilde{h}_{IJ}) =
		\nabla_{\vect{r}_{I}}\overline{W}(\vect{r}_{J}-\vect{r}_{I},\tilde{h}_{IJ})
		=\\
		\vspace{0.2cm}\displaystyle =
		\frac{\nabla_{\vect{r}_{I}}
		W(\vect{r}_{J}-\vect{r}_{I},\tilde{h}_{IJ})}{\nu(\mathbf{r}_{I})}-
		\frac{W(\vect{r}_{J}-\vect{r}_{I},\tilde{h}_{IJ})\nabla_{\vect{r}_{I}}
		\nu(\vect{r}_{I})}{\nu^{2}(\vect{r}_{I})},\\
		\displaystyle
		\nabla_{\vect{r}_{I}}\nu(\vect{r}_{I}) =
		\sum\limits_{\displaystyle\vect{r}_{J}\in\Omega_{\vect{r}_{I},\tilde{h}_{I}}}
		\nabla_{\vect{r}_{I}}W(\vect{r}_{J}-\vect{r}_{I},\tilde{h}_{IJ})V_{\vect{r}}
	\end{array}
\end{equation}
where $\displaystyle\nu(\vect{r}_{I})$ is the specific volume of the 
particle located at the point $\displaystyle\vect{r}_{I}$. Additionally 
two options can be written as
\begin{equation}\label{gradKernelOpt3}
	\displaystyle
	\begin{array}{lcl}
		\vspace{0.2cm}\displaystyle
		\overline{\nabla W}(\vect{r}_{J}-\vect{r}_{I},\tilde{h}_{IJ}) =
		\widetilde{\nabla}_{\vect{r}_{J}}\overline{W}(\vect{r}_{J}-
		\vect{r}_{I},\tilde{h}_{IJ})=
		\frac{\nabla_{\vect{r}_{J}}W(\vect{r}_{J}-\vect{r}_{I},
		\tilde{h}_{IJ})}{\nu(\vect{r}_{I})}
	\end{array}
\end{equation}
where $\displaystyle\nu(\vect{r}_{I})$ is assumed to be a constant 
during the differentiation with respect to $\vect{r}_{J}$ and the 
alternative case is when
$$
\displaystyle\nu(\vect{r}_{I},\vect{r}_{J})=
\sum\limits_{\displaystyle\vect{r}_{J}\in
\Omega_{\vect{r}_{I},\tilde{h}_{I}}}W(\vect{r}_{J}-\vect{r}_{I},\tilde{h}_{IJ})
V_{\vect{r}_{J}},
$$
leading to the following relations
\begin{equation}\label{gradKernelOpt4}
	\displaystyle
	\begin{array}{lcl}
		\vspace{0.2cm}\displaystyle
		\overline{\nabla W}(\vect{r}_{J}-\vect{r}_{I},\tilde{h}_{IJ}) =
		\overline{\nabla}_{\vect{r}_{J}}
		\overline{W}(\vect{r}_{J}-\vect{r}_{I},\tilde{h}_{IJ})=\\
		\vspace{0.2cm}\displaystyle=
		\frac{\nabla_{\vect{r}_{J}}W(\vect{r}_{J}-\vect{r}_{I},
		\tilde{h}_{IJ})}{\nu(\vect{r}_{I},\vect{r}_{J})}-
		\frac{W(\vect{r}_{J}-\vect{r}_{I},\tilde{h}_{IJ})
		\nabla_{\vect{r}_{J}}\nu(\vect{r}_{I},\vect{r}_{J})}{\nu^{2}(\vect{r}_{I},
		\vect{r}_{J})}, \\
		\vspace{0.2cm}\displaystyle
		\nabla_{\vect{r}_{J}}\nu(\vect{r}_{I},\vect{r}_{J}) =
		\nabla_{\vect{r}_{J}}W(\vect{r}_{J}-\vect{r}_{I},\tilde{h}_{IJ})
		V_{\vect{r}_{J}}
	\end{array}.
\end{equation}
Where $\nabla_{\vect{r}_{J}}$ denotes nabla operator with 
respect to $\vect{r}_{J}$ and this index is omitted throughout 
this paper starting from here.

Similar to the MLS method \cite{Breitkopf04}, equations
(\ref{gradKernelOpt2}) and (\ref{gradKernelOpt4}) are two different forms of
"full derivatives". At the same time, equations (\ref{gradKernelOpt1}) and
(\ref{gradKernelOpt3}) are two different forms of "diffuse derivatives". From
these options, it follows that "full derivatives" are connected with the
differentiation of a discrete function and "diffuse derivatives" are connected
with the differentiation of an exact function. Both types of derivatives have
some advantages and disadvantages and the choice depends on the application. 
The impact of different options on numerical results will be shown below.
The SPH method shows good approximation properties in regions where the 
kernel has full support. For particles near free surfaces or boundaries, the SPH 
method shows a poor approximation. Several techniques have been developed 
to address these errors through various correction methods, e.g., by applying 
normalized - corrected meshless methods in the 
derivative approximations \cite{BelytschkoKrongauzDolbowGerlach1998,LiuLiu2003, Lukyanov2007,LukyanovPenkov2007}, 
which requires normalized - corrected definitions of 
the kernel gradient as follows:
\begin{equation}\label{LV_SPH_CK}
	\displaystyle\overline{\nabla^{*}_{\alpha} W} =
	\mathbf{C}_{\alpha\beta}\overline{\nabla_{\beta} W},
\end{equation}
\begin{equation}\label{LV_SPH_CT}
	\displaystyle
	\mathbf{C}_{\alpha\beta} =
	\left[
	\sum\limits_{\vect{r}_{J}\in\Omega_{\vect{r}_{I},\tilde{h}_{I}}}V_{\vect{r}_{J}}
	\left[
	\vect{r}^{\alpha}_{J} -
	\vect{r}^{\alpha}_{I}
	\right]
	\overline{\nabla_{\beta} W}\left(\vect{r}_{J}-\vect{r}_{I}\right)
	\right]^{-1},
\end{equation}
\begin{equation}
	\label{LV_SPH_CT1}
	\displaystyle
	\sum\limits_{\Omega_{\vect{r}_{I},h}}V_{\vect{r}_{J}}
	\left[\vect{r}^{\gamma}_{J}-\vect{r}^{\gamma}_{I}\right]
	\overline{\nabla^{*}_{\alpha} W}\left(\vect{r}_{J}-\vect{r}_{I}, h\right) =
	\delta_{\gamma\alpha},\ \ \ \ \forall \gamma, \alpha;
\end{equation}
where the summation by repeating indexes is assumed 
throughout this paper, $\overline{\nabla^{*}_{\alpha} W}$ is the 
normalized-corrected gradient of the kernel, and $\mathbf{C}_{\alpha\beta}$ 
is the correction symmetric tensor \cite{RandlesLibersky1996}. It was shown 
that the value of the minimum eigenvalue $\lambda^{\mathbf{C}}\left(\vect{r}_{I}\right)$ 
of the matrix $\mathbf{C}^{-1}$ based on the discretized-normalized 
kernel function depends on the particle distribution within the 
domain $\Omega_{\vect{r}_{I},\tilde{h}_{I}} =
\mbox{supp}{W}\left(\vect{r}-\vect{r}_{I},\tilde{h}_{I}\right)$.
When going away from the $\Omega_{\vect{r}_{I},\tilde{h}_{I}}$ domain this
eigenvalue tends theoretically to zero, while inside this domain the eigenvalue
tends theoretically to one. This important information allows determining regions
of the continuum media where free-surfaces are located
\cite{MarroneColagrossiLeTouzeGraziani2010}.

In this paper, the discretizations of the Laplace operator are based on the
gradient of the kernel. Several methods of Laplace discretizations
were proposed \cite{ChenBeraunCarney1999,ChenBeraunJih1999a,BonetKulasegaram2000,ColinEgliLin2006}
using second derivatives of the variable $\mathbf{u}$. However, second-order
derivatives can often be avoided entirely if the PDE is written in a 
weak form. It is important to note that approximations using second-order 
derivatives of the kernel are often noisy and sensitive to the particle 
distribution, particularly for spline kernels of lower orders.

\subsection{SPH symmetrization of smoothing length}

The true particle smoothing length $\tilde{h}_{I}$ may vary both in space and 
time in general. Therefore, in general case, each particle has its own smoothing 
length $\tilde{h}_{I}$. Considering the case where $\tilde{h}_{I}\neq \tilde{h}_{J}$ 
for two different interacting particles $I$ and $J$ and the kernel support based 
on $\tilde{h}_{I}$ and located in $I$ which covers the particle $J$ but the kernel 
support located in $J$ does not cover the particle $I$. In this case, the particle 
$I$ acting on particle $J$ (produces, e.g., a flux or a force) without particle $J$ 
acting on the particle $I$, which leads to a violation of fundamental laws (e.g.,
mass conservation or Newton's third law) for a closed system of particles.
This problem has been resolved by introducing the 
symmetrization of the smoothing length. In this study, the 
following symmetrization option is used:
\begin{equation}\label{Sym_Length}
      \displaystyle\vspace{0.2cm}
      \tilde{h}_{IJ} = \frac{\tilde{h}_{I}+\tilde{h}_{J}}{2}.
\end{equation}
In addition, it is clear that $\tilde{h}_{I}$ has to be defined as
\begin{equation}\label{EffPartSmoothing}
   \displaystyle
   \tilde{h}_{I} = \sup\limits_{J:\vect{r}_{j}\in
   \Omega_{\vect{r}_{I},\tilde{h}_{I}}}\tilde{h}_{IJ}.
\end{equation}
This completes the description of basic properties of SPH method 
allowing to construct all necessary building elements of SPH discretization 
such as list of neighbors, kernel values, and kernel gradients. The following 
sections describe the traditional and newly proposed discretization schemes 
for the Laplace operator.

\subsection{Brookshaw's scheme (1985)}
\label{BrookshawScheme}

Brookshaw proposed \cite{Brookshaw1985} an approximation of the Laplacian 
for an inhomogeneous scalar field $m\left(\vect{r}\right)$,
i.e., $\mathbf{M}^{\alpha\beta}\left(\vect{r}\right) = 
m\left(\vect{r}\right)\delta_{\alpha\beta}, \  \alpha,\beta = 1,\ldots,n;$
that only includes first order derivatives:
\begin{equation}
	\label{Brookshaw1}
	\begin{array}{lcl}
		\displaystyle
		\vspace{0.3cm}
		-\langle\nabla\left(m\left(\vect{r}_{I}\right)\nabla
		\mathbf{u}\left(\vect{r}_{I}\right)\right)\rangle = \\
		\displaystyle 
		\sum\limits_{\Omega_{\vect{r}_{I},\tilde{h}_{I}}}V_{\vect{r}_{J}}
		\left[\mathbf{u}\left(\vect{r}_{J}\right)-\mathbf{u}\left(\vect{r}_{I}\right)
		\right]
		\frac{\left(\vect{r}_{J}-\vect{r}_{I}\right)\cdot\left(m_{J}+m_{I}\right)
		\overline{\nabla W}
		\left(\vect{r}_{J}-\vect{r}_{I},\tilde{h}_{IJ}\right)}{\left\Vert\vect{r}_{J}-\vect{r}_{I}
		\right\Vert^{2}},
	\end{array}
\end{equation}
where $V_{\vect{r}_{J}}$ is the volume of the particle $J$, 
$\left\Vert\bullet\right\Vert$ is the Euclidean norm throughout this paper,
$\mathbf{u}\left(\vect{r}\right)$ is the unknown scalar or vector field (e.g., 
pressure $p$ or velocity $\vect{v}$) $\forall\vect{r}\in\Omega\subset\mathbb{R}^{n}$,
$m_{I}=m\left(\vect{r}_{I}\right)$, $\vect{r}_{I}\in\Omega\subset\mathbb{R}^{n}$
and $m_{J}=m\left(\vect{r}_{J}\right)$, $\vect{r}_{J}\in\Omega\subset\mathbb{R}^{n}$
are the field coefficients.

This scheme can be derived by applying the particle approximation step to 
the right-hand side of \eqref{LV_SPH_ID} with the following assumptions
 \begin{equation}
	 \label{LV_SPH_Brookshaw_Assump1}
	 \displaystyle
	 \nabla \mathbf{u}\left(\vect{r}_{J}\right) \approx
	 \left[\mathbf{u}\left(\vect{r}_{J}\right)-\mathbf{u}\left(\vect{r}_{I}\right)
	 \right]
	 \frac{\left(\vect{r}_{J}-\vect{r}_{I}\right)}{\left\Vert\vect{r}_{J}-
	 \vect{r}_{I}\right\Vert^{2}}.
 \end{equation}
 Some special words need to be said about the mobility approximation, which comes
 in the form
 \begin{equation}
	 \label{LV_SPH_Brookshaw_Assump2}
	 \displaystyle 2m\left(\vect{r}_{J}\right)\approx m_{J}+m_{I}.
 \end{equation}
 The factor 2 is introduced to compensate the factor of 1/2 in the second leading
 term of the Taylor expansion of the relation 
 \eqref{LV_SPH_Brookshaw_Assump1}. Furthermore, the relation 
 \eqref{LV_SPH_Brookshaw_Assump2} allows to capture a 
 heterogeneous mobility field distribution.

\begin{figure}[!htbp]
	\centering
	\includegraphics[height=0.3\textheight]{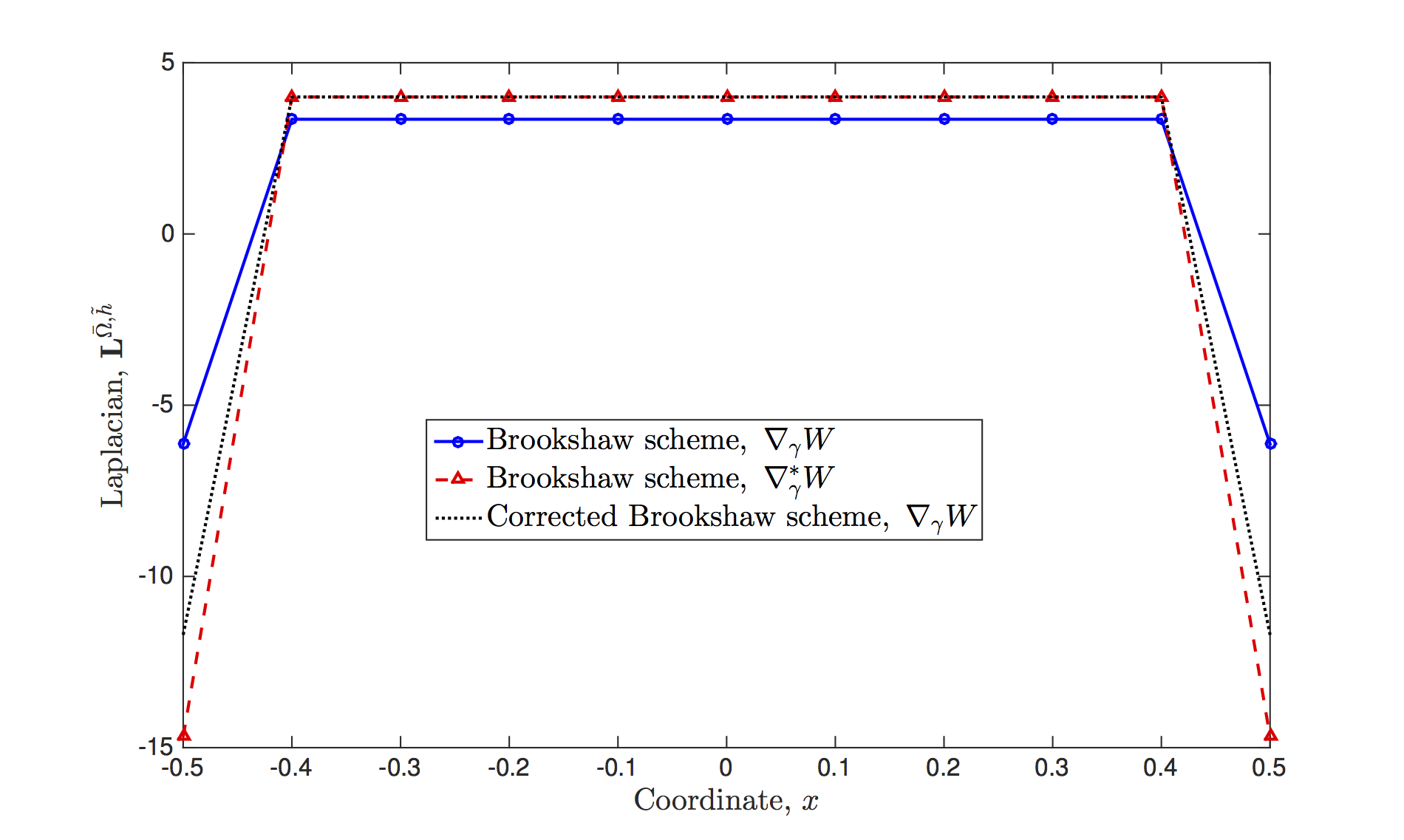}
	\caption{Values for $\nabla^2(x^2+y^2) $ along $y=0$ using Brookshaw's
		approximation with (a) conventional kernel $\nabla_{\gamma}W$, (b)
		corrected kernel $\nabla^{*}_{\gamma}W$ and with the correction
		multiplier $\displaystyle  n\cdot\left[\mathbf{C}_{\alpha\alpha}\right]^{-1}$.
		The numerical domain is a unit square in $\mathbb{R}^{2}$ with the
		center at $a_{i}=0, \ \forall i$ and side length $L=1$. The cubic spline
		\eqref{LV_SPH_K} was used with $\tilde{h}=f\cdot h_{p}$,
	    $h_{p}=0.1$, $f=1.0$.}
	    \label{BrookshawFigure1} 
\end{figure}

The discretization scheme (\ref{Brookshaw1}) is at least
$\mathcal{O}\left(h^\omega\right), \ 0\leq\omega < 2$, 
$h = \sup\limits_{I:\vect{r}_{I}\in\Omega}\tilde{h}_{I}$
order of accuracy in average for any scalar mobility field
$m\left(\vect{r}\right)\in C^{1}\left(\Omega\right), \ m\left(\vect{r}\right)\geq 0$
everywhere within the numerical domain $\Omega$ sufficiently
far away from the boundary $\partial\Omega$.
Using Taylor series expansions about a point $\vect{r}_{I}$,
the following relations can be written:
\begin{equation}
	\label{Brookshaw2}
	\begin{array}{lcl}
		\displaystyle
		\mathbf{u}\left(\vect{r}_{J}\right) =
		\mathbf{u}\left(\vect{r}_{I}\right) +
		\mathbf{u}_{,\alpha}\left(\vect{r}_{I}\right)
		\left[\vect{r}^{\alpha}_{J}-\vect{r}^{\alpha}_{I}\right]
		\displaystyle +
		\frac{1}{2}\mathbf{u}_{,\alpha\gamma}\left(\vect{r}_{I}\right)
		\left[\vect{r}^{\alpha}_{J}-\vect{r}^{\alpha}_{I}\right]
		\left[\vect{r}^{\gamma}_{J}-\vect{r}^{\gamma}_{I}\right] +
		\mathcal{O}\left(h^{3}\right),
	\end{array}
\end{equation}
\begin{equation}
	\label{Brookshaw3}
	\displaystyle
	m\left(\vect{r}_{J}\right) =
	m\left(\vect{r}_{I}\right) +
	m_{,\alpha}\left(\vect{r}_{I}\right)
	\left[\vect{r}^{\alpha}_{J}-\vect{r}^{\alpha}_{I}\right]+
	\mathcal{O}\left(h^{2}\right).
\end{equation}
Substituting relations (\ref{Brookshaw2})-(\ref{Brookshaw3}) into
the scheme (\ref{Brookshaw1}), it leads to the following relations:
\begin{equation}
\label{Brookshaw4}
	\begin{array}{lcl}
		\displaystyle
		\vspace{0.3cm}
		\sum\limits_{\Omega_{\vect{r}_{I},\tilde{h}_{I}}}V_{\vect{r}_{J}}
		\left[\mathbf{u}\left(\vect{r}_{J}\right)-\mathbf{u}\left(\vect{r}_{I}\right)
		\right]
		\frac{\left(\vect{r}_{J}-\vect{r}_{I}\right)\cdot\left(m_{J}+m_{I}\right)
		\cdot
		\overline{\nabla W}
		\left(\vect{r}_{J}-\vect{r}_{I},\tilde{h}_{IJ}\right)}{\left\Vert\vect{r}_{J}-\vect{r}_{I}
		\right\Vert^{2}}=\\
		\vspace{0.3cm}
		\displaystyle=2m\left(\vect{r}_{I}\right)\mathbf{u}_{,\alpha}\left(\vect{r}_{I}
		\right)
		\sum\limits_{\Omega_{\vect{r}_{I},\tilde{h}_{I}}}V_{\vect{r}_{J}}
		\overline{\nabla_{\alpha}W}\left(\vect{r}_{J}-\vect{r}_{I}, \tilde{h}_{IJ}\right)+\\
		\vspace{0.3cm}
		\displaystyle+
		m\left(\vect{r}_{I}\right)\mathbf{u}_{,\alpha\gamma}\left(\vect{r}_{I}\right)
		\sum\limits_{\Omega_{\vect{r}_{I},h}}V_{\vect{r}_{J}}
		\left[\vect{r}^{\alpha}_{J}-\vect{r}^{\alpha}_{I}\right]
		\overline{\nabla_{\gamma}W}\left(\vect{r}_{J}-\vect{r}_{I},\tilde{h}_{IJ}\right)+\\
		\displaystyle+
		m_{,\alpha}\left(\vect{r}_{I}\right)\mathbf{u}_{,\gamma}\left(\vect{r}_{I}
		\right)
		\sum\limits_{\Omega_{\vect{r}_{I},h}}V_{\vect{r}_{J}}
		\left[\vect{r}^{\alpha}_{J}-\vect{r}^{\alpha}_{I}\right]
		\overline{\nabla_{\gamma}W}\left(\vect{r}_{J}-\vect{r}_{I}, \tilde{h}_{IJ}\right)+\mathcal{O}\left(h^{2}\right).
	\end{array}
\end{equation}
Here, the following relation has been used
\begin{equation}
	\label{Brookshaw5}
	\left[\vect{r}^{\alpha}_{J}-\vect{r}^{\alpha}_{I}\right]
	\frac{\left(\vect{r}^{\gamma}_{J}-\vect{r}^{\gamma}_{I}\right)
	\overline{\nabla_{\gamma}W}
	\left(\vect{r}_{J}-\vect{r}_{I}, \tilde{h}_{IJ}\right)}{\left\Vert\vect{r}_{J}-\vect{r}_{I}
	\right\Vert^{2}} =
	\overline{\nabla_{\alpha}W}\left(\vect{r}_{J}-\vect{r}_{I}, \tilde{h}_{IJ}\right),
	\ \ \ \ \forall \alpha.
\end{equation}
The maximum accuracy is achieved when
\begin{equation}
	\label{Brookshaw6}
	\begin{array}{lcl}
		\displaystyle
		\vspace{0.2cm}
		(a) \ \sum\limits_{\Omega_{\vect{r}_{I},\tilde{h}_{I}}}V_{\vect{r}_{J}}
		\overline{\nabla_{\alpha}W}\left(\vect{r}_{J}-\vect{r}_{I}, \tilde{h}_{IJ}\right) = 0,
		\ \ \ \ \forall \alpha,\\
		\displaystyle
		(b) \ \sum\limits_{\Omega_{\vect{r}_{I},\tilde{h}_{I}}}V_{\vect{r}_{J}}
		\left[\vect{r}^{\alpha}_{J}-\vect{r}^{\alpha}_{I}\right]
		\overline{\nabla_{\gamma}W}\left(\vect{r}_{J}-\vect{r}_{I}, \tilde{h}_{IJ}\right) =
		\delta_{\alpha\gamma},
		\ \ \ \ \forall \alpha,\gamma,
	\end{array}
\end{equation}
which is difficult to fulfill simultaneously for different kernel gradients
leading to the overall accuracy $\mathcal{O}\left(h^\omega\right), \ 0\leq\omega < 2$.
In the subsection \ref{SPH_Kernel_Grad}, several options of computing kernel
gradients $\overline{\nabla W}(\vect{r}_{J}-\vect{r}_{I},\tilde{h}_{IJ})$
(i.e., $\nabla_{\gamma}W$, $\nabla_{\alpha}\overline{W}$,
$\overline{\nabla}_{\alpha}\overline{W}$, and
$\widetilde{\nabla}_{\alpha}\overline{W}$) are proposed \eqref{gradKernelOpt1},
\eqref{gradKernelOpt2}, \eqref{gradKernelOpt3}, and \eqref{gradKernelOpt4},
respectively. The kernel gradient \eqref{gradKernelOpt2} satisfies conditions
\eqref{Brookshaw5} and \eqref{Brookshaw6}(a) but not the condition
\eqref{Brookshaw6}(b). At the same time, all corrected options
$\overline{\nabla^{*} W}(\vect{r}_{J}-\vect{r}_{I},\tilde{h}_{IJ})$
(i.e., $\nabla^{*}_{\gamma}W$, $\nabla^{*}_{\alpha}\overline{W}$,
$\overline{\nabla}^{*}_{\alpha}\overline{W}$, and
$\widetilde{\nabla}^{*}_{\alpha}\overline{W}$) satisfy the condition
\eqref{Brookshaw6}(b) but not the conditions
\eqref{Brookshaw5}. One may decide to use
$\overline{\nabla^{*} W}(\vect{r}_{J}-\vect{r}_{I},\tilde{h}_{IJ})$ in the
discretization scheme \eqref{Brookshaw1} which leads to the error with
the leading term
\begin{equation}
	\label{Brookshaw7}
    2m\left(\vect{r}_{I}\right)\mathbf{u}_{,\alpha}\left(\vect{r}_{I}\right)
	\sum\limits_{\vect{r}_{J}\in\Omega_{\vect{r}_{I},\tilde{h}_{I}}}V_{\vect{r}_{J}}
	\left[\vect{r}_{J}-\vect{r}_{I}\right]
	\frac{\left(\vect{r}_{J}-\vect{r}_{I}\right)
	\cdot\overline{\nabla^{*} W}
	(\vect{r}_{J}-\vect{r}_{I},\tilde{h}_{IJ})}{\left\Vert\vect{r}_{J}-\vect{r}_{I}
	\right\Vert^{2}}.
\end{equation}
This leads to the incorporation of the correction factor into Brookshaw's
approximation. A different correction factor has been introduced and
investigated in \cite{Schwaiger2008}. However, the discretization scheme
\eqref{Brookshaw1} with the 
$\overline{\nabla^{*} W}(\vect{r}_{J}-\vect{r}_{I},\tilde{h}_{IJ})$ kernel
gradient is less accurate near the boundary due to the absence of 
the skew symmetric property of the gradient and 
remaining singularity in the error term when $\vect{r}_{J}=\vect{r}_{I}$ 
which can be removed in the conventional form \eqref{Brookshaw4}. 
Alternatively, as was discussed in \cite{Schwaiger2008}, the correction 
multiplier can be introduced in \eqref{Brookshaw1} defined as
$\displaystyle n\cdot\left[\mathbf{C}_{\alpha\alpha}\right]^{-1}$ leading to:
\begin{equation}
	\label{Brookshaw1M}
	\displaystyle
	\vspace{0.3cm}
	\langle\nabla\left(m\left(\vect{r}_{I}\right)\nabla
	\mathbf{u}\left(\vect{r}_{I}\right)\right)\rangle^{*} =
	n\cdot\left[\mathbf{C}_{\alpha\alpha}\right]^{-1}\cdot
	\langle\nabla\left(m\left(\vect{r}_{I}\right)\nabla
	\mathbf{u}\left(\vect{r}_{I}\right)\right)\rangle.
\end{equation}
Figure \ref{BrookshawFigure1} shows Laplacian values for the function
$\nabla^2(x^2+y^2)$ using original Brookshaw's approximation \eqref{Brookshaw1}
with (a) conventional kernel $\nabla_{\gamma}W$, (b) corrected kernel
$\nabla^{*}_{\gamma}W$ and corrected Brookshaw's approximation
\eqref{Brookshaw1M} with the conventional kernel 
$\nabla_{\gamma}W$.

The scheme \eqref{Brookshaw1} is widely (almost unconditionally) used in 
the SPH modeling community. For example, it was used for a thermal conduction 
\cite{Brookshaw1985, ClearyMonaghan1999,JubelgasSpringelDolag2004},
for modeling a viscous diffusion \cite{MorrisFoxZhu1997}, for a vortex spin-down 
\cite{CumminsRudman1999} and Rayleigh-Taylor instability, for simulating Newtonian 
and non-Newtonian flows with a free surface \cite{LeeMoulinecXuVioleauLaurenceStansby2008} 
for the comparison of weakly compressible and truly incompressible algorithms, for 
macroscopic and mesoscopic flows \cite{HuAdams2006}, for a simulation of a solid-fluid 
mixture flow \cite{ZhangKuwabaraSuzukiKawanoMoritaFukuda2009}. 
Recently, it has been used to model electrokinetic flows \cite{PanKimPeregoTartakovskyParks2017}, 
Dam-break problem and Taylor-Green vortex \cite{ZhangHuAdams2017}.
  
There are different numerical SPH schemes used in numerical simulations. High order accuracy 
approximations can also be derived by using the SPH discretization on the higher 
order Taylor series expansion \cite{FulkQuinn1996,LiuLiu2003, Schwaiger2008,FangParriaux2008}. 
However, it is usually required that the discrete numerical schemes can reproduce linear
fields \cite{RandlesLibersky1996, BonetLok1999, Lukyanov2007,LukyanovPenkov2007}
or polynomials up to a given order \cite{LiuJunZhang1995,LiuJun1998}.

\subsection{Schwaiger's scheme (2008)}
\label{SchwaigerScheme}

The correction terms to the Brookshaw formulation which improve 
the accuracy of the Laplacian operator near boundaries were proposed by 
Schwaiger in \cite{Schwaiger2008}:
\begin{equation}
	\label{Schwaiger1}
	\begin{array}{lcl}
		\displaystyle
		\vspace{0.3cm}
		-\frac{n}{\Gamma^{-1}_{kk}}\langle
		\nabla\left(m\left(\vect{r}_{I}\right)\nabla
		\mathbf{u}\left(\vect{r}_{I}\right)\right)
		\rangle = \\
		\vspace{0.3cm}
		\displaystyle\left\lbrace
		\sum_{\Omega_{\vect{r}_{I},\tilde{h}_{I}}}V_{\vect{r}_{J}}
		\left[\mathbf{u}\left(\vect{r}_{J}\right)-\mathbf{u}\left(\vect{r}_{I}
		\right)\right]
		\frac{\left(\vect{r}_{J}-\vect{r}_{I}\right)\cdot\left(m_{J}+m_{I}\right)
		\overline{\nabla W}
		\left(\vect{r}_{J}-\vect{r}_{I}, \tilde{h}_{IJ}\right)}{\left\Vert\vect{r}'-\vect{r}
		\right\Vert^{2}}
		\right\rbrace-\\-
		\displaystyle
		\left\lbrace
		\left[
		\langle\nabla_{\alpha}\left(m\left(\vect{r}_{I}\right)\mathbf{u}
		\left(\vect{r}_{I}\right)\right)\rangle-
		\mathbf{u}\left(\vect{r}_{I}\right)\langle\nabla_{\alpha}m
		\left(\vect{r}_{I}\right)\rangle+
		m\left(\vect{r}_{I}\right)\langle\nabla_{\alpha}
		\mathbf{u}\left(\vect{r}_{I}\right)\rangle
		\right]\mathbf{N}^{\alpha}\right\rbrace,
	\end{array}
\end{equation}
\begin{equation}
	\label{Schwaiger2}
	\displaystyle
	\mathbf{N}^{\alpha}\left(\vect{r}_{I}\right) =
	\left[\sum_{\Omega_{\vect{r}_{I},\tilde{h}_{I}}}V_{\vect{r}_{J}}
	\overline{\nabla_{\alpha} W}\left(\vect{r}_{I}-\vect{r}_{J},\tilde{h}_{IJ}\right)
	\right],
\end{equation}
\begin{equation}
	\label{Schwaiger3}
	\displaystyle
	\langle
	\nabla_{\alpha}\mathbf{u}\left(\vect{r}_{I}\right)
	\rangle =
	\sum_{\Omega_{\vect{r}_{I},\tilde{h}_{I}}}V_{\vect{r}_{J}}
	\left[
		\mathbf{u}\left(\vect{r}_{J}\right)-
		\mathbf{u}\left(\vect{r}_{I}\right)
	\right]
	\overline{\nabla^{*}_{\alpha} W}\left(\vect{r}_{I}-\vect{r}_{J},\tilde{h}_{IJ}\right),
\end{equation}
where $n=1,2,3$ is the spatial dimension and the tensor 
$\Gamma_{\alpha\beta}$ is defined by
\begin{equation}
	\label{Schwaiger5}
	\displaystyle
	\Gamma_{\alpha\beta}\left(\mathbf{r}_{I}\right) =
	\displaystyle\sum_{\Omega_{\vect{r}_{I},\tilde{h}_{I}}}V_{\vect{r}_{J}}
	\frac{\left(\vect{r}^{\gamma}_{J}-\vect{r}^{\gamma}_{I}\right)
	\overline{\nabla_{\gamma}W}\left(\vect{r}_{J}-\vect{r}_{I},\tilde{h}_{IJ}\right)}{\left
	\Vert\vect{r}'-\vect{r}\right\Vert^{2}}
	\left(\vect{r}^{\alpha}_{J}-\vect{r}^{\alpha}_{I}\right)
	\left(\vect{r}^{\beta}_{J}-\vect{r}^{\beta}_{I}\right).
\end{equation}
The gradient $\langle\nabla_{\alpha}\mathbf{u}\left(\vect{r}_{I}\right)\rangle$
is the corrected gradient which can reproduce linear fields
\cite{RandlesLibersky1996,BonetLok1999}.
For multi-dimensional problems, the correction tensor
$\Gamma_{\alpha\beta}\left(\mathbf{r}_{I}\right)$ is a matrix.
If the particle $\mathbf{r}_{I}$ has entire stencil support (i.e.,
the domain support for all kernels 
$W\left(\mathbf{r}_{J}-\mathbf{r}_{I}, \tilde{h}_{IJ}\right)$
is entire and symmetric) then
$\Gamma_{\alpha\beta}\left(\mathbf{r}_{I}\right)\approx\delta_{\alpha\beta}$.

\begin{remark}\label{Remark_Brookshaw1}
	It is important to note that correction tensors $\Gamma_{\alpha\beta}$ and
	$\mathbf{C}^{-1}_{\alpha\beta}$ are the same tensors. Indeed, using the
	following identity:
	\begin{equation}
	\label{Schwaiger6}
	\begin{array}{lcl}
		\displaystyle
		\vspace{0.3cm}
		\left[\mathbf{r}^{\alpha}_{J}-\mathbf{r}^{\alpha}_{I}\right]
		\frac{\left(\mathbf{r}^{\gamma}_{J}-\mathbf{r}^{\gamma}_{I}\right)
		\overline{\nabla_{\gamma} W}
		\left(\vect{r}_{J}-\vect{r}_{I}, \tilde{h}_{IJ}\right)}{\left\Vert\vect{r}_{J}-
		\vect{r}_{I}\right\Vert^{2}} = \\
		\displaystyle=
		\frac{1}{h}\frac{\overline{dW}}{dz}
		\frac{\left[\vect{r}^{\alpha}_{J}-\vect{r}^{\alpha}_{I}\right]}{\left
		\Vert\vect{r}_{J}-\vect{r}_{I}\right\Vert}=
		\overline{\nabla_{\alpha} W}\left(\vect{r}_{J}-\vect{r}_{I}, \tilde{h}_{IJ}\right),
		\ \ \forall \alpha; \ \frac{\overline{dW}}{dz} \le 0,
	\end{array}
	\end{equation}
	where $\displaystyle\frac{\overline{dW}}{dz}$ is computed using either
	conventional $W$ or normalized $\overline{W}$ kernels,
	$\displaystyle z=\left\|\vect{r}_{J}-\vect{r}_{I}\right\|/\tilde{h}$,
	$\forall  \vect{r}_{J},\vect{r}_{I}\in\Omega\subset\mathbb{R}^{n}$, the
	following relation can be established:
	\begin{equation}
		\label{Schwaiger7}
		\begin{array}{lcl}
			\displaystyle
			\vspace{0.3cm}
			\Gamma_{\alpha\beta}\left(\vect{r}_{I}\right) =
			\sum\limits_{\Omega_{\vect{r}_{I},\tilde{h}_{I}}}V_{\vect{r}_{J}}
			\frac{\left(\vect{r}^{\gamma}_{J}-\vect{r}^{\gamma}_{I}\right)
			\overline{\nabla_{\gamma}W}
			\left(\vect{r}_{J}-\vect{r}_{I}, \tilde{h}_{IJ}\right)}{\left\Vert\vect{r}_{J}-
			\vect{r}_{I}\right\Vert^{2}}
			\left(\vect{r}^{\alpha}_{J}-\vect{r}^{\alpha}_{I}\right)
			\left(\vect{r}^{\beta}_{J}-\vect{r}^{\beta}_{I}\right) = \\
			\displaystyle=\sum\limits_{\Omega_{\vect{r}_{I},\tilde{h}_{I}}}V_{\vect{r}_{J}}
			\left[
				\vect{r}^{\alpha}_{J} -
				\vect{r}^{\alpha}_{I}
			\right]
			\overline{\nabla_{\beta} W}\left(\vect{r}_{J}-\vect{r}_{I}, \tilde{h}_{IJ}\right)=
			\mathbf{C}^{-1}_{\alpha\beta}\left(\vect{r}_{I}\right).
		\end{array}
	\end{equation}
\end{remark}
In addition, it is important to note that
\begin{equation}\label{Schwaiger8}
	\displaystyle
	\Gamma_{\alpha\alpha} =
	\sum_{\Omega_{\vect{r}_{I},\tilde{h}_{I}}}V_{\vect{r}_{J}}
	\left(\vect{r}^{\gamma}_{J}-\vect{r}^{\gamma}_{I}\right)\cdot
	\overline{\nabla^{*}_{\gamma}W}\left(\vect{r}_{J}-\vect{r}_{I}, \tilde{h}_{IJ}\right) = n
\end{equation}
in the case of using the corrected gradient and, hence,
$\displaystyle\frac{\Gamma^{-1}_{kk}}{n}=1$. However,
$\overline{\nabla^{*}_{\gamma}W}$ does not satisfy relation \eqref{Schwaiger6}.

For multi-dimensional problems, the correction tensor
$\Gamma_{\alpha\beta}\left(\vect{r}_{I}\right)$ is a matrix. If the particle
$\vect{r}_{I}$ has entire stencil support (i.e., the domain support for all
kernels $W\left(\vect{r}_{J}-\vect{r}_{I},\tilde{h}_{IJ}\right)$ is entire and
symmetric) then
$\Gamma_{\alpha\beta}\left(\vect{r}_{I}\right) \approx\delta_{\alpha\beta}$. 
Unfortunately, $\Gamma_{\alpha\beta}\left(\vect{r}_{I}\right)$ deviates
from $\delta_{\alpha\beta}$ for the provided algorithm and, hence, it is important
to minimize this deviation from $\delta_{\alpha\beta}$ in the new methods.

To calculate coefficients in the scheme (\ref{Schwaiger1})--(\ref{Schwaiger5}) is
a trivial task. However, in general, it should be performed at each Newton-Raphson
iteration in the non-linear case (i.e., $m = m \left(\mathbf{u}\left(\vect{r}_{I}
\right)\right)$). It also requires additional efforts to invert the correction
matrix $\mathbf{A}_{\alpha\beta}$ (inversion of $n\times n$ matrices per
each particle, where $n=1,2,3$ is the spatial dimension) and storage cost of
$\overline{\nabla_{\alpha} W}\left(\mathbf{r}_{J}-\mathbf{r}_{I}, h\right)$,
$\overline{\nabla^{*}_{\alpha} W}\left(\mathbf{r}_{J}-\mathbf{r}_{I}, h\right)$,
and corresponding $\Gamma^{-1}_{\alpha\alpha}=\mathbf{A}^{-1}_{\alpha\alpha}$
per each particle.

\begin{figure}[!htbp]
	\centering
	\includegraphics[height=0.3\textheight]{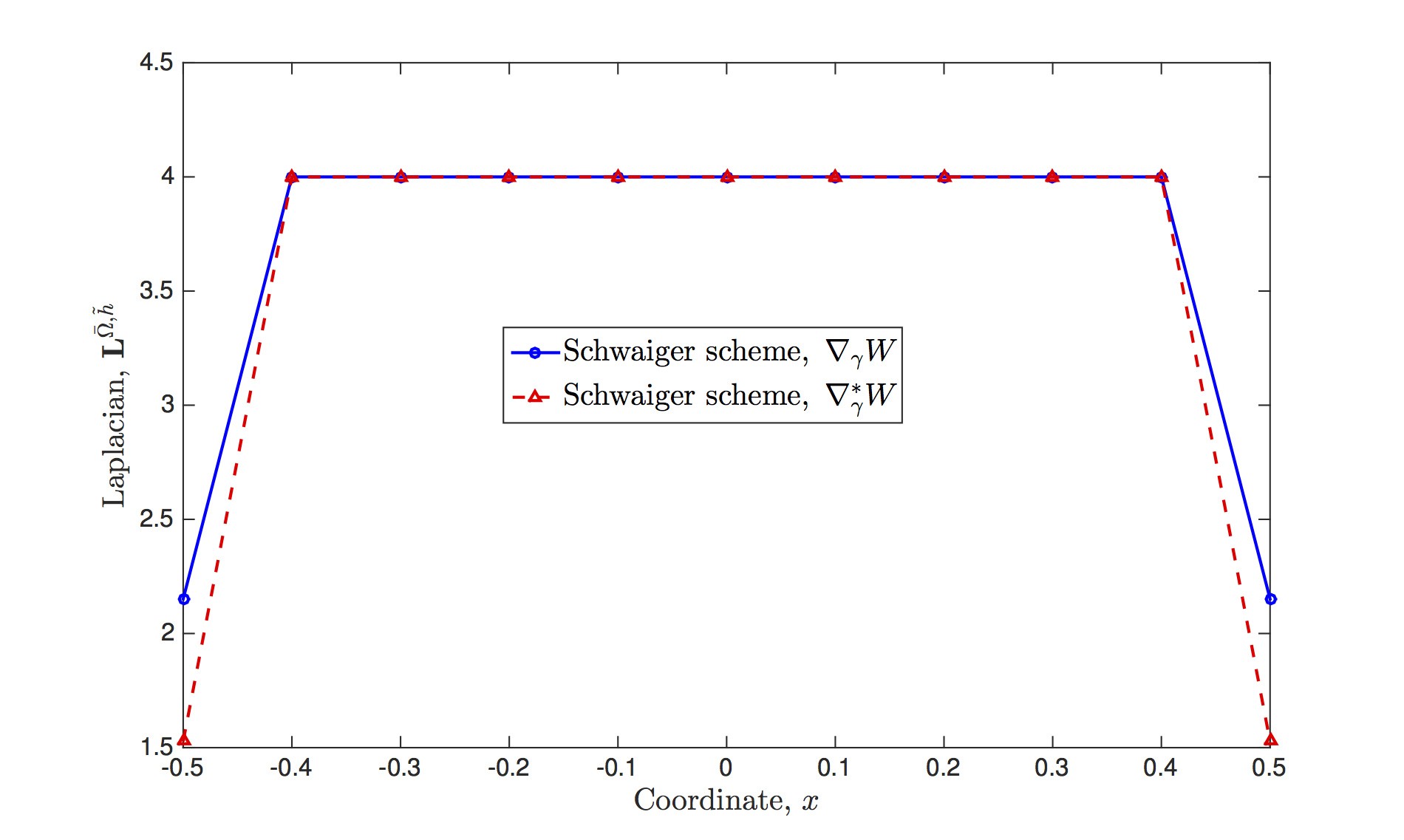}
	\caption{Values for $\nabla^2(x^2+y^2)$ along $y=0$ using Schwaiger's
	approximation with (a) conventional kernel $\nabla_{\gamma}W$ and (b)
	corrected kernel $\nabla^{*}_{\gamma}W$. The numerical domain is a unit
	square in in $\mathbb{R}^{2}$ with the center at $a_{i}=0, \ \forall i$. The 
	cubic spline \eqref{LV_SPH_K} was used with $\tilde{h}=f\cdot h_{p}$, 
	$h_{p}=0.1$, $f=1.0$.}\label{SchwaigerFigure1} 
\end{figure}

Furthermore, additional terms proposed by Schwaiger \cite{Schwaiger2008} 
reduce to
\begin{equation}\label{LV_SPH_Schwaiger4}
   \begin{array}{lcl}
      \displaystyle
      \vspace{0.2cm}
      \left[
      \left\langle
      m\left(\vect{r}_{I}\right)\mathbf{u}\left(\vect{r}_{I}\right)
	  \right\rangle-
      \mathbf{u}\left(\vect{r}_{I}\right)\left\langle m\left(\vect{r}_{I}\right)
	  \right\rangle+
      m\left(\vect{r}_{I}\right)\left\langle\mathbf{u}\left(\vect{r}_{I}\right)
	  \right\rangle
      \right]\mathbf{N} = \\
      \displaystyle=2m
      \left(\vect{r}_{I}\right)\nabla\mathbf{u}\left(\vect{r}_{I}\right)
	  \cdot\mathbf{N}+
	  \mathcal{O}\left(\tilde{h}_{I}^{2}\right)
   \end{array}
\end{equation}
which is the leading term outlined in \eqref{Brookshaw4}. However, if one
uses $\overline{\nabla^{*} W}(\vect{r}_{J}-\vect{r}_{I},\tilde{h}_{IJ})$ in the
first term of the discretization scheme \eqref{Schwaiger1} then the definition
for $\mathbf{N}^{\alpha}$ has to be modified in accordance of \eqref{Brookshaw7}
to maintain the higher order discretization accuracy, for 
example, as:
\begin{equation}\label{LV_SPH_Schwaiger5}
\displaystyle
\widetilde{\mathbf{N}}^{\alpha}=\sum\limits_{\Omega_{\vect{r}_{I},\tilde{h}_{I}}}V_{\vect{r}_{J}}
\left[\vect{r}^{\alpha}_{J}-\vect{r}^{\alpha}_{I}\right]
\frac{\left(\vect{r}_{J}-\vect{r}_{I}\right)
\overline{\nabla^{*} W}\left(\vect{r}_{J}-\vect{r}_{I}, \tilde{h}_{IJ}\right)}{\left
\Vert\vect{r}_{J}-\vect{r}_{I}\right\Vert^{2}},
\ \ \ \ \forall \alpha
\end{equation}
which reduces to the conventional $\mathbf{N}^{\alpha}$ in 
the case of $\overline{\nabla W}(\vect{r}_{J}-\vect{r}_{I},\tilde{h}_{IJ})$ due to
\eqref{Schwaiger6}. Figure \ref{SchwaigerFigure1} shows the Laplacian for the
function $\nabla^2(x^2+y^2) $ using Schwaiger's approximation
(\ref{Schwaiger1})--(\ref{Schwaiger5}) with (a) conventional kernel
$\nabla_{\gamma}W$ and (b) corrected kernel $\nabla^{*}_{\gamma}W$.

Finally, the discretization scheme (\ref{Schwaiger1})--(\ref{Schwaiger5}) is
at least $\mathcal{O}\left(h^\omega\right), \ 1\leq\omega \leq 2$ 
order of accuracy in average for any scalar mobility field
$m\left(\vect{r}\right)\in C^{1}\left(\Omega\right), \ m\left(\vect{r}\right)\geq 0$ 
everywhere within the numerical domain $\Omega\in \mathbb{R}^{n}$ sufficiently 
far away from the boundary $\partial\Omega$.

\section{Meshless Transmissibilities}
\label{MeshlessTransmissibilities}

The well-known two-point flux approximation (TPFA) is a mesh 
dependent numerical scheme used in solving elliptic equation 
(\ref{LV_SPH_Lap}): $\mathbf{L}\left(\mathbf{u}\right) = 0$ with 
the diagonal matrix of coefficients $\mathbf{M}$. The net flow 
rate from a cell $I$ into neighboring cells in this scheme is obtained 
by summing fluxes over the neighboring cells $J$:
\begin{equation}\label{LV_TotalFluxDiscrForm}
	\displaystyle
	\vect{q} =
	\sum\limits_{J} \widetilde{T}_{JI}
	\left[
	\mathbf{u}\left(\vect{r}_{J}\right)-
	\mathbf{u}\left(\vect{r}_{I}\right)
	\right],
	\ \widetilde{T}_{JI} \geq 0,
\end{equation}
where $\widetilde{T}_{JI}$ is the transmissibility between cells $J$ and
$I$, $\vect{q}$ is the total flux through the boundary of the control
volume located at the point $\vect{r}_{I}$. The transmissibility
$\widetilde{T}_{JI}$ defined at an interior face $f$ between cells
$J$ and $I$ is calculated as
\begin{equation}\label{LV_TotalFluxDiscrForm2}
	\displaystyle \widetilde{T}_{JI} =
	\frac{1}{\left[\displaystyle\frac{\left\Vert\vect{r}_{f,J}
	\right\Vert^{2}}{\vect{S}_{f}\mathbf{M}\vect{r}_{f,J}}+
	\frac{\left\Vert\vect{r}_{f,I}\right\Vert^{2}}{\vect{S}_{f}
	\mathbf{M}\vect{r}_{f,I}}\right]},
\end{equation}
where $\vect{r}_{f,J}$ and $\vect{r}_{f,I}$ are the vectors from centers
of cells $J$ and $I$ to the face $f$ respectively, $\vect{S}_{f}$ is the
area vector of the face $f$. In the case of $\mathbf{M}$-orthogonal
mesh, when $\mathbf{M}\vect{S}_{f}$ and $\left[\vect{r}_{J}-\vect{r}_{I}\right]$
are collinear, the expression (\ref{LV_TotalFluxDiscrForm}) reduces to the form
of the central finite difference scheme and approximates the flux with
$\mathcal{O}\left(h^{2}\right)$ order of accuracy for any mobility tensor
field $\mathbf{M}$. The expression (\ref{LV_TotalFluxDiscrForm2}) ensures
that the flux into the adjoining region is continuous \cite{ClearyMonaghan1999}.
The TPFA scheme (\ref{LV_TotalFluxDiscrForm}) is the unconditionally 
monotone scheme.

It is clear that the expression (\ref{Schwaiger1}) cannot be written in the form
(\ref{LV_TotalFluxDiscrForm}) due to terms
$\displaystyle\langle\nabla_{\alpha}
\left(m\left(\vect{r}_{I}\right)\mathbf{u}
\left(\vect{r}_{I}\right)\right)\rangle\mathbf{N}^{\alpha}$ and
$\displaystyle\mathbf{u}\left(\vect{r}_{I}\right)\langle\nabla_{\alpha}m
\left(\vect{r}_{I}\right)\rangle\mathbf{N}^{\alpha}$.
Hence, it is only possible in this case to introduce a definition of a partial
meshless transmissibility between particles $\displaystyle\vect{r}_{J}$ and
$\displaystyle\vect{r}_{I}$ as follows:
\begin{equation}
	\label{ChapIV2_16}
	\begin{array}{lcl}
		\vspace{0.3cm}
		\displaystyle T^{P}\left(\vect{r}_{J}, \vect{r}_{I}\right)= T^{P}_{JI} =
		\frac{\Gamma^{-1}_{\beta\beta}}{n}\times V_{\vect{r}_{J}}\\
		\displaystyle\left\lbrace
		\frac{\left(\vect{r}_{J}-\vect{r}_{I}\right)\cdot\left(m_{J}+m_{I}
		\right)\cdot
		\overline{\nabla W}
		\left(\vect{r}_{J}-\vect{r}_{I}, \tilde{h}_{IJ}\right)}{\left\Vert\vect{r}_{J}-\vect{r}_{I}\right\Vert^{2}}
		-m_{I}\overline{\nabla W}\left(\vect{r}_{J}-\vect{r}_{I}, \tilde{h}_{IJ}\right)
		\mathbf{N}^{\alpha}
		\right\rbrace.
	\end{array}
\end{equation}
It is important to note that transmissibilities $T^{P}_{JI}$ and $\widetilde{T}^{P}_{JI}$
have different physical units. Furthermore, it raises the question wherever the
scheme (\ref{Schwaiger1})--(\ref{Schwaiger5}) is monotone.

Using Taylor series expansions about a point 
$\vect{r}_{I}$, i.e. relations (\ref{Brookshaw2}), (\ref{Brookshaw3}), 
the expression \eqref{LV_SPH_Schwaiger4} can be written keeing higher 
order error terms as:
\begin{equation}
	\label{ChapIV2_18}
	\begin{array}{lcl}
		\displaystyle
		\vspace{0.3cm}
		\left[\left\langle
		\nabla_{\alpha}\left(m\left(\vect{r}_{I}\right)\mathbf{u}
		\left(\vect{r}_{I}\right)\right)
		\right\rangle-\mathbf{u}\left(\vect{r}_{I}\right)
		\left\langle \nabla_{\alpha}m\left(\vect{r}_{I}\right)\right\rangle+
		m\left(\vect{r}_{I}\right)\left\langle\nabla_{\alpha}\mathbf{u}
		\left(\vect{r}_{I}\right)\right\rangle
		\right] = \\
		\vspace{0.3cm}
		\displaystyle
		2m\left(\vect{r}_{I}\right)
		\left(\mathbf{u}_{,\alpha}\left(\vect{r}_{I}\right)+
		\frac{1}{2}\mathbf{u}_{,\omega\gamma}\left(\vect{r}_{I}\right)
		\sum\limits_{\Omega_{\vect{r}_{I},\tilde{h}_{I}}}V_{\vect{r}_{J}}
		\left[\vect{r}^{\omega}_{J}-\vect{r}^{\omega}_{I}\right]
		\left[\vect{r}^{\gamma}_{J}-\vect{r}^{\gamma}_{I}\right]
		\overline{\nabla^{*}_{\alpha} W}
		\left(\vect{r}_{J}-\vect{r}_{I}, \tilde{h}_{IJ}\right)\right)+\\
		+\mathcal{O}\left(\tilde{h}_{I}^{3}\right).
	\end{array}
\end{equation}

Hence, there is an additional term that has not been taken into account in
\eqref{Schwaiger1}--\eqref{Schwaiger5}. The following section describes
an alternative numerical scheme for the heterogeneous Laplace operator. 
Some initial attemps were also made in \cite{LukyanovVuik2017}.

\subsection{ New scheme}
\label{NewScheme}

The correction terms to the Brookshaw \cite{Brookshaw1985}
and Schwaiger \cite{Schwaiger2008} formulations which improve the accuracy 
of the Laplacian operator near boundaries can be done as follows:
\begin{equation}
	\label{NewScheme1}
	\begin{array}{lcl}
		\displaystyle
		\vspace{0.3cm}
		-\frac{n}{\bar\Gamma^{-1}_{\beta\beta}}
		\langle\nabla\left(\mathbf{M}\left(\vect{r}_{I}\right)\nabla \mathbf{u}
		\left(\vect{r}_{I}\right)\right)
		\rangle = \\
		\vspace{0.3cm}
		\displaystyle\left\lbrace
		\sum_{\Omega_{\vect{r}_{I},\tilde{h}_{I}}}V_{\vect{r}_{J}}
		\left[\mathbf{u}\left(\vect{r}_{J}\right)-\mathbf{u}\left(\vect{r}_{I}\right)
		\right]
		\frac{\left(\vect{r}_{J}-\vect{r}_{I}\right)\cdot
		\left(\mathbf{M}_{J}+\mathbf{M}_{I}\right)
		\cdot
		\overline{\nabla W}
		\left(\vect{r}_{J}-\vect{r}_{I}, \tilde{h}_{IJ}\right)}{\left\Vert\vect{r}_{J}-
		\vect{r}_{I}\right\Vert^{2}}
		\right\rbrace-\\-
		\displaystyle
		\left\lbrace
		\mathbf{N}\cdot
		\left(\sum\limits_{\Omega_{\vect{r}_{I},\tilde{h}_{I}}}V_{\vect{r}_{J}}
		\cdot (\mathbf{M}_{J}+\mathbf{M}_{I})\cdot
		\left[\mathbf{u}\left(\vect{r}_{J}\right)-\mathbf{u}
		\left(\vect{r}_{I}\right)\right]
		\overline{\nabla^{*}W}\left(\vect{r}_{J}-\vect{r}_{I}, \tilde{h}_{IJ}\right)
		\right)
		\right\rbrace,
	\end{array}
\end{equation}
where $n=1,2,3$ is the spatial dimension and tensor $\bar\Gamma_{\alpha\beta}$
is defined by
\begin{equation}
	\label{NewScheme2}
	\bar\Gamma_{\alpha\beta}\left(\vect{r}_{I}\right) =
	\left\{
	\begin{array}{lcl}
		\displaystyle
		\vspace{0.3cm}
		\Gamma^{*}_{\alpha\beta}\left(\vect{r}_{I}\right), \ \ \
		\Gamma^{*}_{\alpha\beta}\left(\vect{r}_{I}\right)\neq 0, \\
		\Gamma_{\alpha\beta}\left(\vect{r}_{I}\right), \ \ \
		\Gamma^{*}_{\alpha\beta}\left(\vect{r}_{I}\right) = 0
	\end{array}
	\right.
\end{equation}
where
\begin{equation}
	\label{NewScheme3}
	\begin{array}{lcl}
		\displaystyle
		\vspace{0.3cm}
		\Gamma^{*}_{\alpha\beta}\left(\vect{r}_{I}\right) =
		\displaystyle
		\sum_{\Omega_{\vect{r}_{I},\tilde{h}_{I}}}V_{\vect{r}_{J}}
		\left[\vect{r}^{\beta}_{J}-\vect{r}^{\beta}_{I}\right]
		\overline{\nabla_{\alpha}W}\left(\vect{r}_{J}-\vect{r}_{I},\tilde{h}_{IJ}\right)-\\
		\displaystyle-
		\mathbf{N}^{\gamma}\sum\limits_{\Omega_{\vect{r}_{I},\tilde{h}_{I}}}V_{\vect{r}_{J}}
		\left[\vect{r}^{\alpha}_{J}-\vect{r}^{\alpha}_{I}\right]
		\left[\vect{r}^{\beta}_{J}-\vect{r}^{\beta}_{I}\right]
		\overline{\nabla^{*}_{\gamma} W}\left(\vect{r}_{J}-\vect{r}_{I}, \tilde{h}_{IJ}\right).
	\end{array}
\end{equation}

Following \eqref{NewScheme1}, we only need to compute the trace of the matrix 
$\Gamma^{*}_{\alpha\beta}\left(\vect{r}_{I}\right)$. Furthermore, it is important to 
note the following remark.  	
     \begin{remark}\label{Remark_Gamma_Structure}
		The following relations can be written:
		\begin{equation}
		\label{GammaRelation1}
			\begin{array}{lcl}
				\displaystyle
				\vspace{0.3cm}
				\Gamma^{*}_{\beta\beta}\left(\vect{r}_{I}\right) =
				\displaystyle
				\sum_{\Omega_{\vect{r}_{I},\tilde{h}_{I}}}V_{\vect{r}_{J}}
	            \left\Vert\vect{r}_{J}-\vect{r}_{I}\right\Vert
                \frac{1}{\tilde{h}_{IJ}}\frac{\overline{dW}}{dz}-\\
				\displaystyle-
				\mathbf{N}^{\gamma}\sum\limits_{\Omega_{\vect{r}_{I},\tilde{h}_{I}}}V_{\vect{r}_{J}}
	            \left\Vert\vect{r}_{J}-\vect{r}_{I}\right\Vert^2						
				\overline{\nabla^{*}_{\gamma} W}\left(\vect{r}_{J}-\vect{r}_{I}, \tilde{h}_{IJ}\right), \
				\frac{\overline{dW}}{dz} \le 0.
			\end{array}
		\end{equation}
		\begin{equation}
			\label{EffectiveGradWN}
			\displaystyle
			\left| \overline{\nabla^{*}_{\alpha}W}
			\mathbf{N}^{\alpha}\left(\vect{r}_{I}\right)\right| \le
			\max_{z\in\Upsilon}\left(\frac{\overline{dW}}{dz}\right)^{2}\cdot \tilde{h}
			\cdot\left\Vert C_{\alpha\beta}\right\Vert\cdot V_{\Omega}\cdot
			\left\Vert
			\frac{\sum\limits_{\vect{r}_{\xi}\in\Omega_{\vect{r}_{I},\tilde{h}_{I}}}
			V_{\xi}\cdot\vect{r}_{\xi}}{V_{\Omega_{\vect{r}_{I},\tilde{h}_{I}}}}-\vect{r}_{I}\right\Vert,
		\end{equation}
		where $\displaystyle z=\left\|\vect{r}_{J}-\vect{r}_{I}\right\|/\tilde{h}_{IJ}$, 
		$\displaystyle\Upsilon = \mbox{supp} \ \frac{\overline{dW}}{dz}$,
		$\left\Vert C_{\alpha\beta}\right\Vert$ is the matrix norm of $C_{\alpha\beta}$, 
		$V_{\Omega_{\vect{r}_{I},\tilde{h}_{I}}}$ is the volume of
		$\Omega_{\vect{r}_{I},\tilde{h}_{I}}$. It follows that if 
		$\displaystyle\left\Vert
		\frac{\sum V_{\xi}\vect{r}_{\xi}}{V_{\Omega_{\vect{r}_{I},\tilde{h}_{I}}}}-\vect{r}_{I}
		\right\Vert \le \tilde{h}_{I}$ then there is a parameter $\tilde{h}_{I}$ such that 
		$\Gamma^{*}_{\beta\beta}\left(\vect{r}_{I}\right) \le 0$.
	\end{remark}
The reason for having the correction factor in the form 
\eqref{NewScheme2}--\eqref{NewScheme3} is that 
$\Gamma^{*}_{\alpha\beta}\left(\vect{r}_{I}\right)=0$
in some cases, where particles have the incomplete Kernel support 
(e.g., at the corners and boundaries of the numerical domain).
For multi-dimensional problems, the correction tensor $\bar\Gamma_{\alpha\beta}
\left(\vect{r}_{I}\right)$ is also a matrix. If the particle $\vect{r}_{I}$ has
entire stencil support (i.e., the domain support for all kernels
$W\left(\vect{r}_{J}-\vect{r}_{I}, \tilde{h}_{IJ}\right)$ is completed 
and symmetric) then 
$\bar\Gamma_{\alpha\beta}\left(\vect{r}_{I}\right)\approx \delta_{\alpha\beta}$. 
The proposed correction matrix deviates less from the unit 
matrix compare to (\ref{Schwaiger5}). 
As a result, the discretization scheme 
(\ref{NewScheme1})--(\ref{NewScheme3}) is at least
$\mathcal{O}\left(h^\omega\right), \ 1\leq\omega \leq 2$ order of 
accuracy in average for any scalar mobility field
$m\left(\vect{r}\right)\geq 0, m\left(\vect{r}\right)\in C^{1}\left(\Omega\right)$, 
$\mathbf{M}^{\alpha\beta}\left(\vect{r}\right) = 
m\left(\vect{r}\right)\delta_{\alpha\beta}$
everywhere within the numerical domain $\Omega\in\mathbb{R}^{n}$ 
sufficiently far away from the boundary $\partial\Omega$. The scheme has the 
two-point flux approximation nature and can be written in the form of 
\eqref{TPFA_FORM}, which can be proved using the arguments above. The 
scheme (\ref{NewScheme1})--(\ref{NewScheme3}) is in line with an alternative 
formulation for continuum mechanics called the peridynamic model \cite{Silling2000,KatiyarFosterOuchiSharma2014}, 
which was proposed several years ago.

All presented schemes in this paper do not require exact expressions for the
gradient (i.e., spatial derivatives) of the mobility field
$\nabla_{\gamma}m\left(\vect{r}\right)$ to keep a higher order of accuracy for
any mobility field. Hence, this scheme can be used with the 
discontinuous (or piecewise continuous) mobility field
$m\left(\vect{r}\right)\in L_{2}\left(\Omega\right)$. It is important 
to note that Brookshaw \cite{Brookshaw1985} and Schwaiger \cite{Schwaiger2008} 
schemes can also be written for the diagonal mobility matrix 
$\mathbf{M}^{\alpha\beta}\left(\vect{r}\right)$ by substituting 
$\mathbf{M}^{\alpha\beta}$ into \eqref{Brookshaw1} and
\eqref{Schwaiger1} instead of $m\left(\vect{r}\right)$ and performing 
summation by repeating indices.

\section{Approximation,  Stability, and Monotonicity }
\label{ApproximationStabilityMonotonicity}

The approximation, stability and monotonicity are important properties 
of numerical schemes which provide and quantify the confidence
in the numerical modeling and results from corresponding 
simulations. Therefore, in order to be confident that the 
proposed numerical schemes provide the adequate accuracy 
of the elliptic operator \eqref{LV_SPH_Lap}, several 
numerical analyses to identify the order of approximation, 
stability and monotonicity have been performed. It is important to
recall that all numerical schemes are characterized by two 
length scales: $\tilde{h}=f\cdot h_{p}$ is the radius of 
$\Omega_{\vect{r}, \tilde{h}}=\mbox{supp}W$, and $h_{p}$ is the
inter-particle distance. Hence, while investigating the approximation of the
meshless discretization scheme, it is important to distinguish two cases: (a) the
neighborhood number of particles is fixed $f=\mbox{const}$ with varying the
inter-particle distance, (b) the inter-particle distance is fixed
$h_{p}=\mbox{const}$ with varying the neighborhood number of particles.
The first case will be analyzed by looking at the error defined by
\begin{equation}\label{ErrorDefinition}
\left\Vert E \right\Vert^{L_{2}}_{\Omega} =
\left[
\frac{1}{\sum\limits_{J}V_{\vect{r}_{J}}}
\sum\limits_{J}V_{\vect{r}_{J}}
\left(
\mat{L}
\left[\mathbf{u}\left(\vect{r}_{J}\right)\right] -
\langle\mat{L}\rangle
\left[\mathbf{u}\left(\vect{r}_{J}\right)\right]
\right)^{2}
\right]^{1/2},
\end{equation}
where $\mat{L}\left[\cdot\right]$ is the analytical Laplacian at the 
particle $\vect{r}_{J}$ and $\langle\mat{L}\rangle\left[\cdot\right]$ is the 
approximation of the Laplacian at the particle $\vect{r}_{J}$, 
$\left\Vert E \right\Vert^{L_{2}}_{\Omega}$ is the averaged error 
over the entire domain. In the following paragraph, the numerical 
analysis is performed for various functions, particle 
distributions, and media properties.

Following the work \cite{Schwaiger2008}, the ability of the discretization 
to reproduce the Laplacian was tested for several functions in 
$\mathbb{R}^{n}, \ n=1,2,3$:
\begin{equation}\label{LV_SPH_ConsTest1}
	\displaystyle\vspace{0.2cm}
	(a) \ \vect{u}^{s}\left(\vect{x}\right) = \sum\limits^{n}_{i=1} x^{s}_{i}, \ \ \
	\displaystyle
	(b) \ \vect{u}^{s}_{m}\left(\vect{x}\right) =
	\prod\limits^{n}_{i=1} x^{m_{i}}_{i}, \
	\left|m\right|=s
\end{equation}
where
$\displaystyle m=\left(m_{1},\ldots,m_{n}\right), \ \forall i \ : m_{i}\geq 0$ is
the $n$-dimensional multi-index with the property
$\left|m\right| = \sum\limits^{n}_{i=1} m_{i}$. The reason for selecting these
testing polynomials is as follows. Since the functional space
$L_{p}\left(\Omega\right)$ is separable, the above polynomials form the
everywhere dense subset of the $L_{p}\left(\Omega\right)$. Hence, any function
$\mathbf{u}\in L_{p}$ can be approximated in $L_{p}$ using linear combination of
the above polynomials leading to the following relations:
\begin{equation}\label{LV_SPH_ConsTest2}
\displaystyle\mathbf{u}(\vect{r}) \approx \sum\limits^{\infty}_{k=1}
\sum\limits_{\left|m\right|\leq k}a_{m}\left(\prod\limits^{n}_{i=1}
x^{m_{i}}_{i}\right),
\mat{M}\left(\vect{r}\right) \approx \sum\limits^{\infty}_{k=1}
\sum\limits_{\left|m\right|\leq k}b_{m}\left(\prod\limits^{n}_{i=1}
x^{m_{i}}_{i}\right).
\end{equation}
The approximation error produced by the discretization
schemes and considered in this paper for the polynomials 
$\vect{u}^{s}\left(\vect{x}\right)$ and $\vect{u}^{s}_{m}\left(\vect{x}\right)$ 
gives the information about the error growth for the 
arbitrary function $\mathbf{u}(\vect{r})$. In each test, the homogenous 
and heterogeneous particle distribution varying the smoothing length 
$\tilde{h}$ and inter-particle distance $h_{p}$ are used to study the 
approximation properties of the proposed discretization 
scheme.

\subsection{Isotropic Homogeneous Media}
\label{IsotropicHomogeneousMedia}

Note that in the limit of the homogeneous isotropic 
media (i.e., $\mathbf{M}^{\alpha\beta}\left(\vect{r}\right) = 
m\left(\vect{r}\right)\delta_{\alpha\beta}$, $m\left(\vect{r}\right)=1$) 
and without the source term, the aforementioned 
operator $\mathbf{L}\left(\mathbf{u}\right)$ in \eqref{LV_SPH_Lap} 
reduces to the conventional Laplacian operator (i.e., 
$\mathbf{L}\left(\mathbf{u}\right)\equiv \nabla^{2}\mathbf{u}$).
The domain for the patch test is a unit square similar 
to \cite{Schwaiger2008} for $n=2$:
$$
\displaystyle \Omega =
\left\{\vect{r}=\left\{x_{i}\right\}\in \mathbb{R}^{n} \left|\right.
\ \ \left|x_{i}-2.5\right| \leqslant \frac{1}{2} \ \ \forall i \right\}
$$
with $N=21$ particles in each direction characterized 
by two length spacings: $\tilde{h}=f\cdot h_{p}$, 
$h_{p}=0.05$, $f=1.2$.

\begin{figure}[!htbp]
	\centering
	\hspace{-4.0cm}
	\mbox{
	\begin{minipage}{1.7in}
	     \includegraphics[width=6cm, height=6cm]{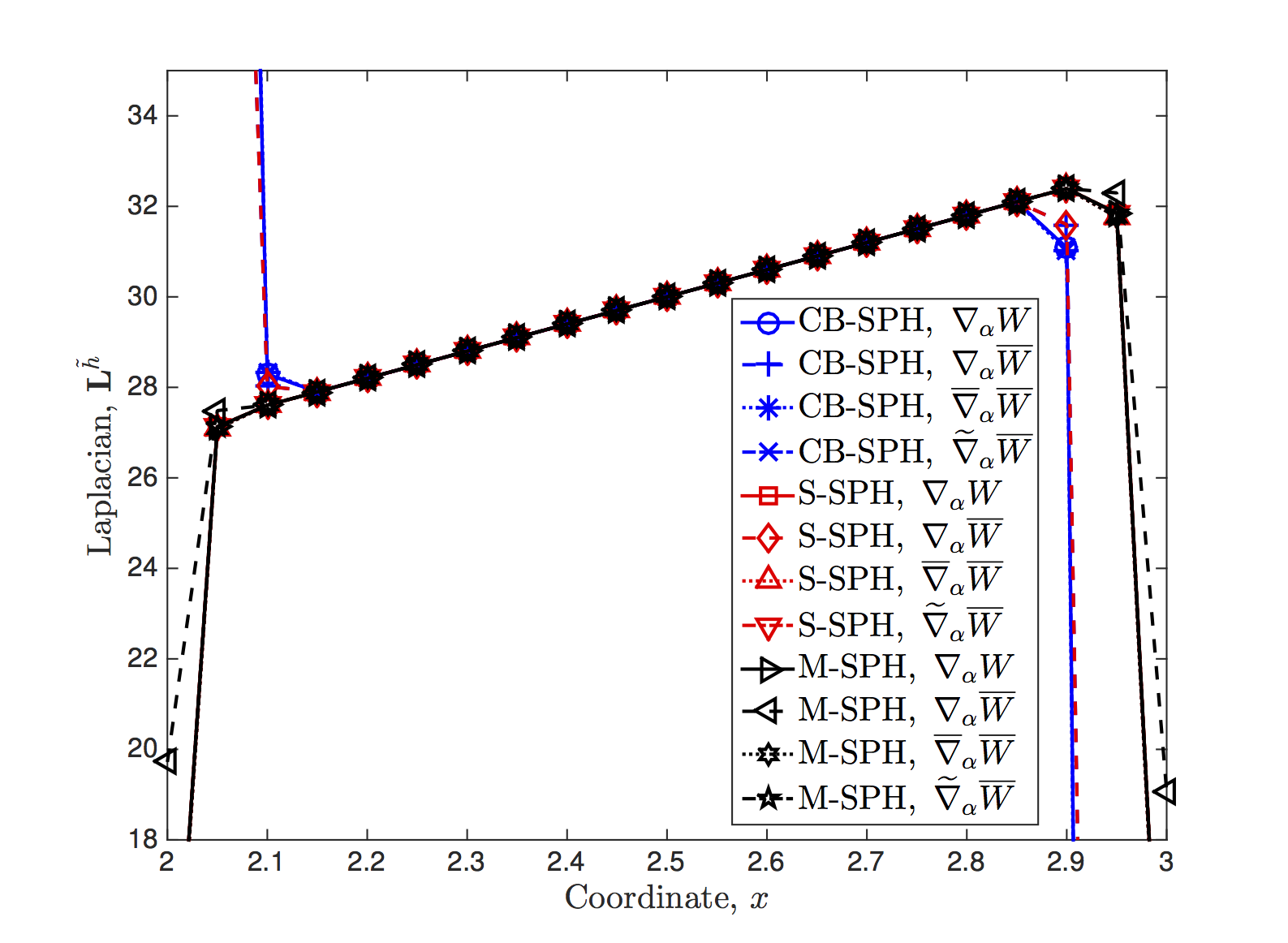}
	\end{minipage}
	\qquad\qquad\qquad
	\begin{minipage}{1.7in}
	    \includegraphics[width=6cm, height=6cm]{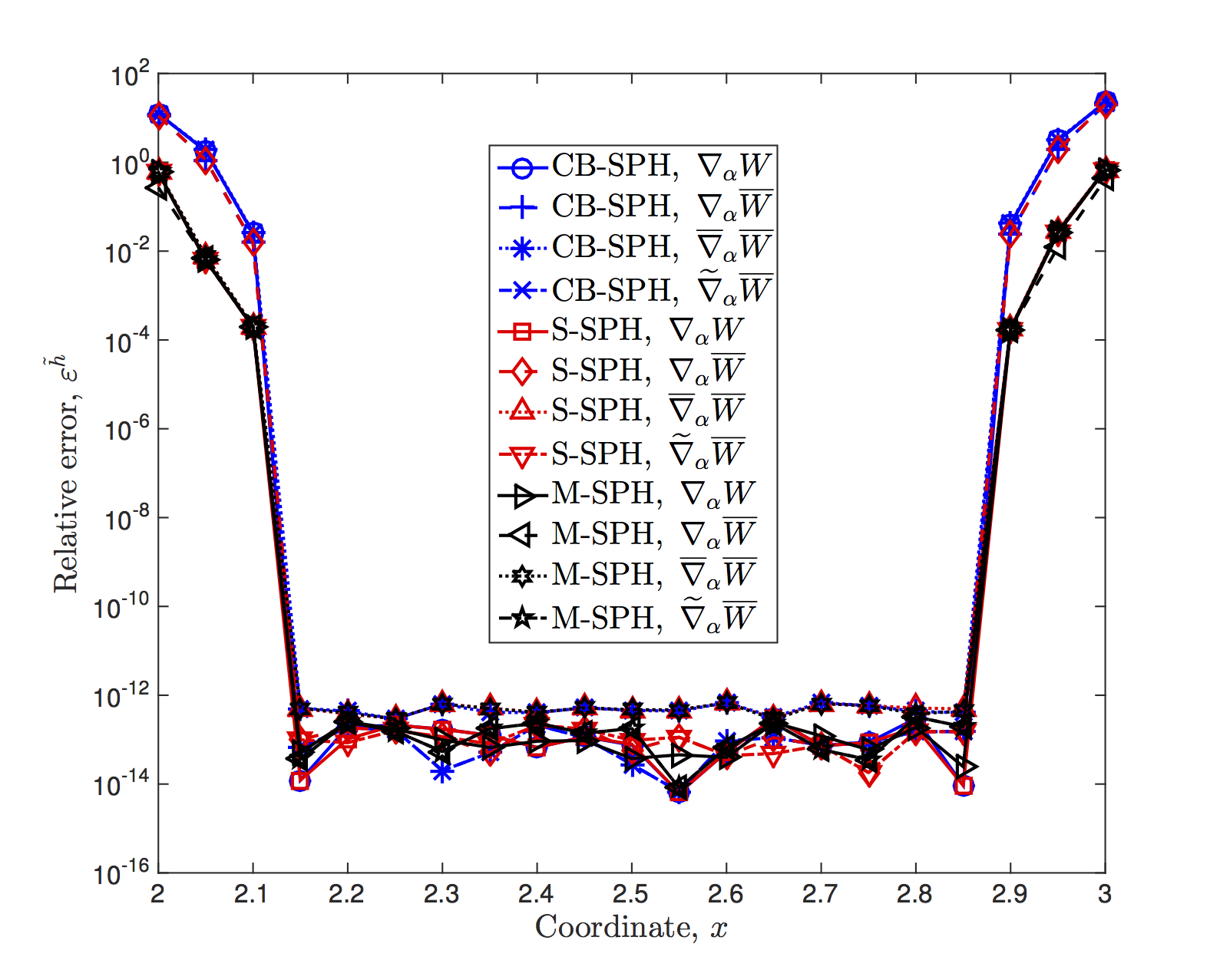}
	\end{minipage}
     }
	\hspace{-3.0cm}     
	\caption{
		Cross-section of the test patch at $y=2.5$. Three SPH 
		approximations of $\nabla^2\left(x^3+y^3\right)$ with different kernel 
		gradients are shown. Corrected Brookshaw's scheme 
		(CB-SPH) is given by \eqref{Brookshaw1M} with the 
		correction multiplier, while  Schwaiger's scheme (S-SPH) is given by 
		\eqref{Schwaiger1}--\eqref{Schwaiger5}. 
		New approximation (M-SPH) considered here is the SPH form 
		\eqref{NewScheme1}--\eqref{NewScheme3}. In this case, the Schwaiger's 
		scheme and new scheme have comparable accuracy at the boundaries and 
		are accurate in the interior to the machine precision. 
		Four different options of computing the kernel 
		gradient (i.e., $\nabla_{\gamma}W$,
		$\nabla_{\alpha}\overline{W}$, and corrected kernel gradients (i.e.,
		$\nabla^{*}_{\alpha}W$, $\nabla^{*}_{\alpha}\overline{W}$,
		$\overline{\nabla}^{*}_{\alpha}\overline{W}$, and
		$\widetilde{\nabla}^{*}_{\alpha}\overline{W}$) are
		shown.
	}
	\label{PatchTest1}
\end{figure}

In all of the following tests, the results are displayed 
along the cross-section $y = 2.5$ in $\mathbb{R}^{n}, \ n=2$ (similar to 
\cite{Schwaiger2008} for $n=2$). In each test, the following discretizations 
are compared: corrected Brookshaw's scheme (CB-SPH) \eqref{Brookshaw1M}, 
Schwaiger's scheme (S-SPH) \eqref{Schwaiger1}--\eqref{Schwaiger5}, and new
proposed scheme (M-SPH) \eqref{NewScheme1}--\eqref{NewScheme3}.
Note that the Schwaiger's scheme was tested against 
several schemes published in \cite{ChenBeraunCarney1999},
\cite{ChaniotisPoulikakosKoumoutsakos2002} 
and show better accuracy; hence, schemes in
\cite{ChenBeraunCarney1999}, \cite{ChaniotisPoulikakosKoumoutsakos2002}
are not considered in this paper. The comparison of different 
schemes starts with test functions of the form  
$\mathbf{u}^{s}\left(\vect{x}\right)$ described by \eqref{LV_SPH_ConsTest1}(a) 
in $\mathbb{R}^{n}, \ n=2$. Plots of the Laplacian approximations and the relative 
errors defined by \eqref{ErrorDefinition} for the case $m=3$ are shown in Figure 
\ref{PatchTest1}.

\begin{figure}[!htbp]
	\centering
	\mbox{
	\begin{minipage}{1.7in}%
		\centering
		\hspace{-4.0cm}
		\includegraphics[width=1.4\linewidth]{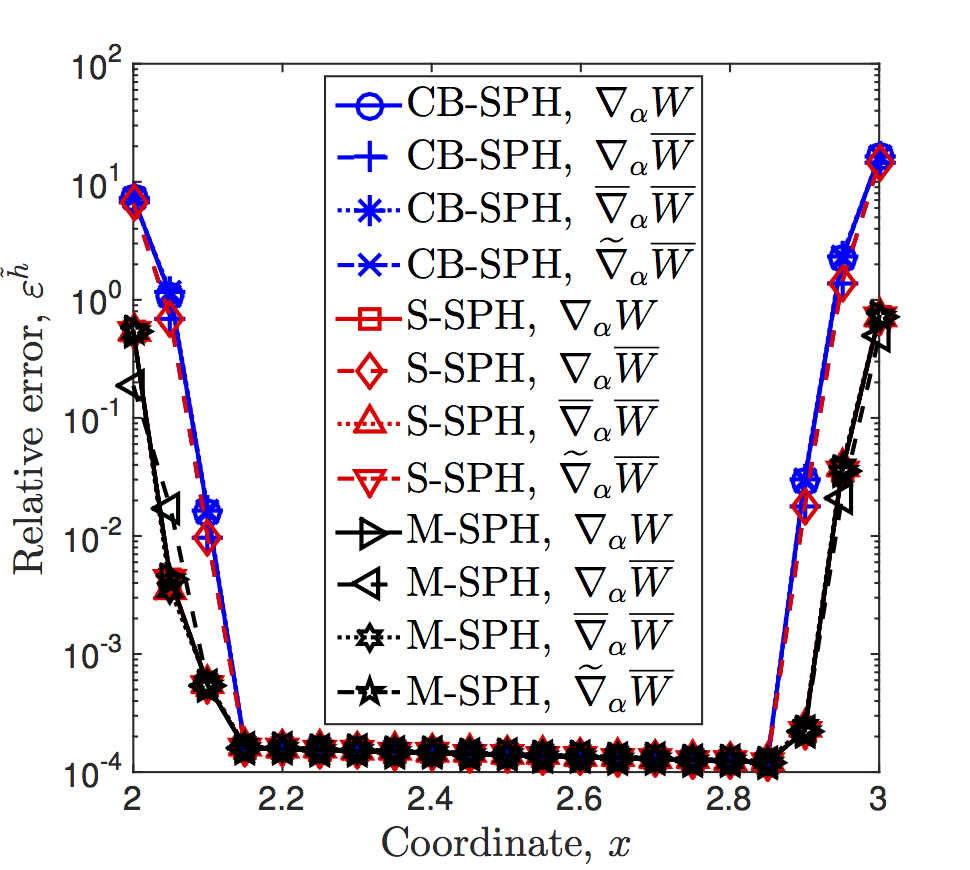}
		\hspace{-4.0cm}
		\caption*{m=4}
		\label{fig:sfig1}
	\end{minipage}%
	\qquad\qquad\qquad\qquad
	\begin{minipage}{1.7in}%
		\centering
		\hspace{-4.0cm}
		\includegraphics[width=1.4\linewidth]{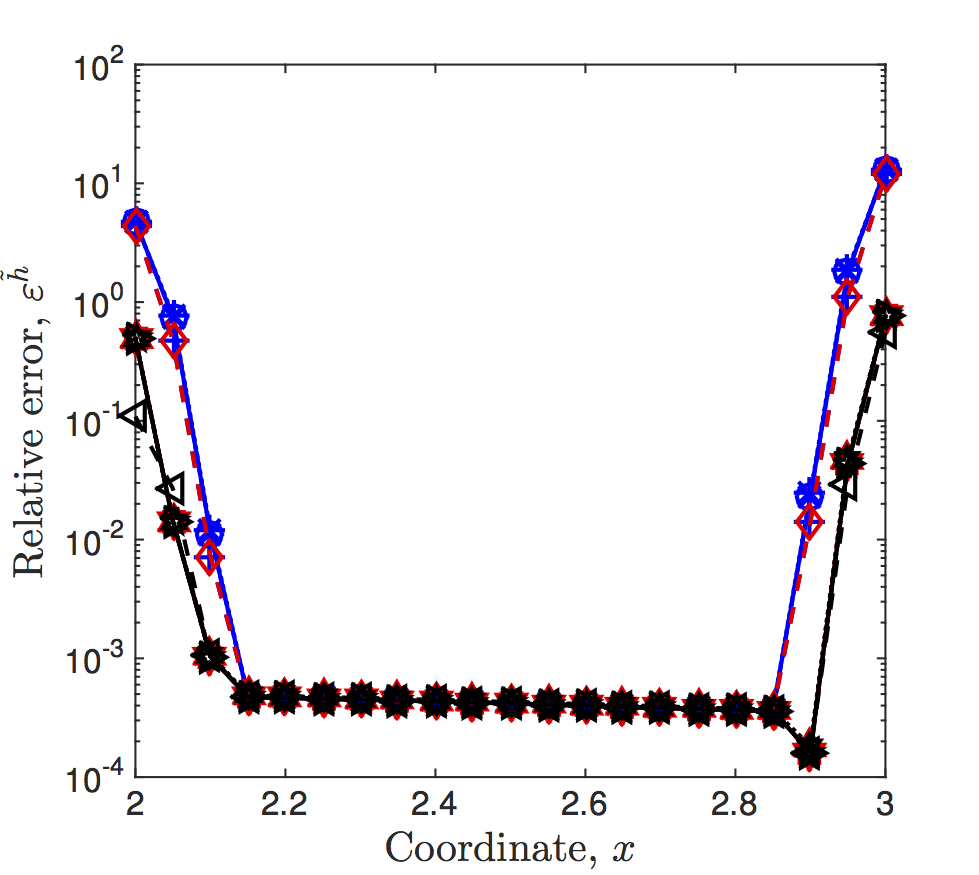}
		\hspace{-4.0cm}
		\caption*{m=5}
		\label{fig:sfig2}
	\end{minipage}%
    }
	\\
	\centering
	\mbox{
	\begin{minipage}{1.7in}%
		\centering
		\hspace{-4.0cm}
		\includegraphics[width=1.4\linewidth]{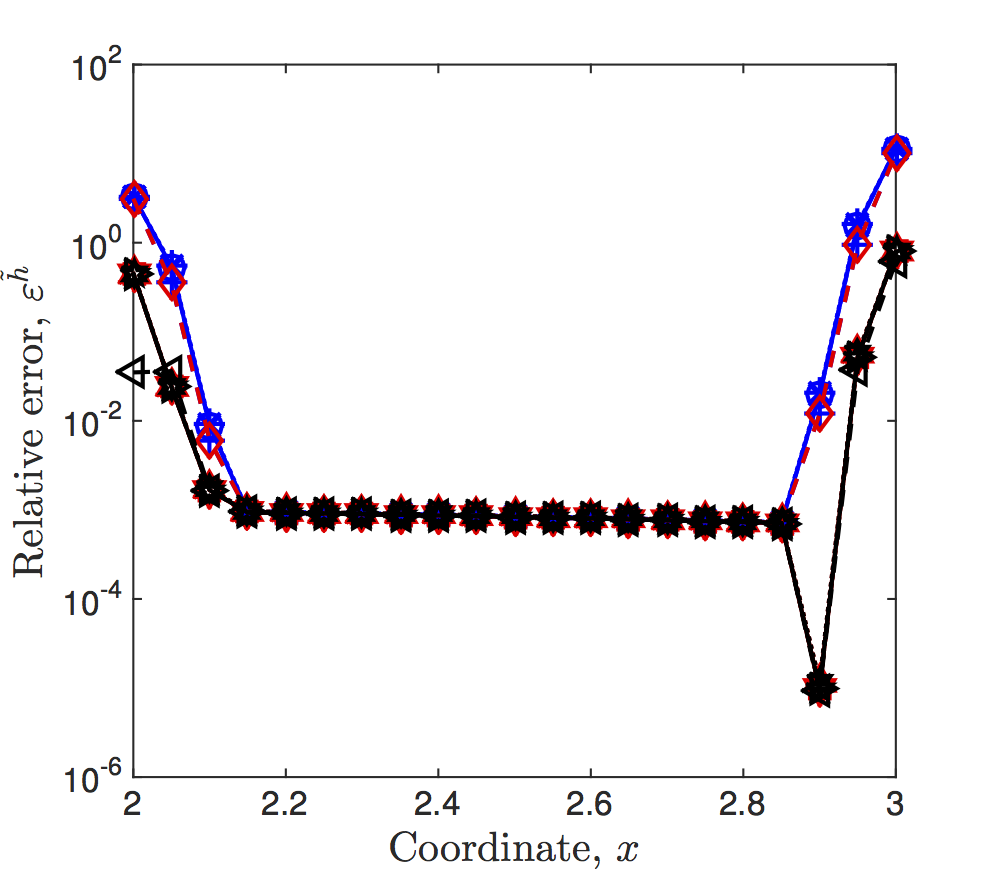}
		\hspace{-4.0cm}
		\caption*{m=6}
		\label{fig:sfig3}
	\end{minipage}%
	\qquad\qquad\qquad\qquad
	\begin{minipage}{1.7in}%
		\centering
		\hspace{-4.0cm}
		\includegraphics[width=1.4\linewidth]{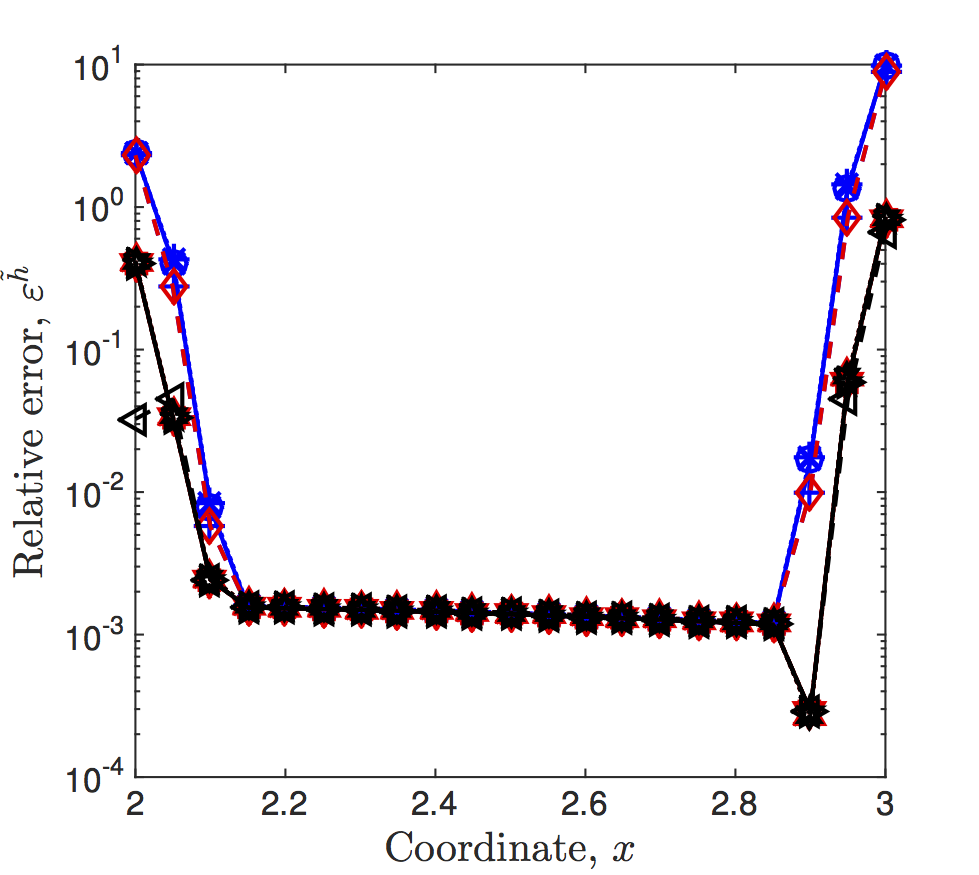}
		\hspace{-4.0cm}
		\caption*{m=7}
		\label{fig:sfig4}
	\end{minipage}%
    }
	\caption{
		The relative errors along $y=2.5$ for each of the discretizations
		used in Figure \ref{PatchTest1} are shown 
		here for the suite of functions $\nabla^2\left(x^m+y^m\right)$ where $m=4,\ldots,7$. 
		Four different options of computing the 
		kernel gradients (i.e., $\nabla_{\gamma}W$, $\nabla_{\alpha}\overline{W}$,
		and corrected kernel gradients (i.e., $\nabla^{*}_{\gamma}W$,
		$\nabla^{*}_{\alpha}\overline{W}$,
		$\overline{\nabla}^{*}_{\alpha}\overline{W}$,
		and  $\widetilde{\nabla}^{*}_{\alpha}\overline{W}$) are
		shown.
	} 
	\label{PatchTest2}
\end{figure}

The new scheme (M-SPH) has the greatest accuracy at the boundary and
is accurate to machine precision $\varepsilon$  in the interior. The error plot
is shown in Figure \ref{PatchTest2} along with test functions with exponents
from $m=4,...,7$. For these functions, the proposed scheme
\eqref{NewScheme1}--\eqref{NewScheme3} is uniformly more
accurate for the various experiments.
The same behavior is observed in 3D where the M-SPH scheme with the
$\nabla_{\alpha}\overline{W}$ is the most accurate scheme. Furthermore, 
the new scheme should be tested for the functions requiring cross-derivatives 
as was reported in \cite{Schwaiger2008}. 
To examine the effect of the  cross-derivative terms, the same suite of tests 
was run with the function $\vect{u}^{s}_{m}\left(\vect{x}\right)$ described by 
\eqref{LV_SPH_ConsTest1}(b) in $\mathbb{R}^{n}, \ n=2,3$.
Relative errors for $(xy)^m$ on a same array as in Figure
\ref{PatchTest1} along $y=2.5$ were computed and again
the proposed scheme \eqref{NewScheme1}--\eqref{NewScheme3} was 
uniformly more accurate for the various experiments.
The behavior of each discretization is similar to that shown in
Figure \ref{PatchTest1}. The proposed new 
scheme (M-SPH) performed nearly as well as the Schwaiger's scheme 
at the boundary. The CB-SPH  and S-SPH forms also perform with 
greater accuracy than all other forms in the interior except in the 
case with the highest exponent.

An additional concern is that although it deviates from the exact 
solution near boundaries, it acquires no off-diagonal terms due to 
the alignment of the array of particles and the boundaries with the
coordinate axes. To test the accuracy of the new approximations 
when there are off-diagonal terms, an array with particles rotated 
$45^{o}$ was used with the test function $(xy)^m$ (similar to 
\cite{Schwaiger2008}). The new proposed scheme formulation 
performs consistently well for lower exponents.

\subsection{Isotropic Heterogeneous Media}
\label{IsotropicHeterogeneousMedia}

\begin{figure}[!htbp]
	\centering
	\hspace{-4.0cm}
	\mbox{
		\begin{minipage}{1.7in}
	     \includegraphics[width=6cm, height=6cm]{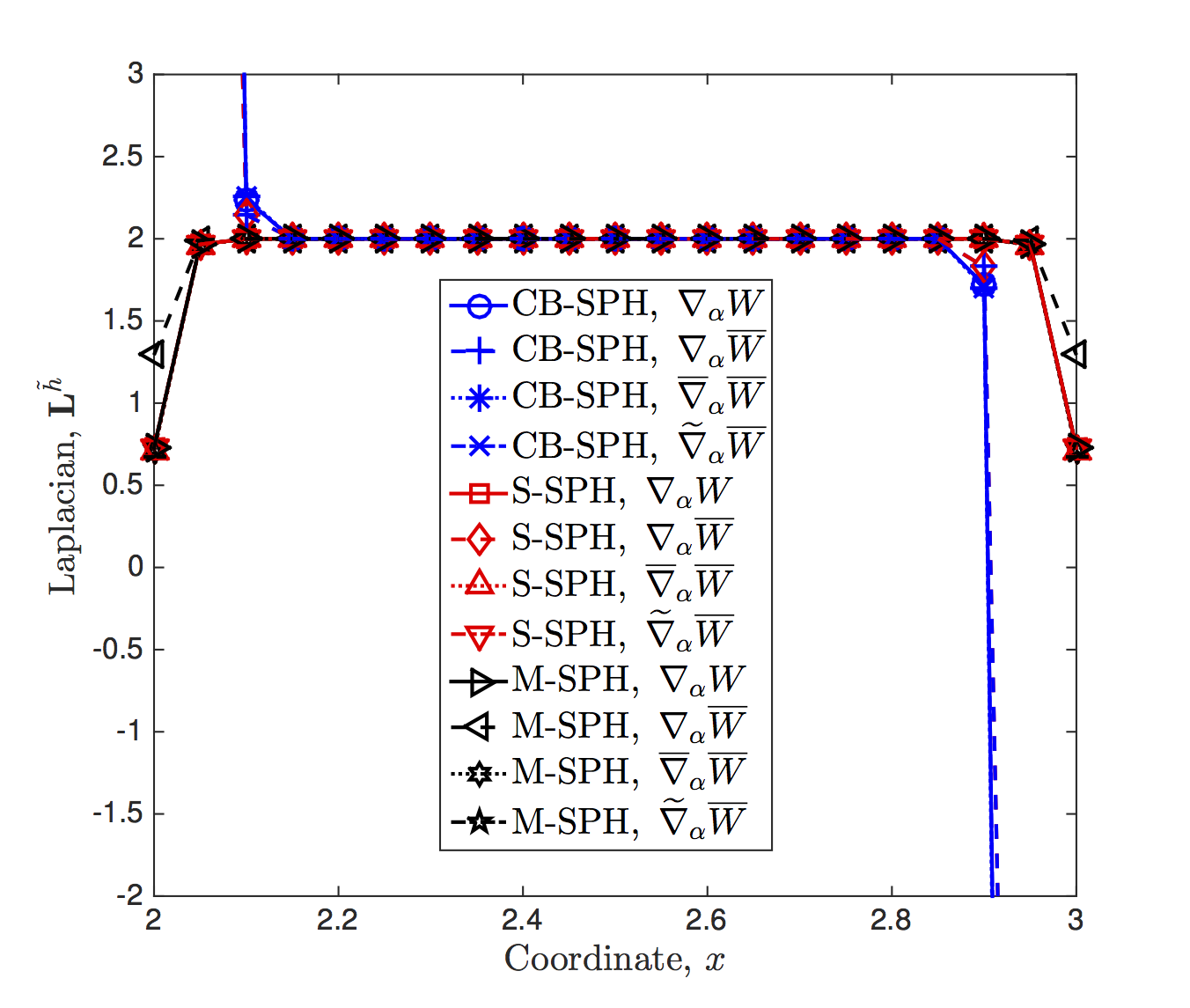}
		\end{minipage}
		\qquad\qquad\qquad
		\begin{minipage}{1.7in}
	     \includegraphics[width=6cm, height=6cm]{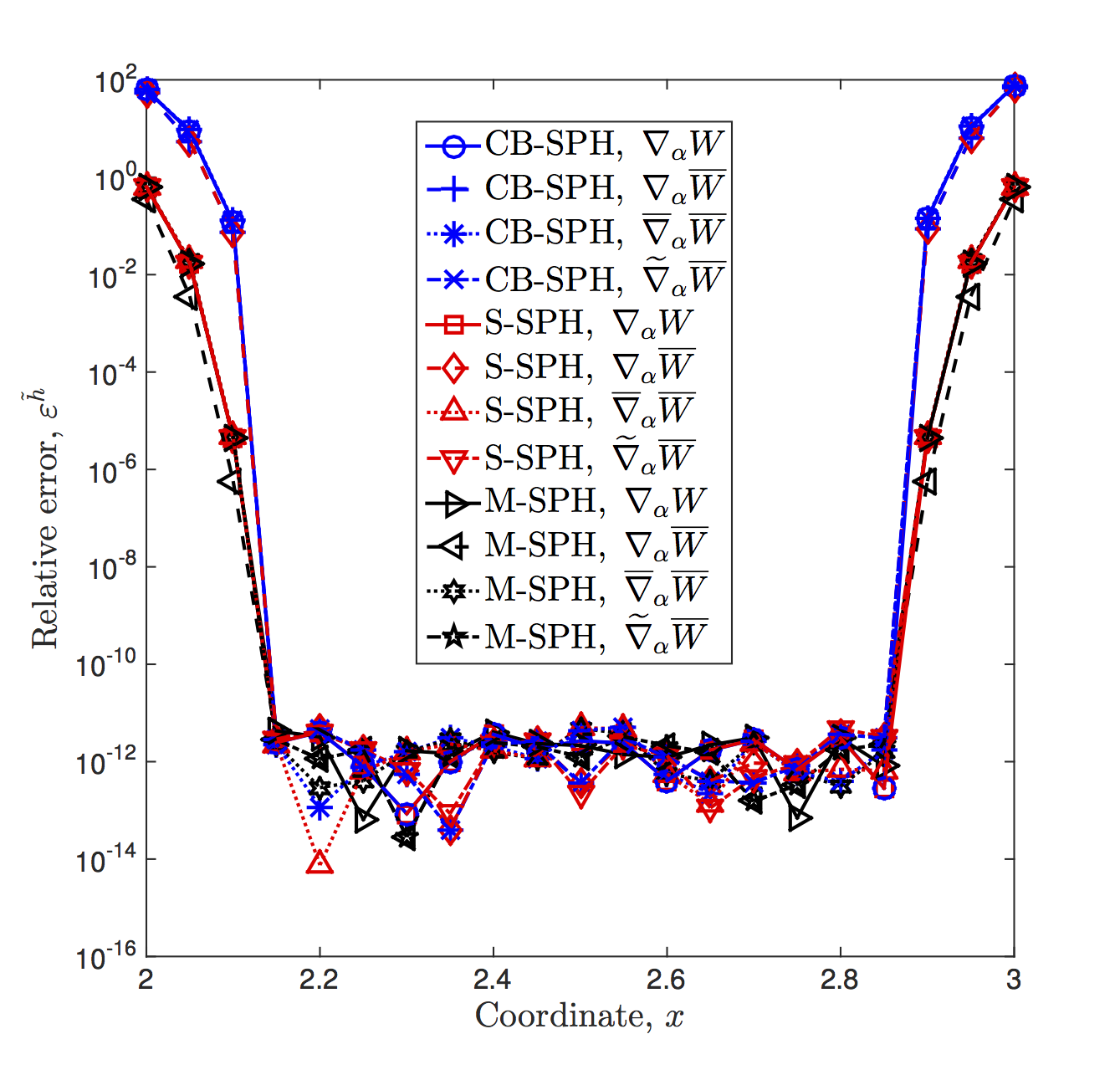}
		\end{minipage}
	}
	\hspace{-3.0cm}	
	\caption{
		Cross-section of the test patch at $y=2.5$. Three SPH 
		approximations of $\nabla(\left(x+y\right)(\nabla\left(x+y\right)))$ with 
		different kernel gradients are shown. The 
		different schemes and kernel gradients are described in the caption of Figure
		\eqref{PatchTest1}.		
	} 
	\label{PatchTestHet1}
\end{figure}

\begin{figure}[!htbp]
	\centering
	\hspace{-1.0cm}
	\mbox{
	\begin{minipage}{1.7in}%
		\centering
		\hspace{-4.0cm}
		\includegraphics[width=1.4\linewidth]{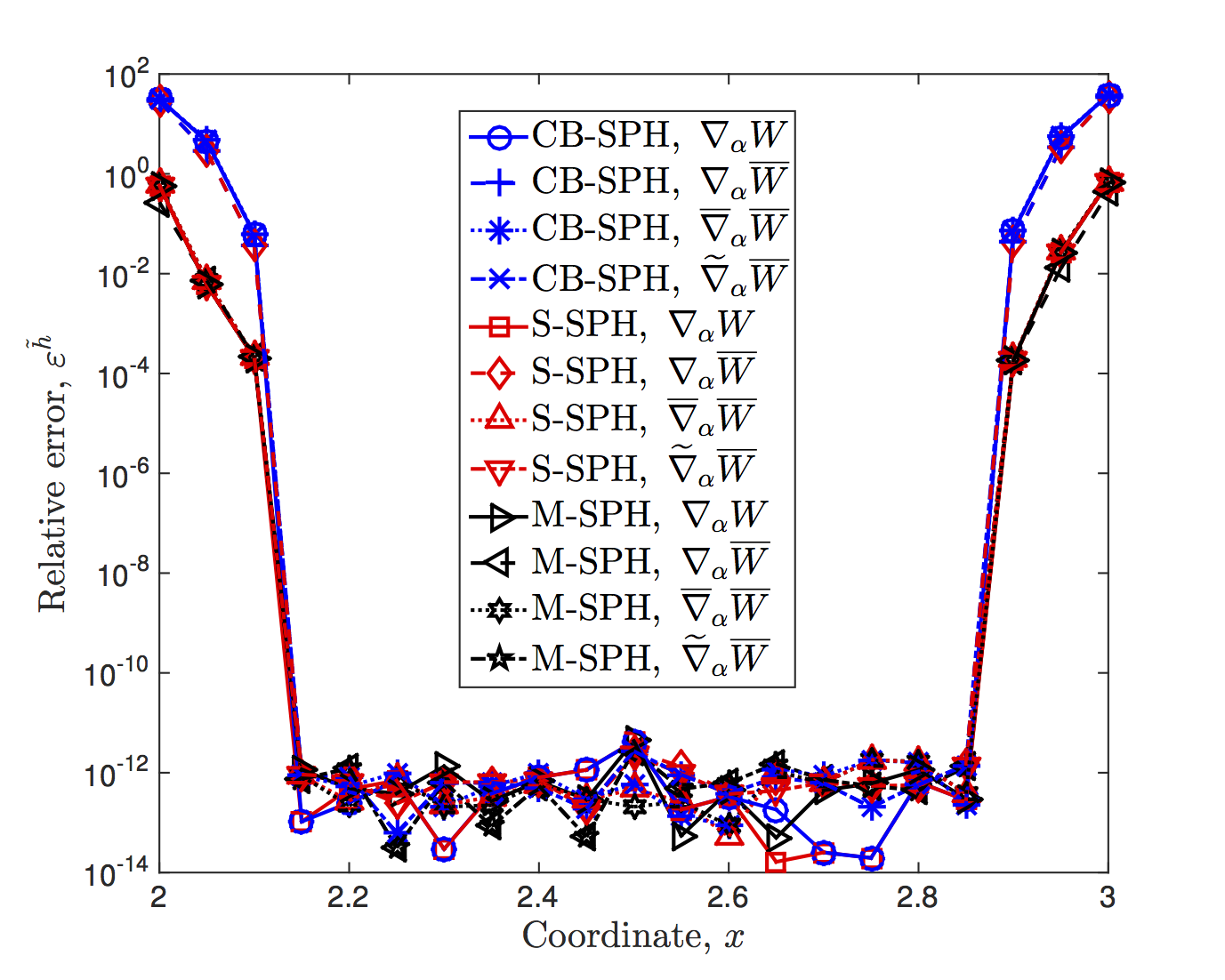}
		\hspace{-4.0cm}
		\caption*{m=2}
		\label{PatchTestHetCase1}
	\end{minipage}%
	\qquad\qquad\qquad\qquad
	\begin{minipage}{1.7in}%
		\centering
		\hspace{-4.0cm}
		\includegraphics[width=1.4\linewidth]{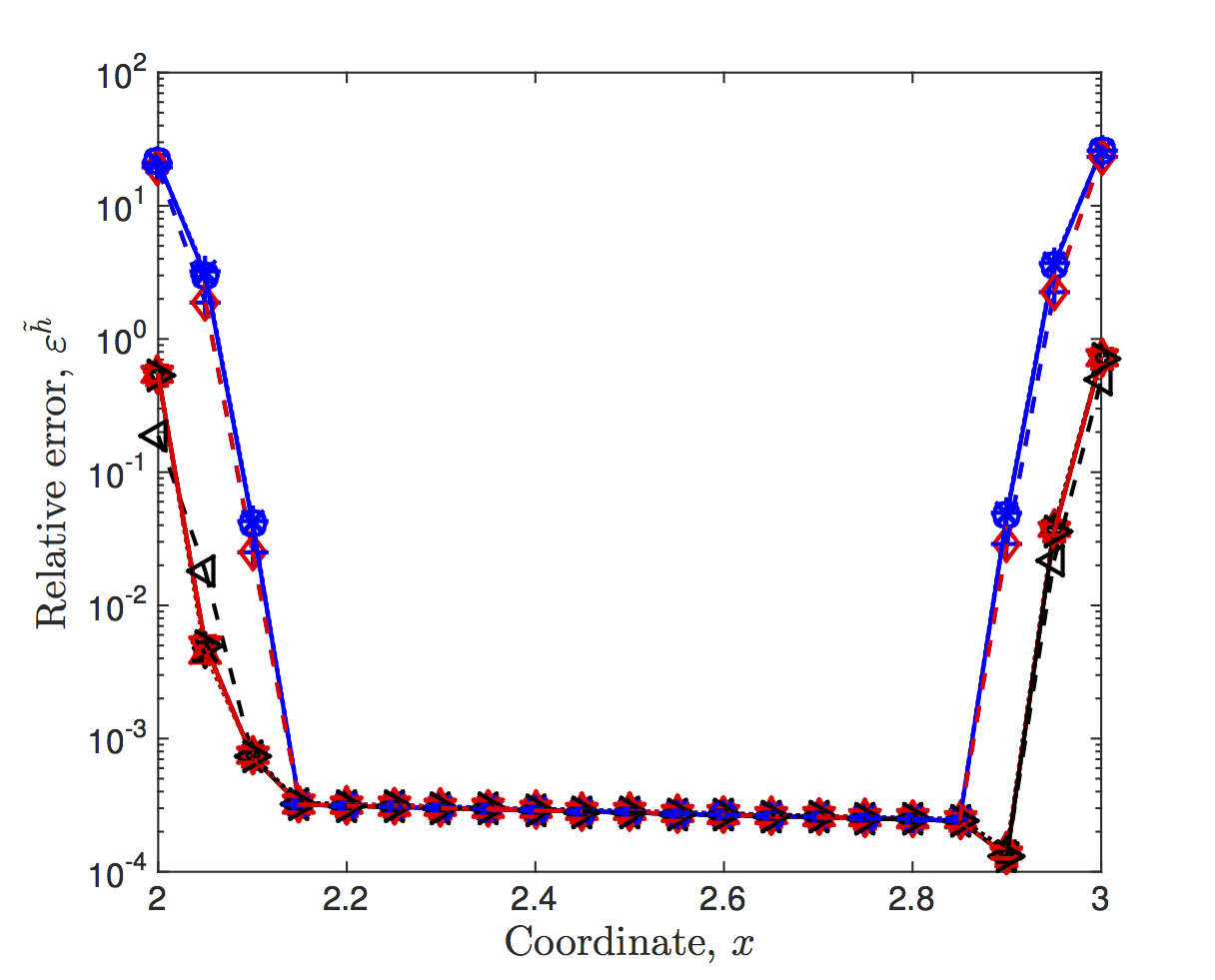}
		\hspace{-4.0cm}
		\caption*{m=3}
		\label{PatchTestHetCase2}
	\end{minipage}%
	}
	\\
	\centering
	\hspace{-1.0cm}
	\mbox{
	\begin{minipage}{1.7in}%
		\centering
		\hspace{-4.0cm}
		\includegraphics[width=1.4\linewidth]{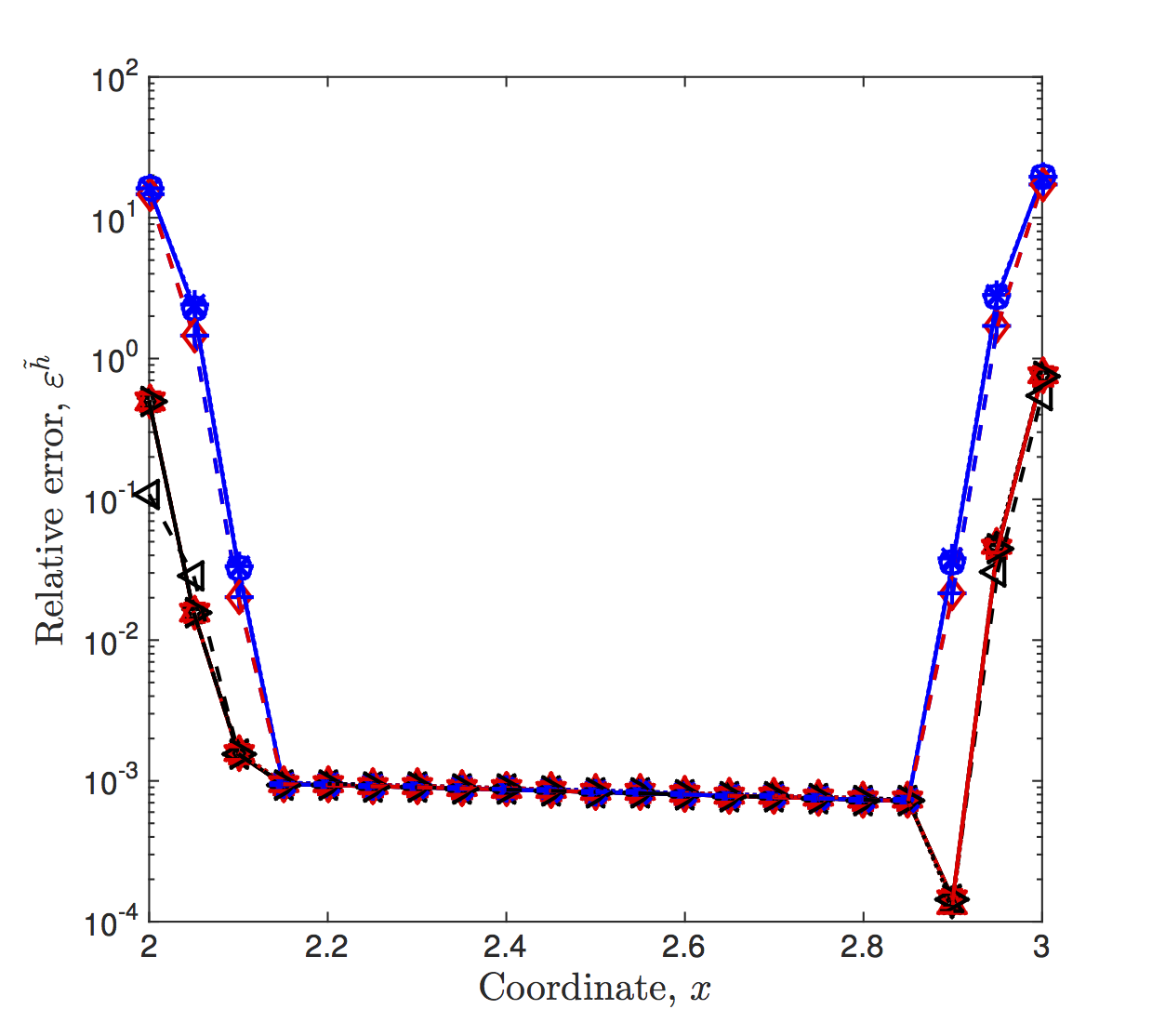}
		\hspace{-4.0cm}
		\caption*{m=4}
		\label{PatchTestHetCase3}
	\end{minipage}%
	\qquad\qquad\qquad\qquad
	\begin{minipage}{1.7in}%
		\centering
		\hspace{-4.0cm}
		\includegraphics[width=1.4\linewidth]{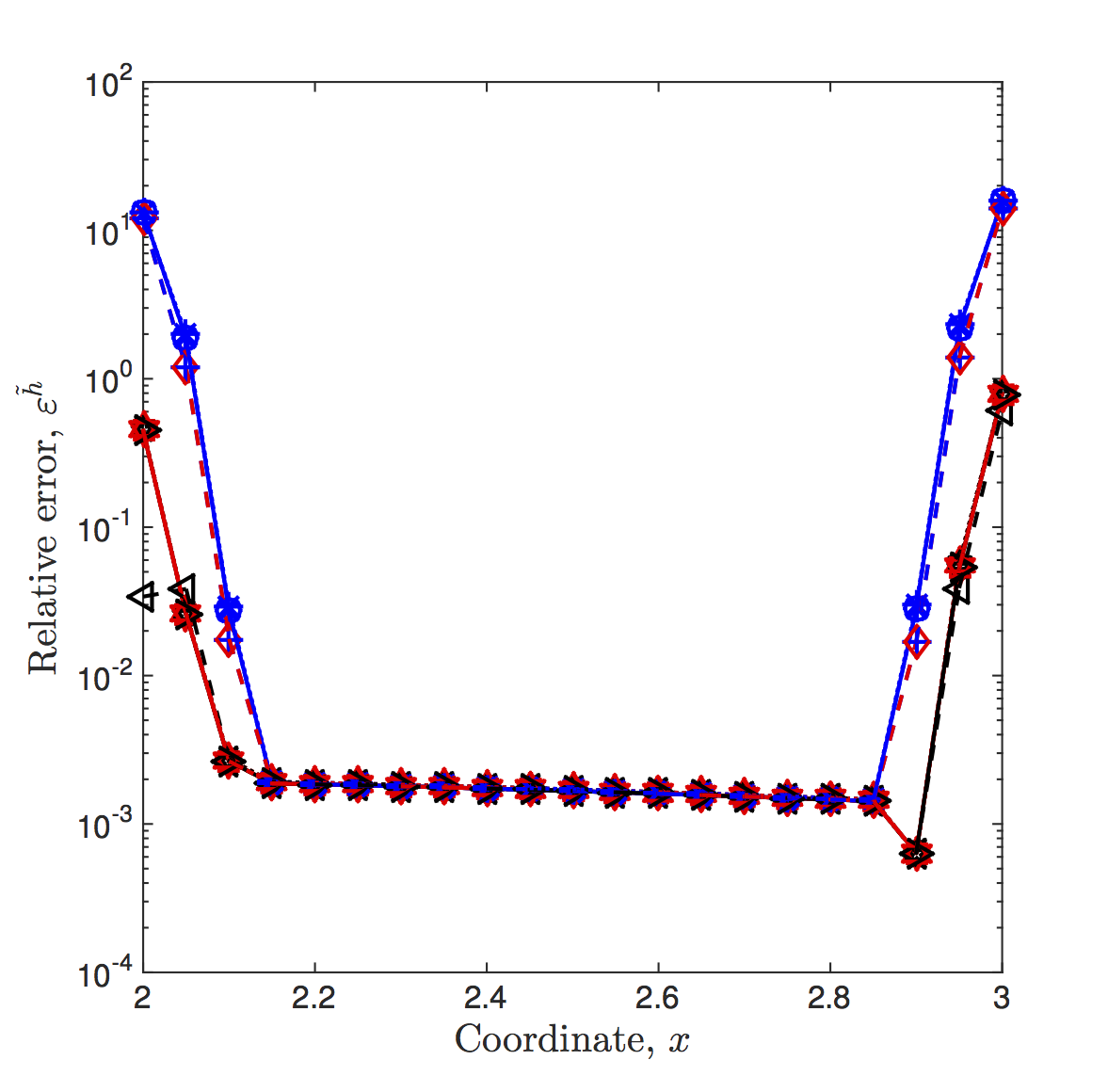}
		\hspace{-4.0cm}
		\caption*{m=5}
		\label{PatchTestHetCase4}
	\end{minipage}%
    }
	\caption{
         The relative errors along $y=2.5$ for each of the discretizations
    	  used in Figure \ref{PatchTest1} are shown here 
    	  only for the suite of functions 
    	  $\displaystyle\nabla\left[\left(x^m+y^m\right)\nabla \left(x+y\right)\right]$ 
    	  where $m=2,\ldots,5$. Four different options of 
	      computing kernel gradients (i.e., $\nabla_{\gamma}W$, 
	      $\nabla_{\alpha}\overline{W}$), and corrected kernel gradients (i.e., 
	      $\nabla^{*}_{\gamma}W$, $\nabla^{*}_{\alpha}\overline{W}$, 
	      $\overline{\nabla}^{*}_{\alpha}\overline{W}$, and 
	      $\widetilde{\nabla}^{*}_{\alpha}\overline{W}$) are shown.} 
	 \label{PatchTestHet2}
\end{figure}

In isotropic heterogeneous media, i.e., $\mathbf{M}^{\alpha\beta}\left(\vect{r}\right) = 
m\left(\vect{r}\right)\delta_{\alpha\beta}$, the Laplace operator 
takes the general form written in \eqref{LV_SPH_Lap}. The numerical domain 
is the same as in the previous section with the same number of particles and
particle length scales: $\tilde{h}=f\cdot h_{p}$, $h_{p}=0.05$, $f=1.2$. Plots 
of the Laplacian approximations and the relative errors defined by 
\eqref{ErrorDefinition} for the case of heterogeneous mobility with $m=1$ 
(see \eqref{LV_SPH_ConsTest2}) are shown in Figure \ref{PatchTestHet1}.
The new scheme (M-SPH) has the best accuracy at the 
boundary and is accurate to machine precision $\varepsilon$  
in the interior. The error plot is shown in Figure \ref{PatchTestHet2} along 
with test functions with exponents from $m=2,...,5$. For these functions, 
the proposed scheme \eqref{NewScheme1}-\eqref{NewScheme3} is uniformly 
more accurate for the different increasing exponents.
The same behavior is observed in 3D where M-SPH scheme with the
$\nabla_{\alpha}\overline{W}$ is the most accurate scheme.

\subsection{von Neumann stability analysis}
\label{NeumannStability}

\begin{figure}[!htbp]
 \mbox{
	\begin{minipage}{1.6in}%
		\centering
		\hspace{-4.0cm}
		\includegraphics[width=4cm, height=4cm]{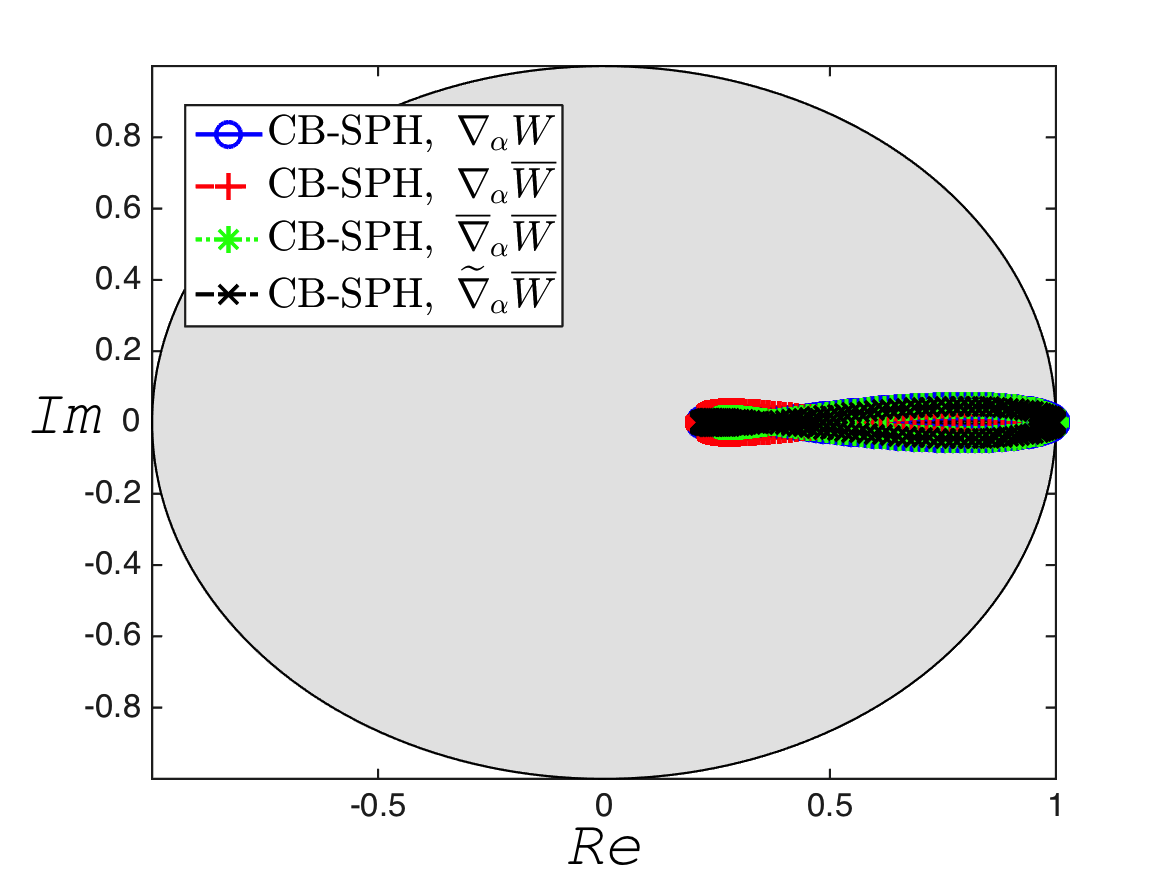}
		\hspace{-4.0cm}
		\caption*{Uniform}
		\label{Stab-CB-RegBox}
	\end{minipage}%
	\begin{minipage}{1.6in}%
		\centering
		\hspace{-4.0cm}
		\includegraphics[width=4cm, height=4cm]{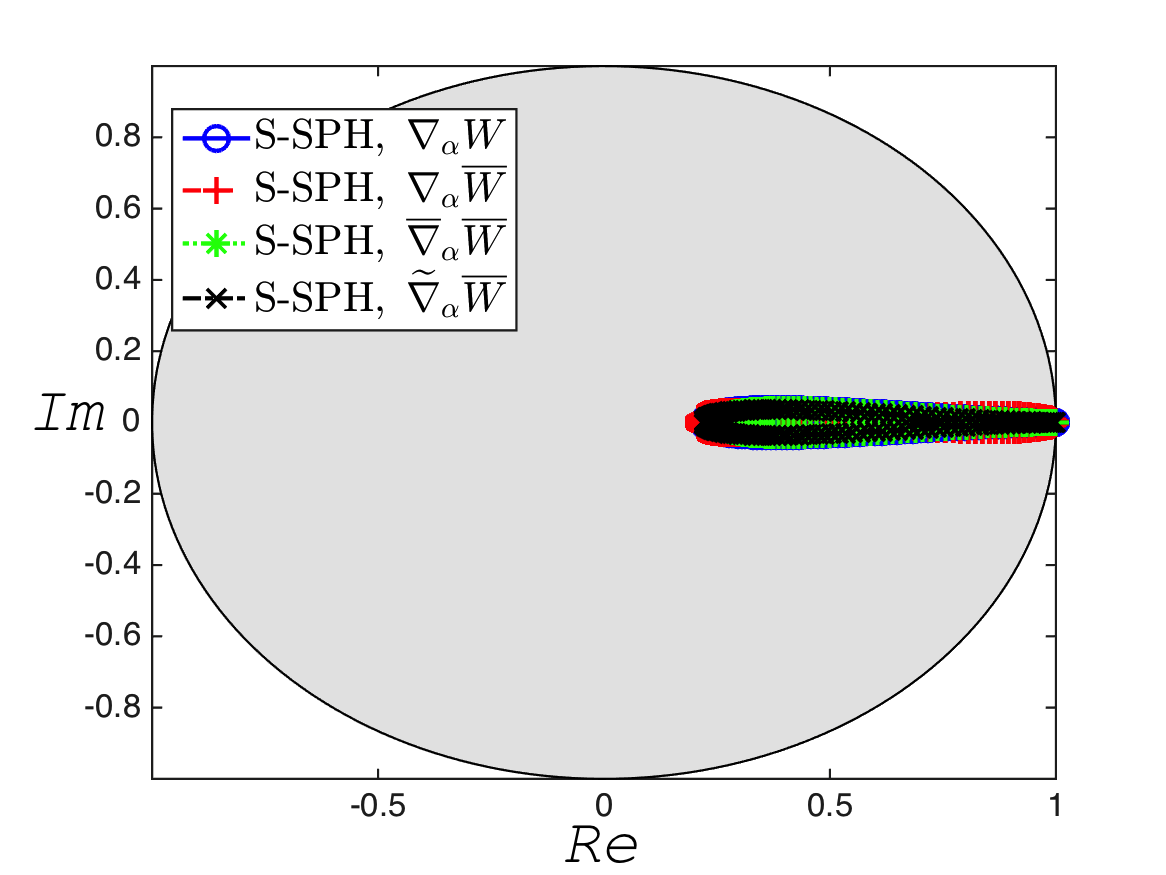}
		\hspace{-4.0cm}
		\caption*{Uniform}
		\label{Stab-S-RegBox}
	\end{minipage}%
	\begin{minipage}{1.6in}%
		\centering
		\hspace{-4.0cm}
		\includegraphics[width=4cm, height=4cm]{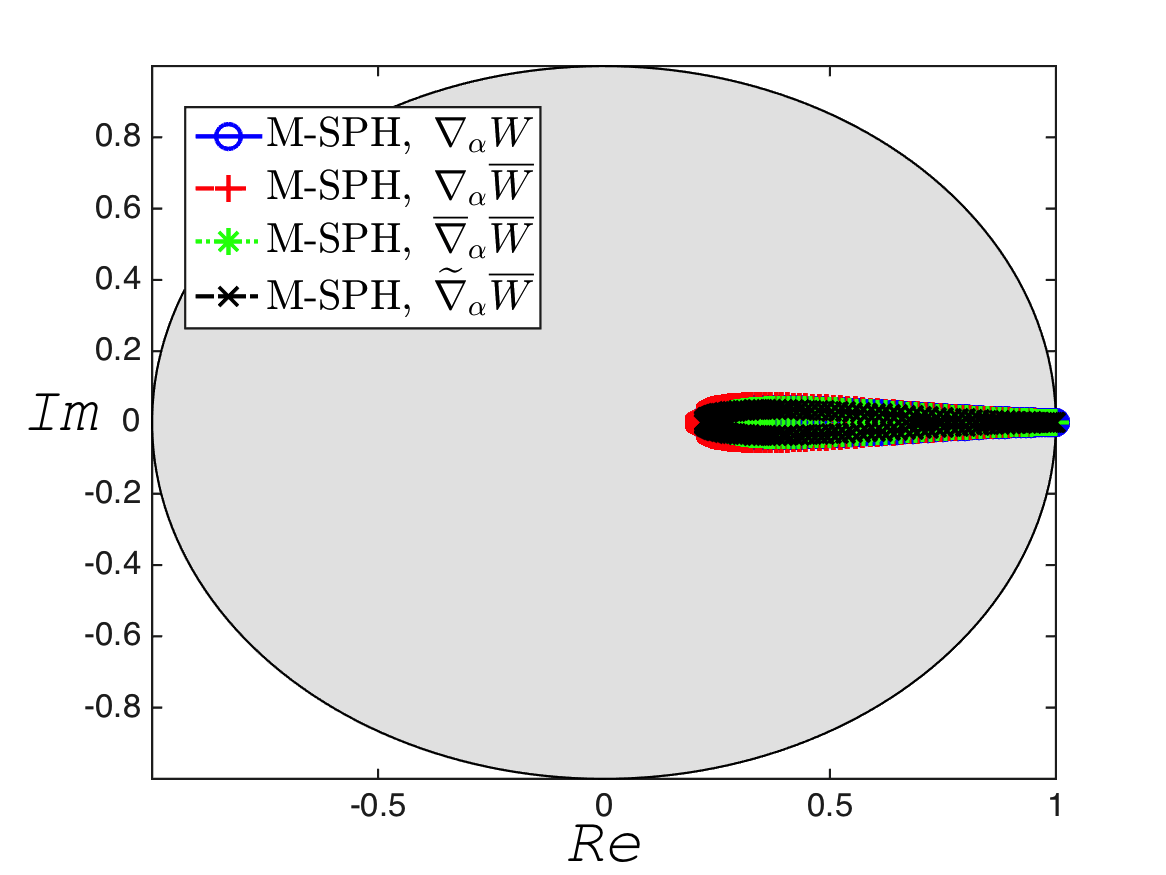}
		\hspace{-4.0cm}
		\caption*{Uniform}
		\label{Stab-M-RegBox}
	\end{minipage}%
    }
	\\
	\mbox{
	\begin{minipage}{1.6in}%
		\centering
		\hspace{-4.0cm}
		\includegraphics[width=4cm, height=4cm]{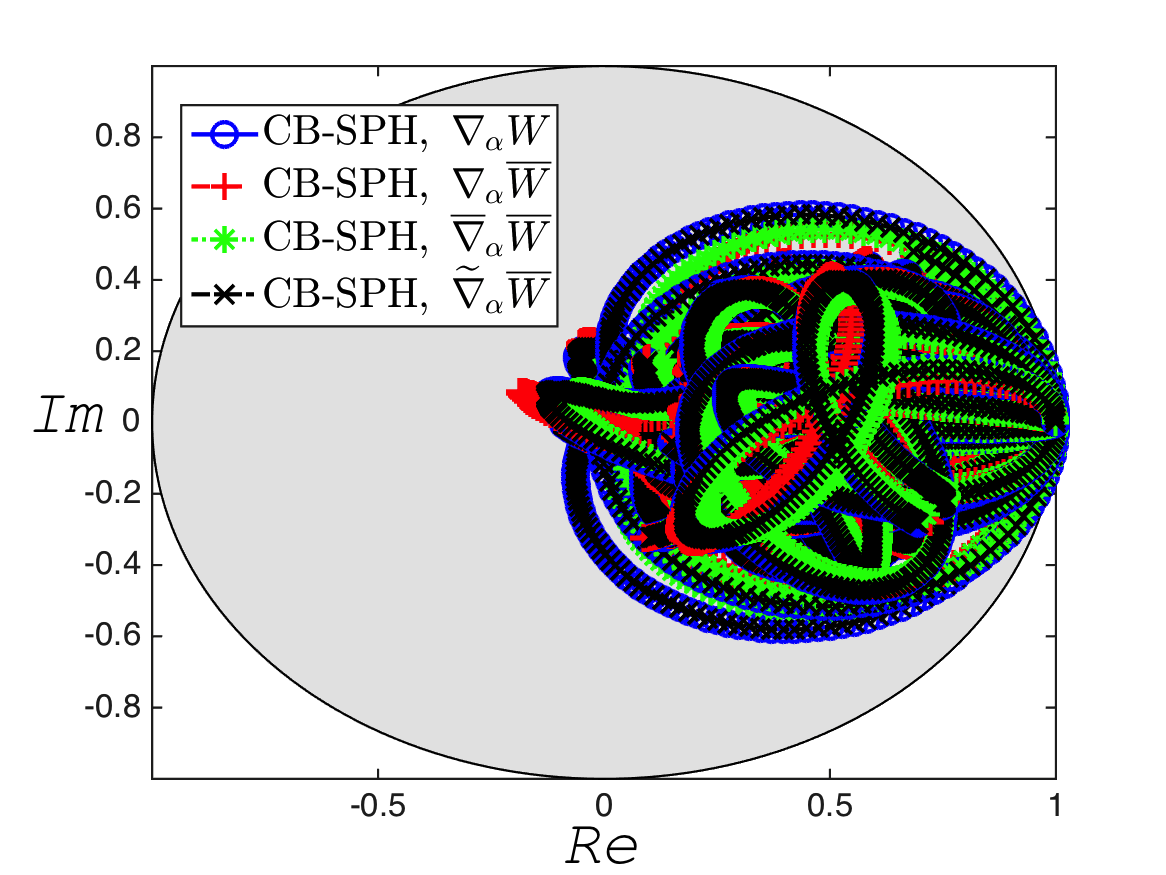}
		\hspace{-4.0cm}
		\caption*{Random}
		\label{Stab-CB-RandBox}
	\end{minipage}%
	\begin{minipage}{1.6in}%
		\centering
		\hspace{-4.0cm}
		\includegraphics[width=4cm, height=4cm]{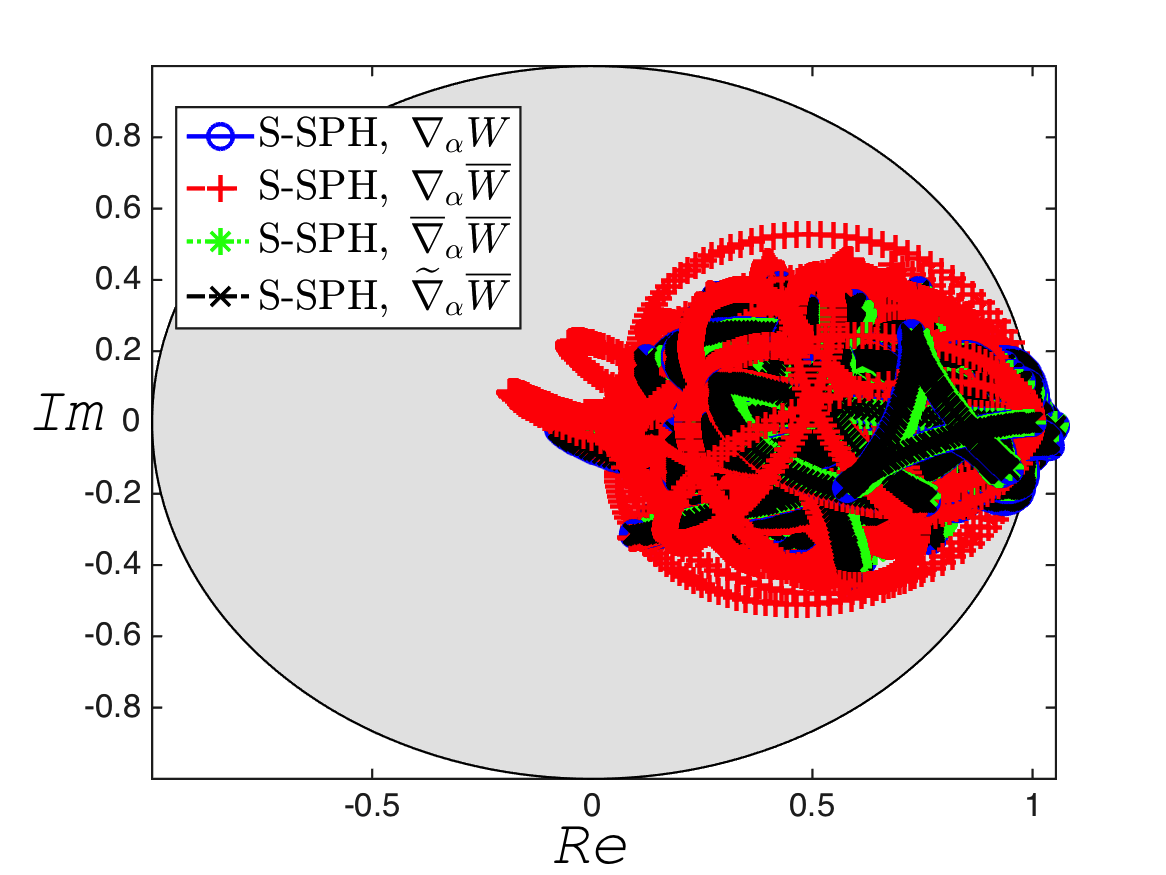}
		\hspace{-4.0cm}
		\caption*{Random}
		\label{Stab-S-RandBox}
	\end{minipage}%
	\begin{minipage}{1.6in}%
		\centering
		\hspace{-4.0cm}
		\includegraphics[width=4cm, height=4cm]{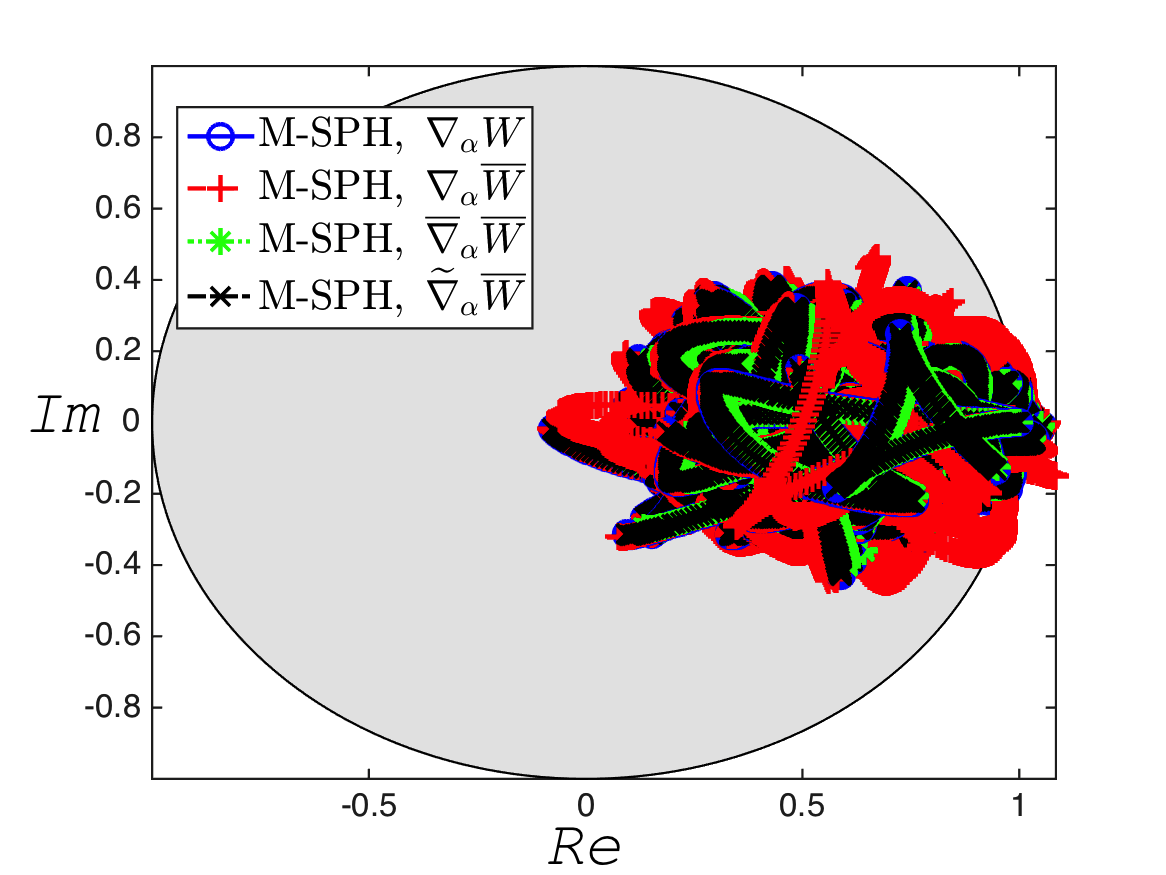}
		\hspace{-4.0cm}
		\caption*{Random}
		\label{Stab-M-RandBox}
	\end{minipage}%
}
	\caption{
		von Neumann growth factor for different discretization schemes. Corrected Brookshaw's scheme 
		(CB-SPH) is given by \eqref{Brookshaw1M} with the 
		correction multiplier, while  Schwaiger's scheme (S-SPH) is given by 
		\eqref{Schwaiger1}--\eqref{Schwaiger5}. New approximation 
		(M-SPH) considered here are the SPH form \eqref{NewScheme1}--\eqref{NewScheme3}. Four different 
		options of computing kernel 
		gradients (i.e., $\nabla_{\gamma}W$, 
		$\nabla_{\alpha}\overline{W}$, and corrected kernel gradients 
		(i.e., $\nabla^{*}_{\gamma}W$,
		$\nabla^{*}_{\alpha}\overline{W}$,
		$\overline{\nabla}^{*}_{\alpha}\overline{W}$,
		and  $\widetilde{\nabla}^{*}_{\alpha}\overline{W}$) are
		shown.}  
	\label{StabilityAnalysisDiffScheme}
\end{figure}

For linear PDEs, there is the Lax-equivalence 
theorem which connects the consistency and stability with 
the convergence. The idea of the von Neumann stability analysis is to study 
the growth of waves $\lambda e^{i\vect{k}\cdot \vect{r}}$ (similar to Fourier 
methods). After applying one of the discretization methods 
above to the Laplace operator \eqref{LV_SPH_Lap}, the following relation can 
be written:
\begin{equation}\label{LV_SPH_Neumann1}
\displaystyle
\frac{\left\{\mathbf{u}\right\}^{n+1}-\left\{\mathbf{u}\right\}^{n}}{\tau}-
\mat{L}^{\tilde{h}}\left[\left\{\mathbf{u}\right\}^{n}\right]
= \left\{\vect{b}\right\},
\end{equation}
where $\tau$ is the iteration parameter (e.g., time step), $n$ is the
iteration index (e.g., time index),
$\left\{\mathbf{u}\right\}^{n}=\left(\vect{u}^{n}_{1},\ldots,\vect{u}^{n}_{N}\right)$
is a vector of all $N$ particle values $\vect{u}^{n}_{i}, \ i=1,\ldots,N$,
$\left\{\vect{b}\right\}=\left(\vect{b}_{1},\ldots,\vect{b}_{N}\right)$ is the
right-hand side vector. The expression \eqref{LV_SPH_Neumann1} may
represent, for example, discretization of the parabolic PDE or an iterative
solver of the linear system of equations arising from discretization of
elliptic boundary value problem. Substituting into the left-hand side of
the relation \eqref{LV_SPH_Neumann1} the following form of perturbation
\begin{equation}\label{LV_SPH_Neumann2}
	\displaystyle
	\vect{u}^{n}_{j}=
	\lambda^{n} e^{\textbf{i}\vect{k}_{j}\cdot\vect{r}_{j}}=
	\lambda^{n}\cdot
	\prod\limits^{n}_{l=1}
	e^{\textbf{i} k^{l}_{j}\cdot x^{l}_{j}}
\end{equation}
leads to the expression for the von Neumann growth factor
subject to \eqref{LV_SPH_Neumann1} and linear Laplace operator:
\begin{equation}\label{LV_SPH_Neumann3}
\displaystyle
\lambda_{j}\left(\tau\right) = 1 +
\tau\cdot e^{-\textbf{i}\vect{k}_{j}
\cdot\vect{r}_{j}} \cdot{\mat{L}}^{\tilde{h}}
\left[\left\{e^{\textbf{i}\vect{k}\cdot\vect{r}}\right\}\right],
\vect{r}\in\mathbb{R}^{n}, \ n=1,2,3;
\end{equation}
where
$
\left\{e^{\textbf{i}\vect{k}\cdot\vect{r}}\right\} =
\left(e^{\textbf{i}\vect{k}_{1}\cdot\vect{r}_{1}},
\ldots,e^{\textbf{i}\vect{k}_{N}\cdot\vect{r}_{N}}\right)
$
and $\lambda_{j}\left(\tau\right)$ is the von Neumann 
growth factor.
For the discretization to be stable, it is required that
$\left|\lambda_{j}\left(\tau\right)\right|\leq 1, \ \forall j$.  
The von Neumann growth factor is 
shown for three main discretization schemes: Corrected
Brookshaw's scheme (CB-SPH) \eqref{Brookshaw1M}, 
Schwaiger's scheme (S-SPH) \eqref{Schwaiger1}--\eqref{Schwaiger5} and 
new approximation (M-SPH)
\eqref{NewScheme1}--\eqref{NewScheme3} for uniform 
and pseudo random particle 
distribution with $\tau = 0.25$ and 
$\mathbf{M}^{\alpha\beta}\left(\vect{r}\right)
=\delta_{\alpha\beta}$ (see, Figures \ref{StabilityAnalysisDiffScheme}).
In case of uniform particle distribution, the von Neumann 
growth factor clusters around the real axis for all schemes and satisfies the 
requirement $\left|\lambda_{j}\left(\tau\right)\right|\leq 1, \ \forall j$. In case of 
pseudo random particle distribution, the 
von Neumann growth factor has both real and imaginary 
parts forming complex shape but satisfying the requirement 
$\left|\lambda_{j}\left(\tau\right)\right|\leq 1, \ \forall j$ almost everywhere 
(i.e., it could be some problems at the boundary particles).
\subsection{Monotonicity and Convergence}
\label{MonotonicityConvergence}

In real life applications, the numerical domains are large and, therefore,
many degree of freedom (i.e., unknowns) is required. The resulting
matrix is usually sparse, but because of fill-in a direct method requires
significant amount of memory and time in general case. 
Hence, iterative solvers are used to overcome these issues. The 
iterative solvers converge only if the matrix satisfies certain properties related to 
the property of the diagonal dominance. However, these properties are only 
sufficient conditions for the method to converge.

The proof of the convergence of linear meshless schemes applied to a linear
elliptic boundary value problem can be done in the following  steps (see,
Bouchon (2007) \cite{Bouchon2007}, Bouchon and Peichl (2007)
\cite{BouchonPeichl2007}, Matsunaga and Yamamoto (2000)
\cite{MatsunagaYamamoto2000}, Thom\'ee (2001) \cite{Thomee2001}:
if it is shown that the truncation error $\varepsilon$ tends to 0 as the maximum
smoothing length $\max\limits_{I\in N}h_{I}$ goes to 0, then the 
linear system $\textbf{A}\delta u =\varepsilon$ that couples the variable error 
$\delta u$ with $\varepsilon$ proves the convergence of the schemes, given that 
the matrix of the discretized meshless operator is monotone. 
A square matrix $\mat{A} = \left(a_{ij}\right)_{1\leq i \leq n, 1\leq j \leq n}\in 
\mathbb{R}^{n\times n}$ is called monotone if $a_{ij}\leq 0 \ \forall i\neq j$, 
$a_{ii} > 0 \ \forall i$ and it is inverse positive $\mat{A}^{-1}\geq 0$.
Furthermore, the monotone schemes do satisfy a discrete maximum 
principle producing solutions without spurious oscillations.

It is clear that the new scheme \eqref{NewScheme1}--\eqref{NewScheme3} 
can be written in the form (\ref{LV_TotalFluxDiscrForm}), where meshless transmissibility 
between particles $\displaystyle\mathbf{r}_{J}$ and $\displaystyle\mathbf{r}_{I}$ can be 
defined as $\displaystyle T\left(\mathbf{r}_{J}, \mathbf{r}_{I}\right)= 
\frac{\bar\Gamma_{\beta\beta}}{n}\left(m_{I}+m_{J}\right)\bar{T}_{IJ}$ with:
\begin{equation}
	\label{FullTransmissibility}
	\begin{array}{lcl}
		\displaystyle
		\bar{T}_{IJ}=
		\displaystyle
		\left\lbrace
		\frac{\left(\vect{r}_{J}-\vect{r}_{I}\right)\cdot\overline{\nabla W}
		\left(\vect{r}_{J}-\vect{r}_{I}, \tilde{h}_{IJ}\right)}{\left\Vert\vect{r}_{J}-
		\vect{r}_{I}\right\Vert^{2}}-
	  \overline{\nabla^{*}_{\alpha}W}\left(\vect{r}_{J}-\vect{r}_{I}, \tilde{h}_{IJ}\right)
		\mathbf{N}^{\alpha}\right\rbrace,
	\end{array}
\end{equation}
where $\mathbf{N}^{\alpha}\left(\vect{r}_{I}\right)$ is defined by 
\eqref{Schwaiger2}. This allows us to formulate the following
remark. 
\begin{remark}\label{Remark_Monotonicity_Convergence}
	Taking into account Remark \ref{Remark_Gamma_Structure} and the following 
	relation:
    \begin{equation}
		\label{EffectiveTransmissibility1}
		\begin{array}{lcl}
		\displaystyle\vspace{0.3cm}
        \frac{\left(\vect{r}_{J}-\vect{r}_{I}\right)\cdot
	    \overline{\nabla W}\left(\vect{r}_{J}-\vect{r}_{I}, h\right)}{\left
		\Vert\vect{r}_{J}-\vect{r}_{I}\right\Vert^{2}}=
	    \frac{1}{z\cdot h^{2}}\frac{\overline{dW}}{dz}\left(z\right) \le 0,
		\end{array}
	\end{equation}
	and the relation \eqref{EffectiveGradWN} which all together lead to 
	the fact that there is a parameter $\tilde{h}_{I}$ such that 
	$T\left(\mathbf{r}_{J}, \mathbf{r}_{I}\right) \ge 0$. This means that  
	the proposed scheme subject to this parameter $\tilde{h}_{I}$ is 
	monotone for the medium with the scalar heterogeneous coefficients.
\end{remark}

\section{Solution of Boundary Value Problems}
\label{SolutionBoundaryValueProblems}

Following the work by \cite{Schwaiger2008}, the
numerical tests for inhomogeneous Dirichlet, Neumann and mixed 
boundary conditions are considered for homogeneous and heterogeneous media 
with the characteristics $\mat{M}\left(\vect{r}\right)$. To illustrate the performance 
of the proposed scheme, a modelling of a single phase steady-state fluid flow in 
fully anisotropic porous media with different type of boundary conditions 
is also presented in this section. The square 2D and 3D domains \eqref{LV_SPH_Box} 
are considered.

The relative error used to quantify the accuracy of the 
proposed schemes in the subsections \eqref{BrookshawScheme}, 
\eqref{SchwaigerScheme}, and \eqref{NewScheme} during numerical simulation 
in this section is given by:
\begin{equation}\label{ChapIV6_1}
\displaystyle\left\Vert E_{R} \right\Vert^{L_{2}}_{\Omega} =
\left[
\frac{1}{\sum\limits_{\Omega}V_{\vect{r}_{J}}}
\sum\limits_{\Omega}V_{\vect{r}_{K}}
\left(\frac{\mathbf{u}\left(\vect{r}_{K}\right) -
\langle \mathbf{u}
\left(\vect{r}_{K}\right)\rangle}{\mathbf{u}\left(\vect{r}_{K}\right)}
\right)^{2}
\right]^{1/2}
\end{equation}
where $\mathbf{u}\left(\vect{r}\right)$ is the analytical or reference solution
field and $\langle \mathbf{u}\left(\vect{r}_{K}\right)\rangle$ is the approximated
solution field. Several numerical results using considered in this 
paper schemes for uniform and pseudo random particle distributions are shown 
in the following sections that confirm the theoretical results from the previous 
sections.

\begin{figure}[!htbp]
	\centering
	\includegraphics[width=4.0in,height=2.5in]{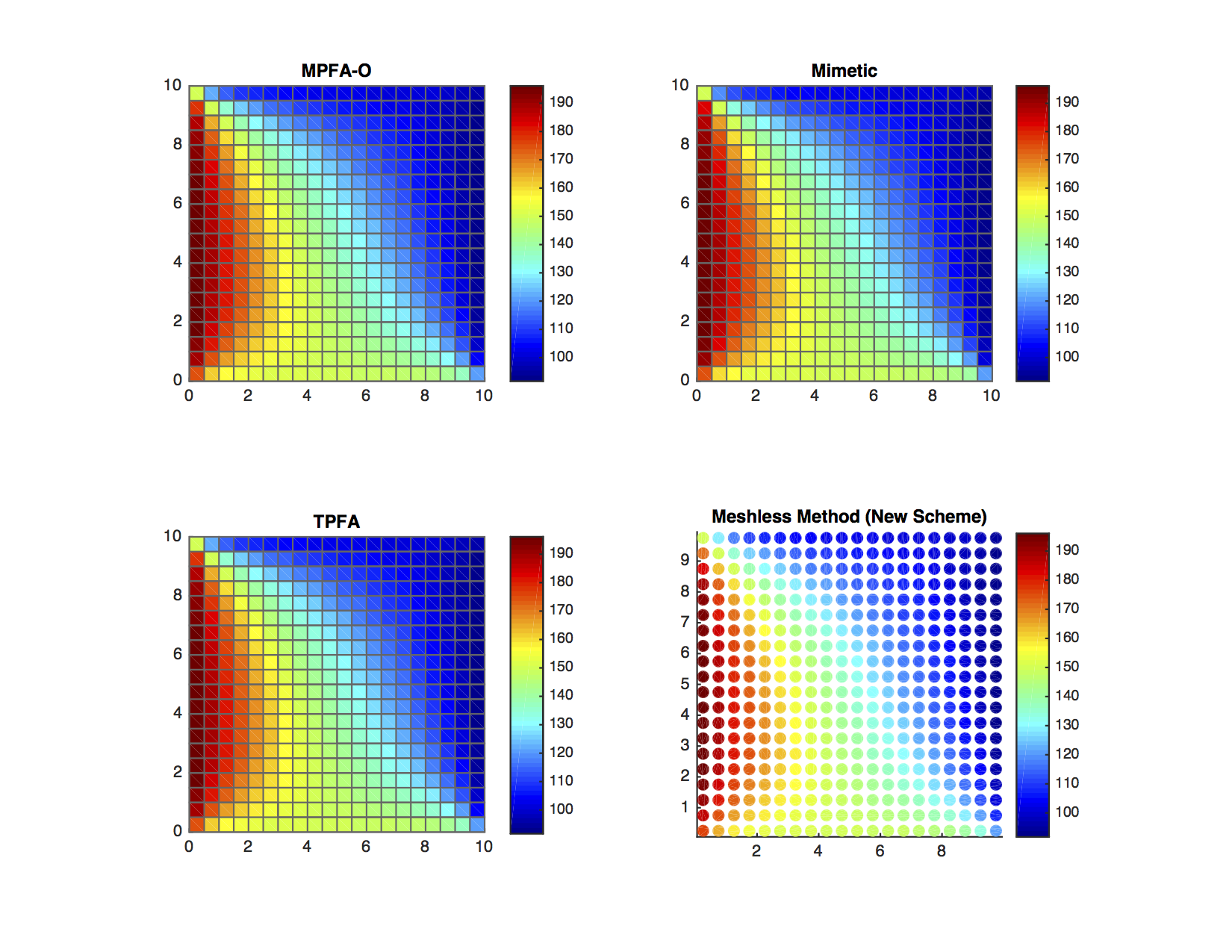}
	\caption{
		Comparison of solutions of the Dirichlet problem for 
		the Laplace equation. The method used are MPFA-O, 
	    Mimetic, TPFA, and Meshless methods (New Method: 
	    \eqref{NewScheme1}--\eqref{NewScheme3}).
	}
	\label{NSchDirichPrblm}
\end{figure}

\subsection{Inhomogeneous Boundary Condition}

Firstly, the homogeneous properties of the porous media 
$\mat{M}\left(\vect{r}\right)=\mat{I}$ are assumed for simplicity of
the derivation. The analytical solution of \eqref{LV_SPH_Lap} 
subject to the assumption that $g\left(\vect{r}\right)\equiv 0$ for 
$\vect{r}\in \Omega\subset \mathbb{R}^{2}$ and constant boundary 
conditions is the following:
\begin{equation}\label{ChapIV6_2}
\begin{array}{lcl}
\vspace{0.3cm}
\displaystyle \mathbf{u}\left(x,y\right)=
\sum\limits_{n=odd}
\left[
\frac{4\psi_{1}}{n\pi \sinh\left(\frac{n\pi H}{L}\right)}
\right]
\sin\left(\frac{n\pi}{L}x\right)
\sinh\left(\frac{n\pi}{L}\left(H-y\right)\right)+\\
\vspace{0.3cm}
\displaystyle+\sum\limits_{n=odd}
\left[
\frac{4\psi_{2}}{n\pi \sinh\left(\frac{n\pi L}{H}\right)}
\right]
\sin\left(\frac{n\pi}{H}y\right)
\sinh\left(\frac{n\pi}{H}x\right)+\\
\vspace{0.3cm}
\displaystyle+\sum\limits_{n=odd}
\left[
\frac{4\psi_{3}}{n\pi \sinh\left(\frac{n\pi L}{H}\right)}
\right]
\sin\left(\frac{n\pi}{H}y\right)
\sinh\left(\frac{n\pi}{H}x\right)+\\
\displaystyle+\sum\limits_{n=odd}
\left[
\frac{4\psi_{4}}{n\pi \sinh\left(\frac{n\pi L}{H}\right)}
\right]
\sin\left(\frac{n\pi}{H}y\right)
\sinh\left(\frac{n\pi}{H}\left(L-x\right)\right),
\end{array}
\end{equation}
where $\psi_{1}$ is the boundary at $y = 0$, $\psi_{2}$ at $x = L$, $\psi_{3}$
at $y = H$, and $\psi_{4}$ at $x = 0$.

Tables \ref{tab:ConvLaplace2DU1}, \ref{tab:ConvLaplace2DU2} 
shows the convergence rate of different schemes for uniform 
particle distribution. The solutions with MPFA and Mimetic schemes were obtained 
using MATLAB Reservoir Simulation Toolbox (MRST) \cite{mrst2016}. It is clear from 
the  convergence results that schemes \eqref{Schwaiger1}--\eqref{Schwaiger5} and 
\eqref{NewScheme1}--\eqref{NewScheme3} are identical and they overperform the 
scheme (\ref{Brookshaw1M}). The convergence rate as was predicted theoretically 
is at least $\mathcal{O}\left(h^\omega\right), \ 1\leq\omega \leq 2$.
\begin{table}[!htbp]
	\begin{center}
		\caption{The error of convergence for different schemes (uniform 
		particle distribution) and $f=0.5005$.}
		\label{tab:ConvLaplace2DU1}       
		\begin{tabular}{llllll}
			\hline\noalign{\smallskip}	
			DoF  & MPFA  & Mimetic & Brookshaw \eqref{Brookshaw1M} & Schwaiger \eqref{Schwaiger1}--\eqref{Schwaiger5} & New Scheme \eqref{NewScheme1}--\eqref{NewScheme3}  \\
			\noalign{\smallskip}\hline\noalign{\smallskip}
			25     & $1.318\cdot 10^{-1}$ & $1.240\cdot 10^{-1}$  & $1.314\cdot 10^{0}$  & $6.608\cdot 10^{-2}$  &  $6.608\cdot 10^{-2}$  \\
			100    & $3.296\cdot 10^{-2}$ & $3.078\cdot 10^{-2}$  & $3.925\cdot 10^{-1}$ & $1.865\cdot 10^{-2}$  &  $1.865\cdot 10^{-2}$  \\
			400    & $8.233\cdot 10^{-3}$ & $7.685\cdot 10^{-3}$  & $1.105\cdot 10^{-1}$ & $4.714\cdot 10^{-3}$  &  $4.714\cdot 10^{-3}$  \\
			1600   & $2.058\cdot 10^{-3}$ & $1.920\cdot 10^{-3}$  & $3.025\cdot 10^{-2}$ & $1.179\cdot 10^{-3}$  &  $1.179\cdot 10^{-3}$  \\
			6400   & $5.318\cdot 10^{-4}$ & $4.968\cdot 10^{-4}$  & $8.152\cdot 10^{-3}$ & $3.213\cdot 10^{-4}$  &  $3.213\cdot 10^{-4}$  \\
			25600  & $3.007\cdot 10^{-4}$ & $2.955\cdot 10^{-4}$  & $2.193\cdot 10^{-3}$ & $2.807\cdot 10^{-4}$  &  $2.807\cdot 10^{-4}$  \\
			\noalign{\smallskip}\hline
		\end{tabular}
	\end{center}
\end{table}
\begin{table}[!htbp]
	\begin{center}
		\caption{The error of convergence for different schemes (uniformed particle distribution) and $f=1.001$.}
		\label{tab:ConvLaplace2DU2} 
		\begin{tabular}{llll}
			\hline\noalign{\smallskip}
			DoF  & Brookshaw \eqref{Brookshaw1M} & Schwaiger \eqref{Schwaiger1}--\eqref{Schwaiger5} & New Scheme \eqref{NewScheme1}--\eqref{NewScheme3}  \\			
			\noalign{\smallskip}\hline\noalign{\smallskip}
			25     & $2.312\cdot 10^{-1}$ & $1.091\cdot 10^{-1}$ & $1.091\cdot 10^{-1}$ \\
			100    & $1.120\cdot 10^{-1}$ & $2.698\cdot 10^{-2}$ & $2.698\cdot 10^{-2}$ \\
			400    & $3.998\cdot 10^{-2}$ & $6.737\cdot 10^{-3}$ & $6.737\cdot 10^{-3}$ \\
			1600   & $1.262\cdot 10^{-2}$ & $1.684\cdot 10^{-3}$ & $1.684\cdot 10^{-3}$ \\
			6400   & $3.737\cdot 10^{-3}$ & $4.418\cdot 10^{-4}$ & $4.418\cdot 10^{-4}$ \\
			25600  & $1.097\cdot 10^{-3}$ & $2.919\cdot 10^{-4}$ & $2.919\cdot 10^{-4}$ \\
			\noalign{\smallskip}\hline	
		\end{tabular}
	\end{center}
\end{table}
Tables \ref{tab:ConvLaplace2DIR1}, \ref{tab:ConvLaplace2DIR2} shows the convergence 
rate for different schemes and pseudo random distribution. The reason 
for starting from $f=0.5005$ is that we want to force the same number of non-zero entries 
in the matrix as for MPFA or Memetic methods. This procedure may lead to a badly 
conditioned matrix since $\overline{\nabla W}\left(\vect{r}_{J}-\vect{r}_{I}, h\right)\ll 1$. 
This distribution of particles was generated by perturbing the regularly 
distributed particles using a uniform random variable varying between $+10\%$ and $-10\%$ 
of the maximum smoothing length $h_{max}=1.001$. Since random realizations are used to 
compute errors, the statistical data about 30 realizations are used to compute mean values 
and standard deviations. This data is presented in Tables 
\ref{tab:ConvLaplace2DIR1}, \ref{tab:ConvLaplace2DIR2}. The positive 
trends of the scheme remain the same as for the case of uniform particle distributions. 
The scheme \eqref{Brookshaw1M} has higher error compare to schemes
\eqref{Schwaiger1}--\eqref{Schwaiger5} and  \eqref{NewScheme1}--\eqref{NewScheme3}. 
However, the dispersion of the approximation error is higher
in the scheme \eqref{NewScheme1}--\eqref{NewScheme3}.
\begin{table}[!htbp]
	\begin{center}
		\caption{The error of convergence for different schemes (random 
		particle distribution) and $f=0.6006$.}
		\label{tab:ConvLaplace2DIR1} 
		\begin{tabular}{llll}
			\hline\noalign{\smallskip}
			DoF  & Brookshaw \eqref{Brookshaw1M} & Schwaiger \eqref{Schwaiger1}--\eqref{Schwaiger5} & New Scheme \eqref{NewScheme1}--\eqref{NewScheme3}  \\			
			\noalign{\smallskip}\hline\noalign{\smallskip}
			25     & $1.102\cdot 10^{+0}\pm  1.966\cdot 10^{-2}$ & $1.940\cdot 10^{-1}\pm 3.954\cdot 10^{-2}$ & $3.707\cdot 10^{-1}\pm 1.319\cdot 10^{-1}$ \\
			100    & $3.361\cdot 10^{-1}\pm 4.923\cdot 10^{-3}$  & $5.511\cdot 10^{-2}\pm 1.205\cdot 10^{-2}$ & $9.074\cdot 10^{-2}\pm 2.668\cdot 10^{-2}$ \\
			400    & $9.576\cdot 10^{-2}\pm 8.705\cdot 10^{-4}$  & $1.376\cdot 10^{-2}\pm 2.440\cdot 10^{-3}$ & $2.299\cdot 10^{-2}\pm 6.972\cdot 10^{-3}$ \\
			1600   & $2.640\cdot 10^{-2}\pm 2.551\cdot 10^{-4}$  & $4.157\cdot 10^{-3}\pm 7.485\cdot 10^{-4}$ & $5.912\cdot 10^{-3}\pm 1.289\cdot 10^{-3}$ \\
			6400   & $7.185\cdot 10^{-3}\pm 5.982\cdot 10^{-5}$  & $1.044\cdot 10^{-3}\pm 1.465\cdot 10^{-4}$ & $1.406\cdot 10^{-3}\pm 3.090\cdot 10^{-4}$ \\
			25600  & $1.944\cdot 10^{-3}\pm 1.451\cdot 10^{-5}$  & $3.858\cdot 10^{-4}\pm 5.887\cdot 10^{-5}$ & $4.960\cdot 10^{-4}\pm 8.425\cdot 10^{-5}$ \\
			\noalign{\smallskip}\hline
		\end{tabular}
	\end{center}
\end{table}
\begin{table}[!htbp]
	\begin{center}
		\caption{The error of convergence for different schemes (random particle 
		distribution) and $f=1.2012$.}
		\label{tab:ConvLaplace2DIR2}       
		\begin{tabular}{llll}
			\hline\noalign{\smallskip}
			DoF  & Brookshaw \eqref{Brookshaw1M} & Schwaiger \eqref{Schwaiger1}--\eqref{Schwaiger5} & New Scheme \eqref{NewScheme1}--\eqref{NewScheme3}  \\			
			\noalign{\smallskip}\hline\noalign{\smallskip}
			25    & $3.578\cdot 10^{-1}\pm 2.879\cdot 10^{-2}$ & $2.162\cdot 10^{-1}\pm 2.800\cdot 10^{-2}$ & $1.842\cdot 10^{-1}\pm 8.871\cdot 10^{-3}$ \\
			100   & $1.261\cdot 10^{-1}\pm 4.859\cdot 10^{-3}$ & $5.696\cdot 10^{-2}\pm 5.518\cdot 10^{-3}$ & $4.596\cdot 10^{-2}\pm 1.999\cdot 10^{-3}$ \\
			400   & $4.224\cdot 10^{-2}\pm 1.092\cdot 10^{-3}$ & $1.493\cdot 10^{-2}\pm 1.546\cdot 10^{-3}$ & $1.157\cdot 10^{-2}\pm 5.691\cdot 10^{-4}$ \\
			1600  & $1.299\cdot 10^{-2}\pm 2.031\cdot 10^{-4}$ & $3.934\cdot 10^{-3}\pm 3.969\cdot 10^{-4}$ & $2.889\cdot 10^{-3}\pm 1.451\cdot 10^{-4}$ \\
			6400  & $3.862\cdot 10^{-3}\pm 3.452\cdot 10^{-5}$ & $1.037\cdot 10^{-3}\pm 8.443\cdot 10^{-5}$ & $7.507\cdot 10^{-4}\pm 4.060\cdot 10^{-5}$ \\
			25600 & $1.133\cdot 10^{-3}\pm 2.125\cdot 10^{-5}$ & $4.087\cdot 10^{-4}\pm 5.405\cdot 10^{-5}$ & $3.475\cdot 10^{-4}\pm 3.814\cdot 10^{-5}$ \\
			\noalign{\smallskip}\hline
		\end{tabular}
	\end{center}
\end{table}
Figure \ref{NSchDirichPrblm} shows the numerical solution of the boundary value problem
obtained using different discretization methods. The solution was obtained using 40
particles/cells in each direction.

\subsection{Inhomogeneous Mixed Boundary Condition Test}

The general steady-state solution for a Dirichlet condition
along the base of the plate and Neumann conditions elsewhere is given in 
\cite{Schwaiger2008} by
\begin{equation}\label{ChapIV7_1}
	\displaystyle \mathbf{u}\left(x,y\right)= \psi_{1} + \psi_{3}y + \sum\limits^{\infty}_{n=1}
	\frac{2\psi_{2}\cosh\left(\lambda_{n}x\right) + 
	2\psi_{4}\cosh\left(\lambda_{n}\left(L-x\right)\right)}{H\lambda^{2}_{n}
	\sinh\left(\lambda_{n}L\right)}
	\sinh\left(\lambda_{n}y\right),
\end{equation}
where $\psi_{1}$ is the boundary value of $\mathbf{u}$ at $y = 0$,
$\psi_{2}$ is the flux at $x = L$, $\psi_{3}$ is the flux at $y = H$,
and $\psi_{4}$ is the flux at $x = 0$,  and 
$\displaystyle\lambda_{n} = \frac{\left(2n-1\right)\pi}{2H}$.
The following parameters (dimensionless) were chosen:$\psi_{1}=150$, 
$\psi_{3}=150$, $\psi_{4}=200$, $\psi_{2}=90$. Again, for uniform and 
regular particles distribution, MPFA method reduces to the TPFA. Similar 
to the previous section, the regular and pseudo random 
particle distribution are considered to access the convergence properties 
of the proposed schemes. 
\begin{remark}\label{Remark_Neumann_BC}
	Let $\Gamma_{N}=\varnothing$ and $\Gamma_{D}=\varnothing$, 
	$g\left(\vect{r}\right) \geq 0, \ \vect{r}\in\Omega$,
	$g_{N}\left(\vect{r}\right) \geq 0, \ \vect{r}\in\Gamma_{N}$, 
	$g_{D}\left(\vect{r}\right) \geq 0, \ \vect{r}\in\Gamma_{D}$
	and the solution of 
	\begin{equation}
	\label{LV_SPH_Lap1}
		\displaystyle
		\vspace{0.2cm}
		\mathbf{L}\left(\mathbf{u}\right) = 0 , \ \mbox{or} \ \ \
		\nabla\left(\mathbf{M}\left(\vect{r}\right)\nabla \mathbf{u}
		\left(\vect{r}\right)\right) = g\left(\vect{r}\right), \
		\forall \vect{r}\in\Omega\subset\mathbb{R}^{n}, 
	\end{equation}	
	exists then it can be discretized in all internal particles by the following schemes
	Brookshaw \eqref{Brookshaw1M},  Schwaiger \eqref{Schwaiger1}--\eqref{Schwaiger5}, and 
	New Scheme \eqref{NewScheme1}--\eqref{NewScheme3} and compounded 
	with the following condition for all $\vect{r}\in\Gamma_{N}$:
	\begin{equation}\label{ChapIV5_4}
	\langle\vect{u}^{\alpha}\vect{n}_{\alpha}\rangle = 
	\vect{n}_{\alpha}\mat{M}^{\alpha\beta}\left(\vect{r}\right)
	\sum\limits_{\Omega_{\vect{r},\tilde{h}}}V_{\vect{r}_{J}}
	\left[
	\mathbf{u}\left(\vect{r}_{J}\right) -
	\mathbf{u}\left(\vect{r}_{I}\right)
	\right]
	\overline{\nabla^{*}_{\beta} W}\left(\vect{r}_{J}-\vect{r}_{I}, h\right)  
	\end{equation}
	where $\vect{n}_{\alpha}$ are the component of the external normal to the boundary 
	$\vect{r}\in\Gamma_{N}$.
\end{remark}

The following smoothing multiplication factors $f=0.5005$, $f=1.001$, and 
$f=2.002$  are considered. Tables \ref{tab:ConvLaplace2DUN1}, 
\ref{tab:ConvLaplace2DUN2} show the convergence rate for different schemes 
using the uniform and pseudo random particle 
distribution. The  convergence rate as was predicted theoretically is at least 
$\mathcal{O}\left(h^\omega\right), \ 1\leq\omega < 2$. Interestingly, 
the scheme (\ref{Brookshaw1M}) shows a very bad convergence rate for $f=0.5005$.
This explains by the fact that we have one particle from each side in the kernel support 
almost next to the boundary of the kernel support leading to 
$\overline{\nabla W}\left(\vect{r}_{J}-\vect{r}_{I}, h\right)\ll 1$ and, hence, to a badly
conditioned matrix. This does not observed for other schemes due to the normalization 
coefficient. As a result, the scheme (\ref{Brookshaw1M}) should be used with more than 
one particle in the compact support in each direction leading to a larger bandwidth in the 
matrix.      

\begin{table}[!htbp]
	\begin{center}
		\caption{The error of convergence for different schemes (uniform 
		particle distribution) and $f=0.5005$.}
		\label{tab:ConvLaplace2DUN1}       
		\begin{tabular}{llllll}
			\hline\noalign{\smallskip}	
			DoF  & MPFA  & Mimetic & Brookshaw \eqref{Brookshaw1M} & Schwaiger 
			\eqref{Schwaiger1}--\eqref{Schwaiger5} & New Scheme \eqref{NewScheme1}--\eqref{NewScheme3}  \\			
			\noalign{\smallskip}\hline\noalign{\smallskip}	
			25    &  $1.489\cdot 10^{-1}$ & $2.651\cdot 10^{-2}$ & $1.489\cdot 10^{6}$ & $8.244\cdot 10^{-1}$ &$8.244\cdot 10^{-1}$  \\
			100   & $2.015\cdot 10^{-2}$ & $4.424\cdot 10^{-3}$ & $9.471\cdot 10^{5}$ &$2.178\cdot 10^{-1}$ &$2.178\cdot 10^{-1}$ \\
			400   &  $2.727\cdot 10^{-3}$ & $6.889\cdot 10^{-4}$ &$5.378\cdot 10^{5}$ &$5.550\cdot 10^{-2}$ &$5.550\cdot 10^{-2}$ \\
			1600  & $3.666\cdot 10^{-4}$ & $1.024\cdot 10^{-4}$ &$2.859\cdot 10^{5}$ & $1.398\cdot 10^{-2}$ &$1.398\cdot 10^{-2}$ \\
			6400 &  $4.951\cdot 10^{-5}$ & $1.522\cdot 10^{-5}$ & $1.473\cdot 10^{5}$ & $3.504\cdot 10^{-3}$ &  $3.504\cdot 10^{-3}$   \\
			25600 &  $7.309\cdot 10^{-6}$ & $2.920\cdot 10^{-6}$ & $7.474\cdot 10^{4}$ & $8.770\cdot 10^{-4}$ & $8.770\cdot 10^{-4}$  \\
			\noalign{\smallskip}\hline
		\end{tabular}
	\end{center}
\end{table}
\begin{table}[!htbp]
	\begin{center}
		\caption{The error of convergence for different schemes (uniform particle 
		distribution) and $f=1.001$.}
		\label{tab:ConvLaplace2DUN2}       
		\begin{tabular}{llll}
			\hline\noalign{\smallskip}	
			DoF  & Brookshaw \eqref{Brookshaw1M} & Schwaiger \eqref{Schwaiger1}--\eqref{Schwaiger5} & New Scheme \eqref{NewScheme1}--\eqref{NewScheme3}  \\						
			\noalign{\smallskip}\hline\noalign{\smallskip}	
			25         & $1.163\cdot 10^{0}$     &  $1.078\cdot 10^{-1}$      &  $1.078\cdot 10^{-1}$      \\
			100       & $1.208\cdot 10^{0}$     &  $1.488\cdot 10^{-2}$     &  $1.488\cdot 10^{-2}$     \\
			400       & $8.291\cdot 10^{-1}$    &  $2.244\cdot 10^{-3}$    &  $2.244\cdot 10^{-3}$     \\
			1600     & $4.831\cdot 10^{-1}$    &  $4.215\cdot 10^{-4}$     &  $4.215\cdot 10^{-4}$     \\
			6400     & $2.606\cdot 10^{-1}$   &  $9.858\cdot 10^{-5}$    &  $9.858\cdot 10^{-5}$     \\
			25600   & $1.353\cdot 10^{-1}$    &  $2.502\cdot 10^{-5}$    &  $2.502\cdot 10^{-5}$       \\
			\noalign{\smallskip}\hline	
		\end{tabular}
	\end{center}
\end{table}
Tables \ref{tab:ConvLaplace2DIRN1}, \ref{tab:ConvLaplace2DIRN2}  
show the convergence rate for different
schemes with pseudo random particle distribution. 
The pseudo random distribution of particles was again generated 
by perturbing the regularly distributed particles using the uniform 
random variable varying between $+10\%$ and $-10\%$ of the maximum smoothing 
length $h_{max}=1.001$. Similar to the above case, the statistical 
data about 30 realizations are used to compute mean values and standard deviations. 
These data are presented in Tables \ref{tab:ConvLaplace2DIRN1}, \ref{tab:ConvLaplace2DIRN2}. 
The general trends of the scheme prosperities remain the same as for uniform particle 
distributions. The scheme (\ref{Brookshaw1M}) has higher error 
compared to schemes \eqref{Schwaiger1}--\eqref{Schwaiger5} 
and \eqref{NewScheme1}--\eqref{NewScheme3}.
However, the dispersion of the approximation error is higher in the scheme
\eqref{NewScheme1}--\eqref{NewScheme3}.
\begin{table}[!htbp]
	\begin{center}
		\caption{The error of convergence for different schemes (random particle 
		distribution) and $f=0.6006$.}
		\label{tab:ConvLaplace2DIRN1}       
		\begin{tabular}{llll}
			\hline\noalign{\smallskip}
			DoF  & Brookshaw \eqref{Brookshaw1M} & Schwaiger \eqref{Schwaiger1}--\eqref{Schwaiger5} & New Scheme \eqref{NewScheme1}--\eqref{NewScheme3}  \\						
			\noalign{\smallskip}\hline\noalign{\smallskip}
			25            &  $5.608\cdot 10^{1}\pm  1.019\cdot 10^{0}$ &   $8.815\cdot 10^{-1}\pm 4.248\cdot 10^{-1}$   &  $1.242\cdot 10^{0}\pm 5.058\cdot 10^{-1}$    \\
			100          &  $3.769\cdot 10^{1}\pm 5.708\cdot 10^{-1}$   &   $3.454\cdot 10^{-1}\pm 1.977\cdot 10^{-1}$   &  $4.361\cdot 10^{-1} \pm  2.425\cdot 10^{-1}$  \\
			400          &  $2.154\cdot 10^{1}\pm 1.757\cdot 10^{-1}$  &   $8.454\cdot 10^{-2}\pm 5.111\cdot 10^{-2}$   &  $1.525\cdot 10^{-1}\pm 1.195\cdot 10^{-1}$  \\
			1600        &  $1.159\cdot 10^{1}\pm 7.457\cdot 10^{-2}$   &   $2.245\cdot 10^{-2}\pm 1.055\cdot 10^{-2}$   &  $4.563\cdot 10^{-2}\pm 2.328\cdot 10^{-2}$   \\
			6400        &  $6.027\cdot 10^{0}\pm 2.312\cdot 10^{-2}$   &   $1.104\cdot 10^{-2}\pm 7.523\cdot 10^{-3}$   &  $1.848\cdot 10^{-2}\pm 9.276\cdot 10^{-3}$  \\
			25600      &  $3.054\cdot 10^{0}\pm 1.174\cdot 10^{-2}$    &   $4.123\cdot 10^{-3}\pm 2.059\cdot 10^{-3}$   &  $5.879\cdot 10^{-3}\pm 4.253\cdot 10^{-3}$  \\			
			\noalign{\smallskip}\hline				
		\end{tabular}
	\end{center}
\end{table}
\begin{table}[!htbp]
	\begin{center}
		\caption{The error of convergence for different schemes (random particle 
		distribution) and $f=1.2012$.}
		\label{tab:ConvLaplace2DIRN2}       
		\begin{tabular}{llll}
			\hline\noalign{\smallskip}
			DoF  & Brookshaw \eqref{Brookshaw1M} & Schwaiger \eqref{Schwaiger1}--\eqref{Schwaiger5} & New Scheme \eqref{NewScheme1}--\eqref{NewScheme3}  \\						
			\noalign{\smallskip}\hline\noalign{\smallskip}
			25       & $5.177\cdot 10^{-1}\pm 7.312\cdot 10^{-2}$   & $3.143\cdot 10^{-1}\pm 7.543\cdot 10^{-2}$     &  $4.302\cdot 10^{-1}\pm 1.461\cdot 10^{-1}$   \\
			100      & $5.160\cdot 10^{-1}\pm 4.035\cdot 10^{-2}$   & $5.959\cdot 10^{-2}\pm 4.713\cdot 10^{-3}$    &  $1.106\cdot 10^{-1}\pm 1.739\cdot 10^{-2}$ \\
			400     & $4.349\cdot 10^{-1}\pm 1.715\cdot 10^{-2}$   & $1.515\cdot 10^{-2}\pm 1.973\cdot 10^{-3}$    &  $3.015\cdot 10^{-2}\pm 1.045\cdot 10^{-2}$   \\
			1600    & $2.745\cdot 10^{-1}\pm 9.059\cdot 10^{-3}$   & $3.926\cdot 10^{-3}\pm 6.865\cdot 10^{-4}$   &  $1.133\cdot 10^{-2}\pm 5.486\cdot 10^{-3}$  \\
			6400   & $1.508\cdot 10^{-1}\pm 2.822\cdot 10^{-3}$    & $1.595\cdot 10^{-3}\pm 8.144\cdot 10^{-4}$   &  $3.187\cdot 10^{-3}\pm 1.438\cdot 10^{-3}$  \\
			25600 & $7.958\cdot 10^{-2}\pm 1.146\cdot 10^{-3}$  & $5.354\cdot 10^{-4} \pm 2.036\cdot 10^{-4}$  &  $9.201\cdot 10^{-4} \pm 5.201\cdot 10^{-4}$       \\
			\noalign{\smallskip}\hline
		\end{tabular}
	\end{center}
\end{table}
Results of the numerical solution are compared with
the series solution in Figure 11.
\begin{figure}[!htbp]
	\centering
	\includegraphics[width=4.0in,height=2.5in]{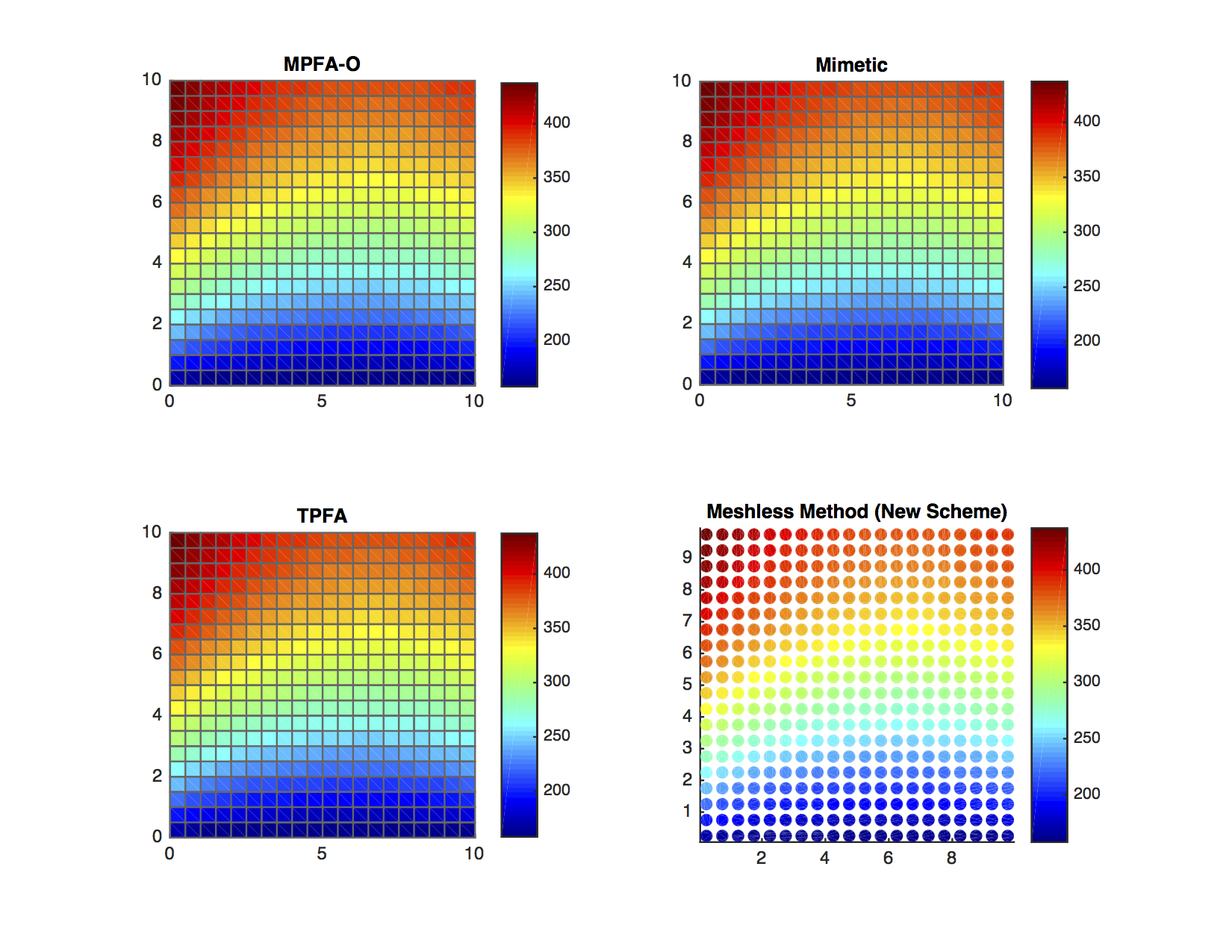}
	\caption{
		Comparison of solutions of the inhomogeneous mixed boundary problems for
		the Laplace equations. The method used are MPFA-O, Mimetic, TPFA, and Meshless
		methods (New Method: \eqref{NewScheme1}--\eqref{NewScheme3}).
	}       
	\label{MixedNSchDirichPrblm}
\end{figure}

\subsection{SPE10}
\begin{figure}
	\centering
	\hspace{-2.0cm}	
	\mbox{
		\begin{minipage}{1.6in}%
			\centering
			\hspace{-7.0cm}
	        \includegraphics[width=3.0in,height=2.0in]{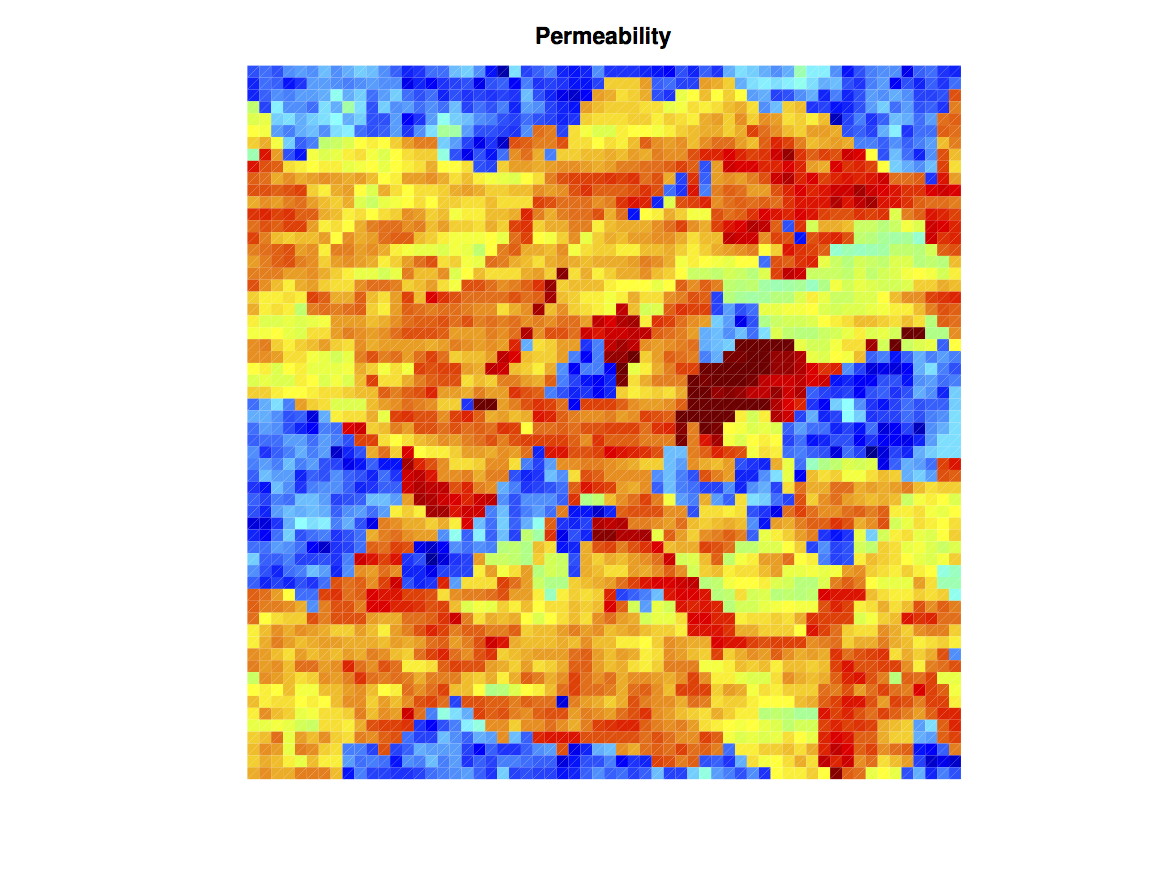}			
			\hspace{-7.0cm}
			\caption*{Case (a)}
			\label{SPE10Layer85_60by60}
		\end{minipage}%
		\qquad
		\qquad
		\qquad
		\qquad		
		\begin{minipage}{1.6in}%
			\centering
			\hspace{-7.0cm}
	        \includegraphics[width=3.0in,height=1.85in]{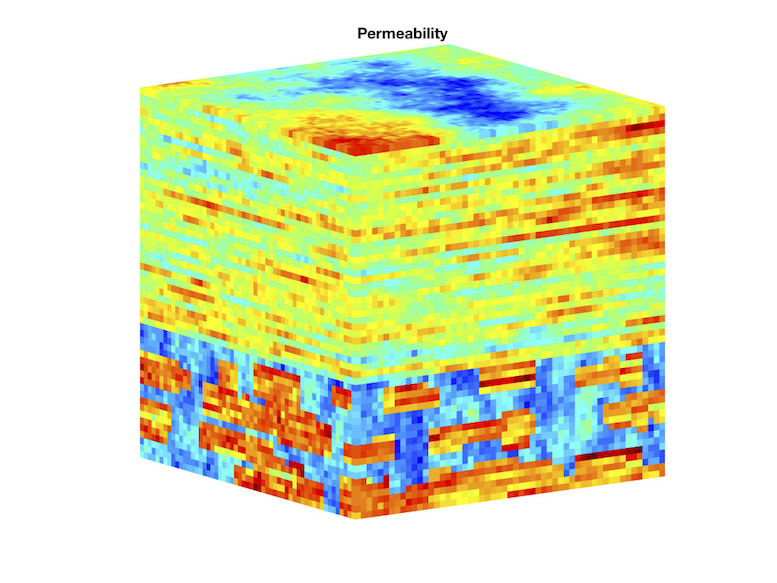}	
			\hspace{-7.0cm}
			\caption*{Case (b)}
			\label{SPE10Permeability3D_60by60}
		\end{minipage}%
		}
	   \caption{The lognormal permeability field in the SPE10 benchmark test.  
	   	Case (a) $ 60 \times 60 $cells of 85 layer.  Case (b) subsection of the 
	   	SPE10 model defined by  $60 \times 60 \times 60$ cells.}
	\label{SPE10Fig}
\end{figure}

\begin{figure}
	\mbox{
		\begin{minipage}{3in}%
			\centering
			\hspace{-4.0cm}
			\includegraphics[width=5cm, height=4cm]{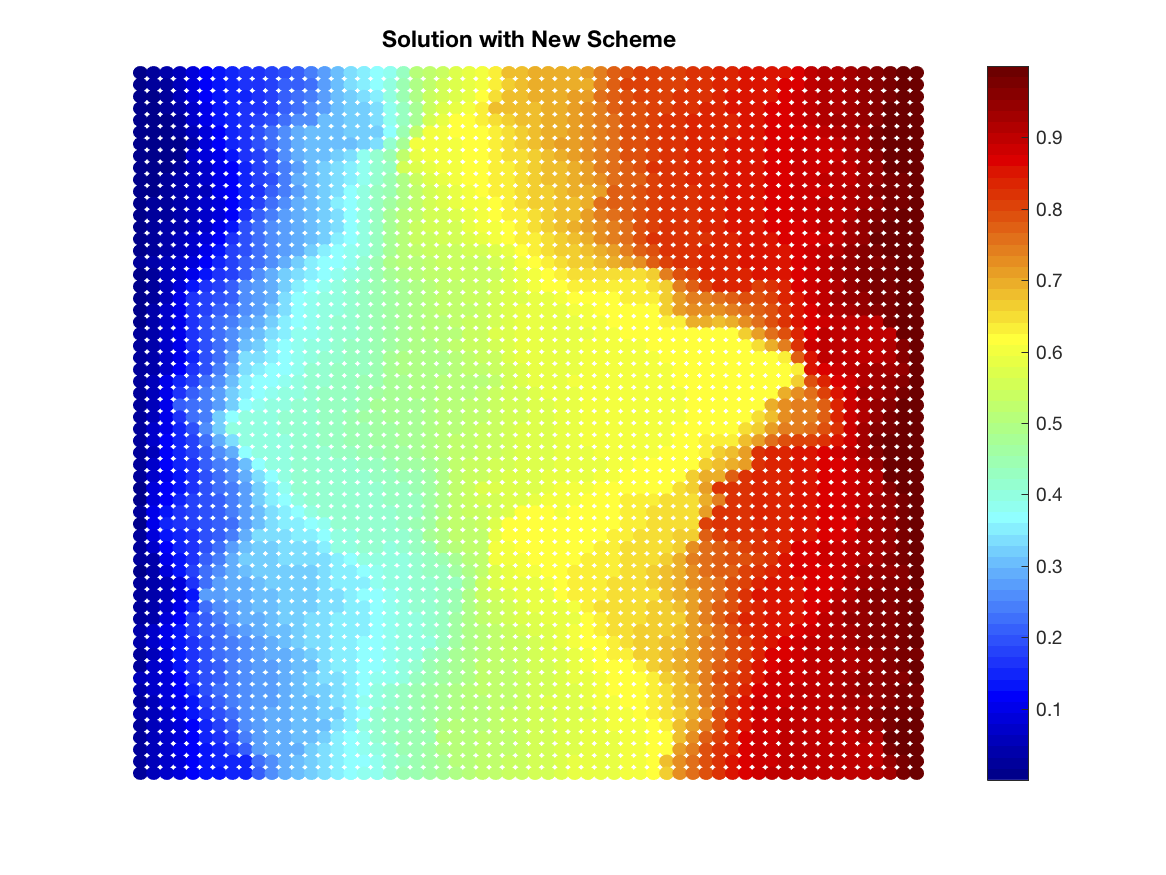}
			\hspace{-3.0cm}
			\caption*{Numerical Solution (2D Case)}
			\label{SPE10LVScheme2DSol}
		\end{minipage}%
		\begin{minipage}{3in}%
			\centering
			\hspace{-4.0cm}
			\includegraphics[width=5cm, height=4cm]{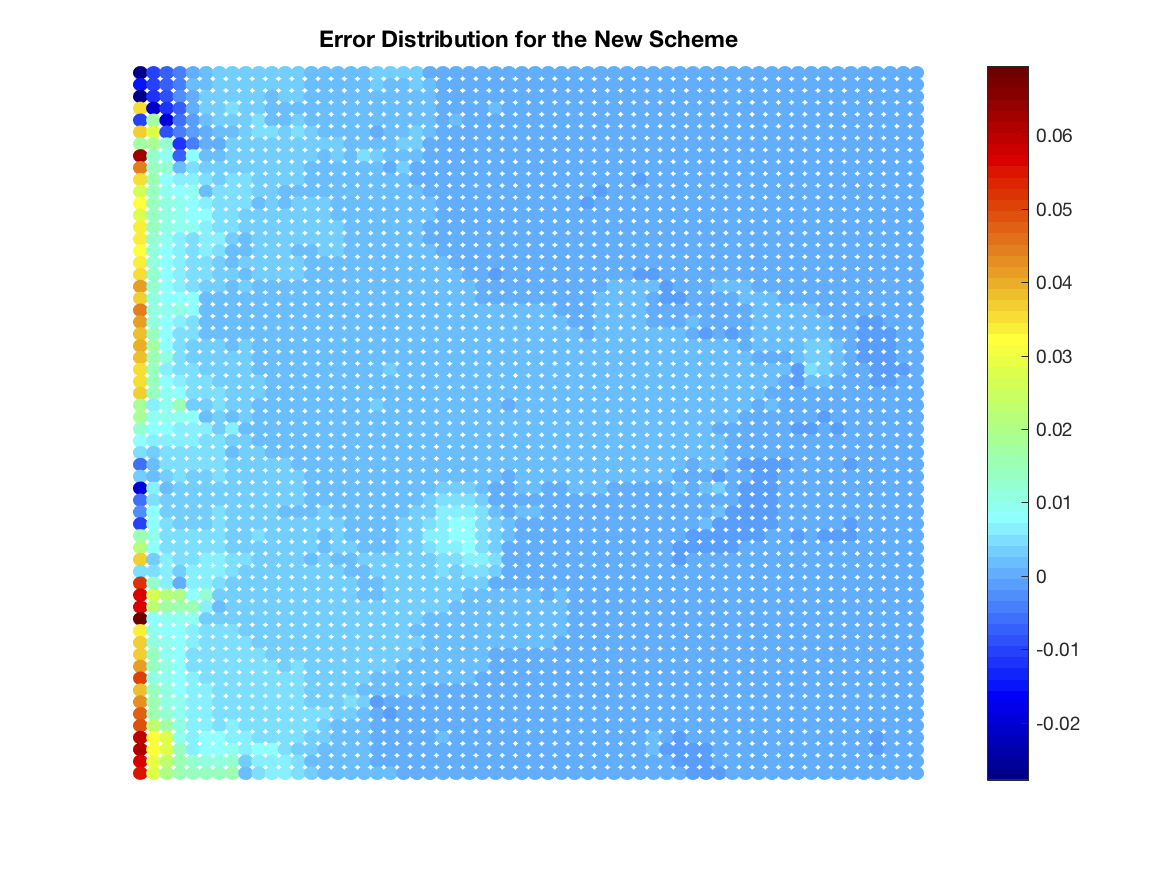}
			\hspace{-3.0cm}
			\caption*{Relative Error \% (2D Case) }
			\label{SPE10LVErrorScheme2DSol}
		\end{minipage}%
	}
	\\
	\mbox{
		\begin{minipage}{3in}%
			\centering
			\hspace{-4.0cm}
			\includegraphics[width=5cm, height=4cm]{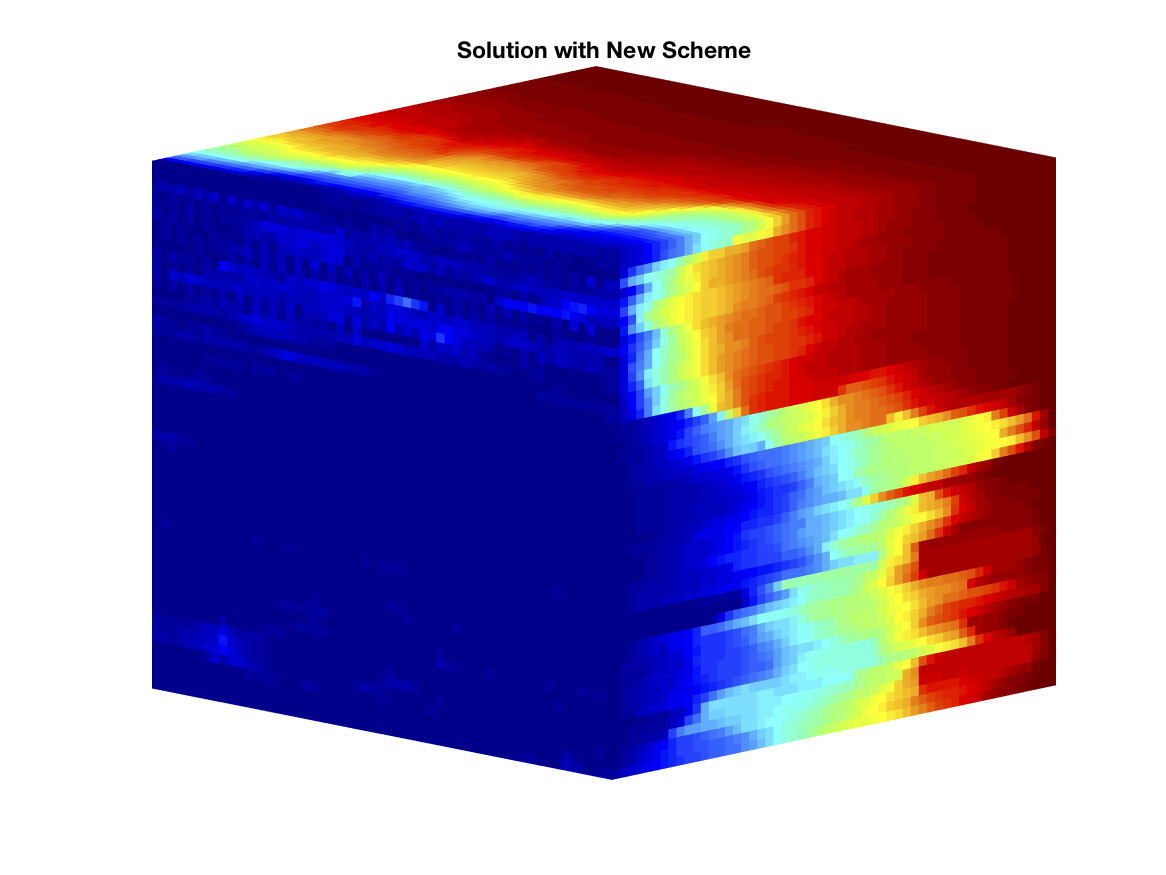}
			\hspace{-3.0cm}
			\caption*{Numerical Solution (3D Case)}
			\label{SPE10LVScheme2DSol3D}
		\end{minipage}%
		\begin{minipage}{3in}%
			\centering
			\hspace{-4.0cm}
			\includegraphics[width=5cm,height=4cm]{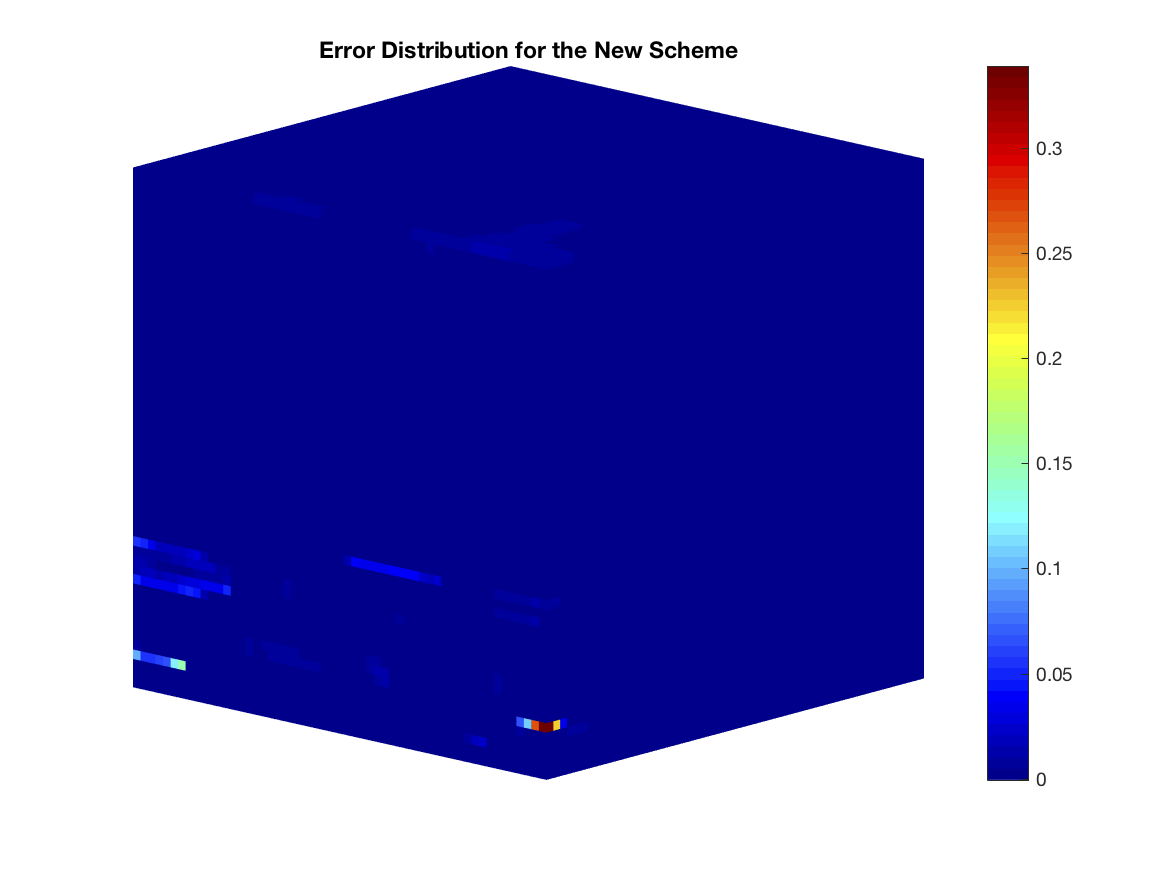}
			\hspace{-3.0cm}
			\caption*{Relative Error \% (3D Case) }
			\label{SPE10LVErrorScheme2DSol3D}
		\end{minipage}%
	}
	\caption{Numerical solution obtained with the new scheme for the SPE10 cases. The relative 
    error distribution is also shown for 2D and 3D cases, where error is computed using 
    \eqref{ChapIV6_1} and numerical solution based on TPFA. }  
	\label{Total_Numerical_Solution}
\end{figure}

\begin{figure}[!htbp]
	\centering
	\includegraphics[width=4.0in,height=2.0in]{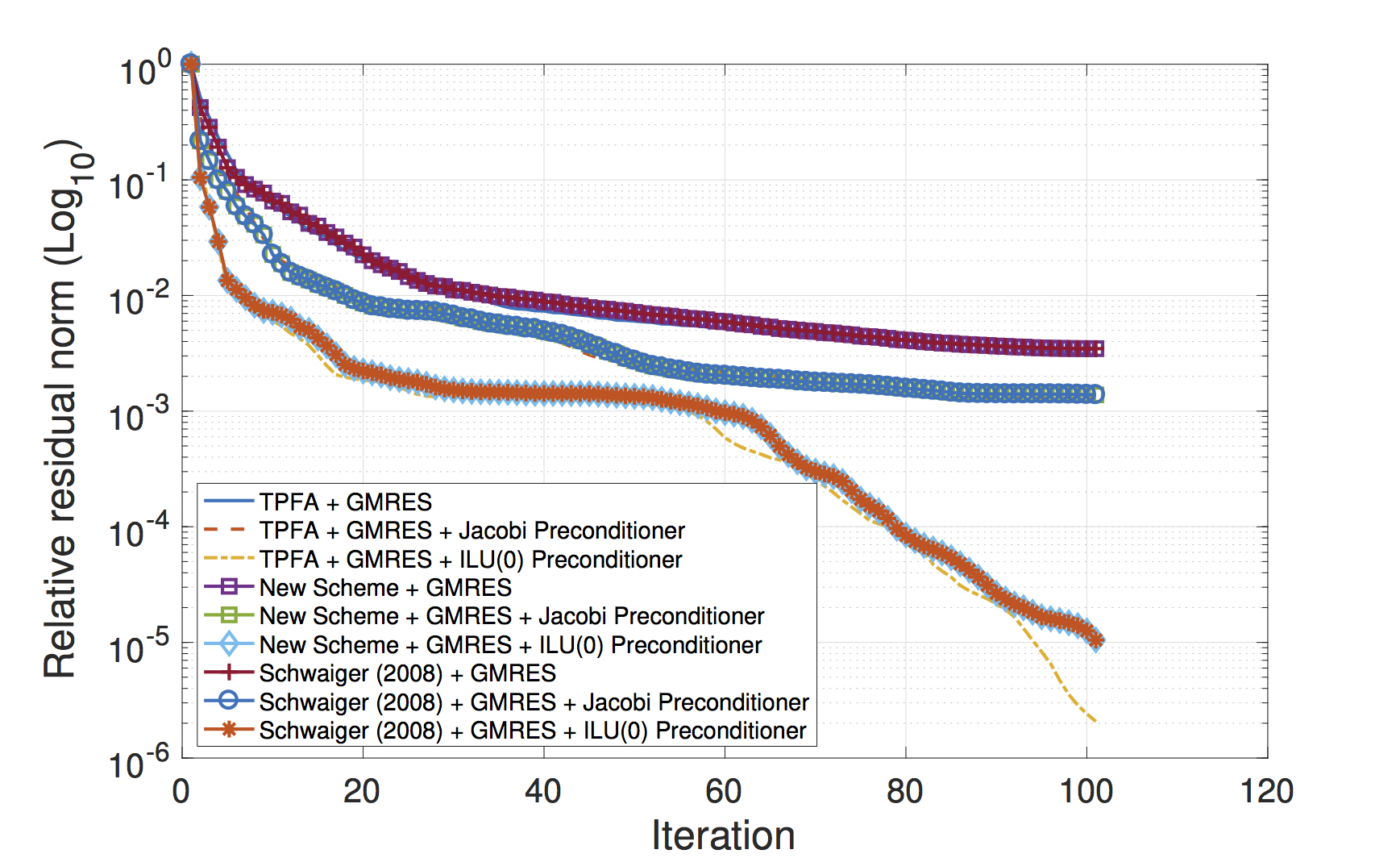}
	\caption{Comparison of the convergence rate of the GMRES iterative method with and without 
	preconditioner for different numerical discretization methods (TPFA, New Scheme 
	\eqref{NewScheme1}--\eqref{NewScheme3}, Schwaiger \eqref{Schwaiger1}--\eqref{Schwaiger5}). The 
	linear system of equations is built using 2D case.}
	\label{IterativeConvergenceRateSPE10Layer85_60by60}
\end{figure}

\begin{figure}[!htbp]
	\centering
	\includegraphics[width=4.0in,height=2.0in]{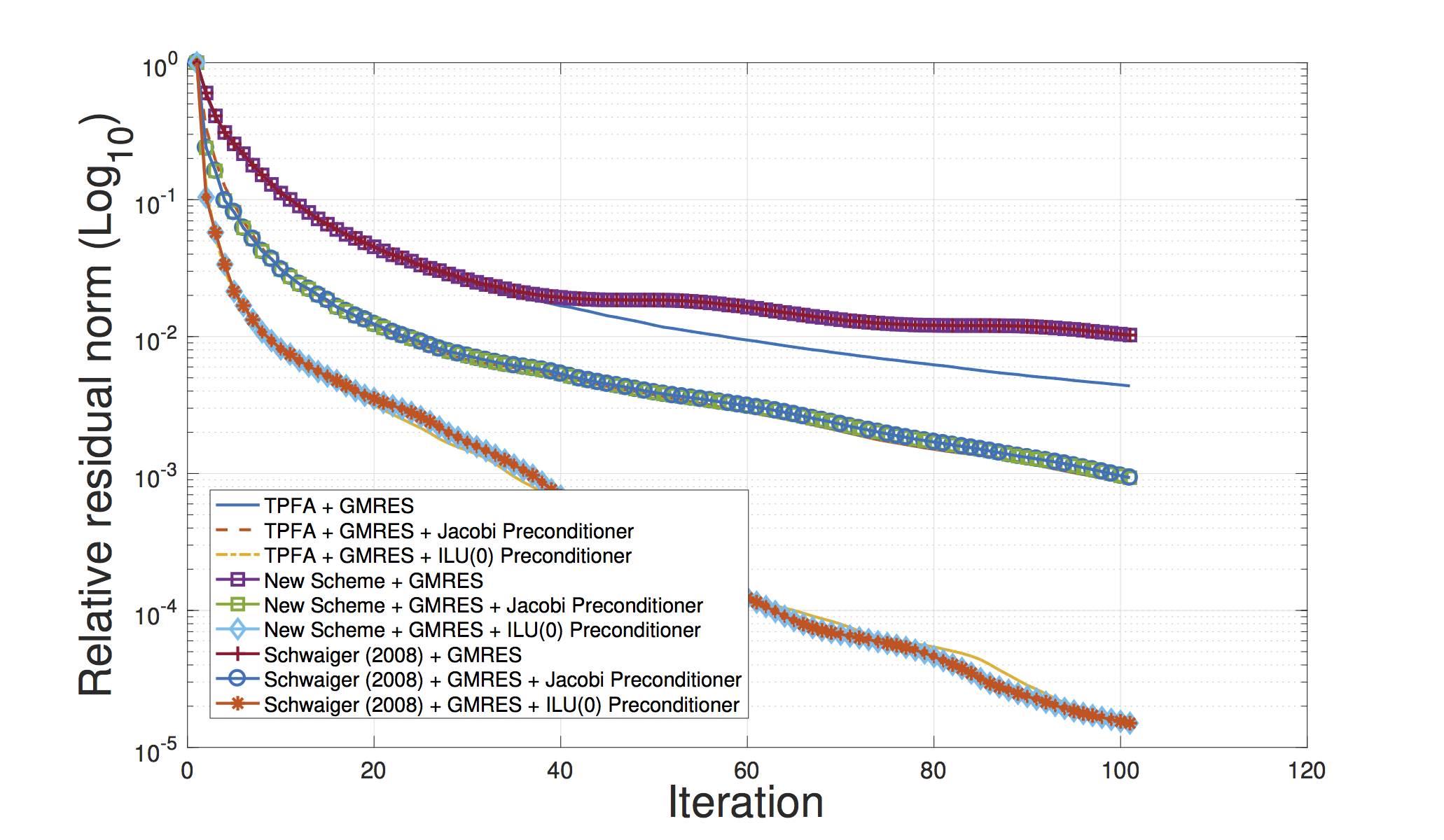}
	\caption{Comparison of the convergence rate of the GMRES iterative method with and 
	without preconditioner for different numerical discretization methods (TPFA, New Scheme \eqref{NewScheme1}--\eqref{NewScheme3}, Schwaiger \eqref{Schwaiger1}--\eqref{Schwaiger5}). 
	The linear system of equations is built using 3D case.}	
    \label{IterativeConvergenceRateSPE103D_60by60}
\end{figure}

In addition, we investigate the accuracy of the numerical 
schemes using a well-known SPE10 benchmark \cite{SPE10} 
with the Laplace equation:
\begin{equation}
\label{LV_SPH_Lap_Eq}
\begin{array}{lcl}
\displaystyle
\vspace{0.2cm}
\mathbf{L}\left(\mathbf{u}\right) = -
\nabla\left(\mathbf{M}\left(\vect{r}\right)\nabla \mathbf{u}
\left(\vect{r}\right)\right) =0, \
\forall \vect{r}\in\Omega\subset\mathbb{R}^{n}, \\
\mathbf{M}^{\alpha\beta}\left(\vect{r}\right) = K^{\alpha\beta}\left(\vect{r}\right), \
K^{\alpha\beta}\left(\vect{r}\right) = 0 \ \forall \alpha\neq\beta;
\ \alpha, \beta = 1,\ldots,n;
\end{array}
\end{equation}
where $\mathbf{u}\left(\vect{r}\right)$ is the unknown pressure field, 
$K^{\alpha\beta}\left(\vect{r}\right)$ is the diagonal permeability field.
The original model contains 85 layers, where layers from 1 to 35 
have smooth permeability with lognormal distribution and  layers 36 to 85 
have channelized formations that are considered to be significantly more 
challenging for numerical simulations. The subset of this model 
is defined by the global Cartesian 
indices $I$, $J$, $K$. The 85 layer was 
used  as a permeability field for 2D simulation with Cartesian indices 
$I=1:60$ and $J=1:60$. In 3D, the subsection Cartesian indices $I=1:60$, 
$J=1:60$, and $K=1:60$. 
The permeability fields $K^{11}$ for both cases are  shown in 
Figure \ref{SPE10Fig}.
The boundary conditions correspond to the unit pressure drop over the entire
domain in $J$-direction (i.e., $y_{min}=0$ and $y_{max}=1$). The numerical results 
using new scheme for the SPE10 cases are presented in Figure \ref{Total_Numerical_Solution}. 
The relative error distribution is also shown for 2D and 3D cases, where 
the error is computed using \eqref{ChapIV6_1} 
and numerical solution based on TPFA.
Figures \ref{IterativeConvergenceRateSPE10Layer85_60by60} and 
\ref{IterativeConvergenceRateSPE103D_60by60}  compare the 
convergence rates of the various discretizations subject to different
preconditioners. 
 Furthermore, if we notice the convergence comparison for this test in Figures 
 \ref{IterativeConvergenceRateSPE10Layer85_60by60} and 
 \ref{IterativeConvergenceRateSPE103D_60by60}, the proposed new method 
 offers similar convergence as compared to the Schwaiger's scheme 
 \eqref{Schwaiger1}--\eqref{Schwaiger5}). The reason for such convergence is the 
 similarity between condition numbers of the linear system of 
 equations resulting from these methods.

Another interesting observation in this numerical test is the small values of the 
relative error between the TPFA solution and proposed new method for both 2D 
and 3D cases. This observation exhibits the efficiency of the proposed meshless 
discretization scheme, which explains the higher accuracy of linear reproduction.

\section{Conclusion}
\label{Summary}

In this paper, the new stable SPH discretization of the elliptic 
operator for heterogeneous media is proposed. The scheme has the two-point 
flux approximation nature and can be written in the form of \eqref{TPFA_FORM}. 
Using this structure, it was possible to make some theoretical monotonicity analysis 
(see, Remarks 2 and 3), which is difficult to perform for other schemes (e.g., Schwaiger's 
method). Furthermore, it follows from the Taylor's series analysis that proposed 
scheme is the optimum one in this class \eqref{TPFA_FORM} for a diagonal matrix 
of the operator coefficients. In addition, the proposed scheme allows to apply 
upwinding strategy during the solution of nonlinear PDEs. 

The new scheme is based on a gradient approximation commonly 
used in thermal, viscous, and pressure projection problems and can be extended to 
include higher-order terms in the appropriate Taylor series. The proposed 
new scheme is combined with mixed corrections which ensure 
linear completeness. The mixed correction utilizes Shepard Functions in combination 
with a correction to derivative approximations. Incompleteness of the kernel support 
combined with the lack of consistency of the kernel interpolation in conventional 
meshless method results in fuzzy boundaries. In the presented 
meshless method, the domain boundary conditions and 
internal field variables are approximated with the default 
accuracy of the method. The resulting new scheme not 
only ensures first order accuracy $\mathcal{O}(h^{\alpha}), \ 1\leq\alpha \leq2$, 
where $h$ denotes the maximum particle spacing, but also minimize the impact 
of the particle deficiency (kernel support incompleteness) problem. 

Furthermore, different Kernel gradients and their impact on the property of 
the scheme and accuracy are discussed. The model was tested by solving an
inhomogeneous Dirichlet and mixed boundary value problems 
for the Laplacian equation with good accuracy confirming our theoretical results.
The accuracy of Schwaiger's scheme and new scheme is the same for 
homogeneous particle distribution and different for the distorted particles. The new 
scheme takes into account all terms related to Hessian of the unknown function 
(see, relation \eqref{ChapIV2_18}). 
 
The stability analysis shows that von Neumann growth factor 
has both real and imaginary parts forming complex shape for general 
particle distribution but satisfying the stability requirement 
almost everywhere. The paper also discusses the monotonicity and convergence 
properties of the new proposed scheme and demonstrates that there is a parameter 
$h$ such that the proposed new scheme is unconditionally monotone with the 
scalar heterogeneous media.         

As was previously mentioned in the introduction, several
methods have been proposed to address the difficulties involved in calculating
second-order derivatives with SPH. In contrast to the present formulation, many
of these methods achieve a high accuracy through fully calculating the Hessian
or requiring that the discrete equations exactly reproduce quadratic functions.
The primary attraction of the present method is that it provides a weak
formulation for Darcy's law which can be of use in further
development of meshless methods.  The SPH model was previously used to
model three-dimensional miscible flow and transport in porous media with complex
geometry, and we are planning to use this model in future work for large (field)
scale simulation of transport in porous media with general permeability distributions.

   \bibliographystyle{unsrt}
   \bibliography{LukyanovVuik_JCP}{}

\end{document}